\documentclass{amsart}
\usepackage{amssymb, latexsym}
\usepackage{hyperref}

\allowdisplaybreaks

\newcommand{\Z}{{\mathbb{Z}}}
\newcommand{\C}{{\mathbb{C}}}
\newcommand{\R}{{\mathbb{R}}}

\let\Re=\undefined\DeclareMathOperator*{\Re}{Re}
\let\Im=\undefined\DeclareMathOperator*{\Im}{Im}

\DeclareMathOperator*{\supp}{supp}

\DeclareMathOperator{\osc}{osc}

\newcommand{\fdt}{{\frac 4{d-2}}}

\newcommand{\tdt}{\frac {2d}{d-2}}
\newcommand{\dtt}{\frac {d-2}2}
\newcommand{\sct}{\frac {2(d+2)}{d-2}}
\newcommand{\ird}{ I\times\R^d}

\newcommand{\eps}{{\varepsilon}}
\newcommand{\propagate}{{e^{it\Delta}}}

\newcommand{\gnj}{g_n^j}
\newcommand{\tnj}{t_n^j}

\newcommand{\wnJ}{w_n^J}
\newcommand{\twnJ}{\tilde w_n^J}
\newcommand{\vnj}{v_n^j}
\newcommand{\unJ}{u_n^J}

\theoremstyle{plain}
\newtheorem{theorem}{Theorem}
\newtheorem{proposition}[theorem]{Proposition}
\newtheorem{lemma}[theorem]{Lemma}
\newtheorem{corollary}[theorem]{Corollary}
\newtheorem{conjecture}[theorem]{Conjecture}
\theoremstyle{definition}
\newtheorem{definition}[theorem]{Definition}
\newtheorem{remark}[theorem]{Remark}

\newcounter{smalllist}

\newenvironment{CI}{\begin{list}{{\ $\bullet$\ }}{%
\setlength{\topsep}{0mm}\setlength{\parsep}{0mm}\setlength{\itemsep}{0mm}%
\setlength{\labelwidth}{0mm}\setlength{\itemindent}{0mm}\setlength{\leftmargin}{0mm}%
\setlength{\labelsep}{0mm}}}{\end{list}}

\numberwithin{equation}{section}
\numberwithin{theorem}{section}

\begin{document}

\title[Focusing energy-critical NLS]{The focusing energy-critical nonlinear Schr\"odinger equation in dimensions\\five and higher}
\author{Rowan Killip}
\address{University of California, Los Angeles, and Institute for Advanced Study, Princeton, NJ}
\author{Monica Visan}
\address{Institute for Advanced Study, Princeton, NJ}

\begin{abstract}
We consider the focusing energy-critical nonlinear Schr\"odinger equation $iu_t+\Delta u = - |u|^{\frac4{d-2}}u$ in dimensions
$d\geq 5$. We prove that if a maximal-lifespan solution $u:I\times\R^d\to \C$
obeys $\sup_{t\in I}\|\nabla u(t)\|_2<\|\nabla W\|_2$, then it is global and scatters both forward and backward in time. Here
$W$ denotes the ground state, which is a stationary solution of the equation.  In particular, if a solution has both energy and
kinetic energy less than those of the ground state $W$ at some point in time, then the solution is global and scatters.  We also
show that any solution that blows up with bounded kinetic energy must concentrate at least the kinetic energy of the ground
state.  Similar results were obtained by Kenig and Merle in \cite{Evian, kenig-merle} for spherically symmetric initial data and
dimensions $d=3,4,5$.
\end{abstract}

\maketitle

\tableofcontents

\section{Introduction}
We consider the initial-value problem for the focusing energy-critical nonlinear Schr\"odinger equation in dimension $d\geq 3$,
\begin{equation}\label{nls}
\begin{cases}
\ iu_t+\Delta u=F(u)\\
\ u(t_0)=u_0\in \dot H^1_x(\R^d),
\end{cases}
\end{equation}
where the nonlinearity is given by $F(u)=- |u|^{\fdt}u$.  As indicated in the title, our main results are for dimensions $d\geq 5$;
nevertheless, many of our arguments remain valid in dimensions three and four.

The name `energy-critical' refers to the fact that the scaling symmetry
\begin{equation}\label{scaling}
u(t,x) \mapsto  u_\lambda(t,x):= \lambda^{\frac{d-2}2} u( \lambda^2 t, \lambda x)
\end{equation}
leaves both the equation and the energy invariant.  The energy of a solution is defined by
\begin{equation*}
E(u(t)) := \int_{\R^d} \bigl(\tfrac 12 |\nabla u(t,x)|^2 - \tfrac{d-2}{2d} |u(t,x)|^{\frac{2d}{d-2}}\bigr)\, dx
\end{equation*}
and is conserved under the flow; see Theorem~\ref{T:local} below.  We refer to the gradient term in the formula above as the
\emph{kinetic energy} and to the second term as the \emph{potential energy}.  Note that the potential energy is negative,
which expresses the focusing nature of the nonlinearity.

\begin{definition}[Solution] A function $u: I \times \R^d \to \C$ on a non-empty time interval $t_0\in I \subset \R$
is a \emph{solution} (more precisely, a strong $\dot H^1_x(\R^d)$ solution) to \eqref{nls} if it lies in the class $C^0_t \dot
H^1_x(K \times \R^d) \cap L^{2(d+2)/(d-2)}_{t,x}(K \times \R^d)$ for all compact $K \subset I$, and obeys the Duhamel
formula
\begin{align}\label{old duhamel}
u(t_1) = e^{i(t_1-t_0)\Delta} u(t_0) - i \int_{t_0}^{t_1} e^{i(t_1-t)\Delta} F(u(t))\ dt
\end{align}
for all $t_1 \in I$.  We refer to the interval $I$ as the \emph{lifespan} of $u$. We say that $u$ is a \emph{maximal-lifespan
solution} if the solution cannot be extended to any strictly larger interval. We say that $u$ is a \emph{global solution} if $I
= \R$.
\end{definition}

The condition that $u$ is in $L^{2(d+2)/(d-2)}_{t,x}$ locally in time is natural for several reasons (see
Theorem~\ref{T:local}): (i) By Strichartz inequality, all solutions to the linear problem lie in this space. (ii) Local
solutions exist in this space. (iii) Finiteness of this norm on the maximal interval of existence implies that the solution is
global and scatters both forward and backward in time.  (iv)  A posteriori, one can show that the locally
$L^{2(d+2)/(d-2)}_{t,x}$ solution is in fact the only solution belonging to $C^0_t \dot H^1_x(I \times \R^d)$.

In view of point (iii) above, it is natural to define the \emph{scattering size} of a solution to \eqref{nls}
on a time interval $I$ by
$$
S_I(u):= \int_I \int_{\R^d} |u(t,x)|^{\frac{2(d+2)}{d-2}}\, dx \,dt.
$$

Associated to the notion of solution is a corresponding notion of blowup.  As we will see in Theorem~\ref{T:local}, this
precisely corresponds to the impossibility of continuing the solution.

\begin{definition}[Blowup]\label{D:blowup}
We say that a solution $u$ to \eqref{nls} \emph{blows up forward in time} if there exists a time $t_1 \in I$ such that
$$ S_{[t_1, \sup I)}(u) = \infty$$
and that $u$ \emph{blows up backward in time} if there exists a time $t_1 \in I$ such that
$$ S_{(\inf I, t_1]}(u) = \infty.$$
\end{definition}

The local theory for \eqref{nls} was worked out by Cazenave and Weissler \cite{cwI}.  They performed a fixed-point argument to
construct local-in-time solutions for arbitrary initial data in $\dot H^1_x(\R^d)$; however, as is the case with critical
equations, the resulting time of existence depends on the profile of the initial data and not merely on its $\dot H^1_x$-norm.
They also constructed global solutions from small initial data.  Unconditional uniqueness was subsequently proved
by Cazenave, \cite[Proposition 4.2.5]{cazenave:book}.  We summarize their results in the theorem below.

\begin{theorem}[Local well-posedness, \cite{cwI, cazenave:book}]\label{T:local}
Given $u_0\in\dot H_x^1(\R^d)$ and $t_0\in \R$, there exists a unique maximal-lifespan solution $u:\ird \rightarrow \C$ to
\eqref{nls} with initial data $u(t_0)=u_0$.  This solution also has the following properties:

\begin{CI}
\item (Local existence) $I$ is an open neighbourhood of $t_0$.
\item (Energy conservation) The energy of $u$ is conserved, that is, $E(u(t))=E(u_0)$ for all $t\in I$.
\item (Blowup criterion) If $\sup(I)$ is finite, then $u$ blows up forward in time; if $\inf(I)$ is finite,
then $u$ blows up backward in time.
\item (Scattering) If $\sup(I)=+\infty$ and $u$ does not blow up forward in time, then $u$ scatters forward in time, that is,
there exists a unique $u_+ \in \dot H^1_x(\R^d)$ such that
\begin{equation}\label{like u+}
\lim_{t \to +\infty} \| u(t)-e^{it\Delta} u_+ \|_{\dot H^1_x(\R^d)} = 0.
\end{equation}
Conversely, given $u_+ \in \dot H^1_x(\R^d)$ there is a unique solution to \eqref{nls} in a neighbourhood of infinity
so that \eqref{like u+} holds.
\item (Small data global existence) If $\|\nabla u_0\|_2$ is sufficiently small (depending on $d$), then $u$ is a global solution
which does not blow up either forward or backward in time.  Indeed, in this case $S_\R(u)\lesssim \|\nabla u_0\|_2^{2(d+2)/(d-2)}$.
\item (Unconditional uniqueness) If $\tilde u \in C^0_t \dot H^1_x(J \times \R^d)$ with $t_0\in J$, obeys \eqref{old duhamel}
and $\tilde u(t_0)=u_0$, then $J\subseteq I$ and $\tilde u \equiv u$ throughout $J$.
\end{CI}
\end{theorem}

A variant of the local well-posedness theorem above is the following lemma.

\begin{lemma}[Stability, \cite{TV}]\label{L:stab} Let $d\geq 3$.  For every $E,L > 0$ and $\eps > 0$ there exists $\delta > 0$
with the following property:  Suppose $\tilde u: I \times \R^d \to \C$ is an approximate solution to \eqref{nls} in the sense
that
\begin{equation*}
\bigl\| \nabla \bigl[i\tilde u_t + \Delta \tilde u - F(\tilde u)\bigr] \bigr\|_{L_{t,x}^{\frac{2(d+2)}{d+4}}(I\times\R^d)} \leq \delta
\end{equation*}
and also obeys
\begin{align*}
\|\tilde u\|_{L_t^\infty \dot H^1_x(\ird)}\leq E \quad \text{and} \quad S_I(\tilde u)\leq L.
\end{align*}
If $t_0 \in I$ and $u_0 \in \dot H^1_x(\R^d)$ are such that
$$
\|\tilde u(t_0)-u_0\|_{\dot H^1_x(\R^d)} \leq \delta,
$$
then there exists a solution $u: I \times \R^d \to \C$ to \eqref{nls} with $u(t_0) = u_0$ such that
\begin{align*}
\|\tilde u -u\|_{L_t^\infty \dot H^1_x(\ird)} + S_I(\tilde u - u) \leq \eps.
\end{align*}
\end{lemma}

\begin{remark}
The result in \cite{TV} is slightly more general; we merely stated the version we will use. Lemma~\ref{L:stab} implies the
existence and uniqueness of maximal-lifespan solutions in Theorem~\ref{T:local}.  It also proves that the solution depends
uniformly continuously on the initial data (on bounded sets), which was missing from \cite{cwI,cazenave:book}. See also
\cite{ckstt:gwp, RV} for earlier results in dimensions three and four.
\end{remark}

The defocusing case, that is, $F(u)=|u|^{\frac{4}{d-2}}u$, has received a lot of attention.
It is known that all $\dot H^1_x$ initial data lead to global solutions with finite scattering size.
Indeed, this was proved by Bourgain \cite{borg:scatter}, Grillakis \cite{grillakis}, and
Tao \cite{tao: gwp radial} for spherically symmetric initial data, and by Colliander--Keel--Staffilani--Takaoka--Tao
\cite{ckstt:gwp}, Ryckman--Visan \cite{RV}, and Visan \cite{thesis:art, Monica:thesis} for arbitrary initial data.

In the focusing case, things are more subtle.  From Theorem~\ref{T:local}, we see that all maximal-lifespan solutions with
sufficiently small kinetic energy are global and scatter.  However,
$$
W(t,x) = W(x):=\frac 1{(1+\frac{|x|^2}{d(d-2)})^{\dtt}}\in \dot H^1_x(\R^d)
$$
is a stationary solution to \eqref{nls}, that is, $W$ solves the nonlinear elliptic equation
\begin{align*}
\Delta W + |W|^{\frac4{d-2}}W=0.
\end{align*}
(See Appendix~\ref{A: W} for further properties of the ground state $W$, including its connection to the sharp version of
Sobolev embedding.) In particular, $W$ is a solution to \eqref{nls} that blows up both forward and backward in time in the sense
of Definition~\ref{D:blowup}. It is believed that $W$ has minimal kinetic energy among all blowup solutions.  More precisely, we
have

\begin{conjecture}\label{conj}
Let $d\geq 3$ and let $u:I\times\R^d\to\C$ be a solution to \eqref{nls}.  If
$$
E_* :=  \sup_{t\in I} \|\nabla u(t)\|_2 < \|\nabla W\|_2,
$$
then
$$
\int_I \int_{\R^d} |u(t,x)|^{\frac{2(d+2)}{d-2}}\, dx\, dt \leq C(E_*) < \infty.
$$
\end{conjecture}

An analogous conjecture is believed to hold in the mass-critical case.  There, the role of $W$ is played by the ground state
$Q$, which is also a maximizer in the sharp Gagliardo--Nirenberg inequality.  Indeed, Weinstein \cite{weinstein} first realized
the importance of $Q$ as a minimal-mass blowup example, albeit in the finite energy case.  The mass-critical conjecture for
$L_x^2$ initial data and dimensions $d\geq 2$ was recently settled in the spherically symmetric case; see \cite{KTV} and
\cite{KVZ}.

Conjecture~\ref{conj} was verified in dimensions $d=3,4,5$ for spherically symmetric initial data by Kenig and Merle \cite{Evian, kenig-merle}.
In this paper, we verify the conjecture in dimensions $d\geq 5$ without any further symmetry assumptions on the solution.
More precisely, we derive

\begin{theorem}[Spacetime bounds]\label{T:main}
Let $d\geq 5$ and let $u:I\times\R^d\to\C$ be a solution to \eqref{nls}.  If
$$
E_* := \sup_{t\in I} \|\nabla u(t)\|_2 < \|\nabla W\|_2,
$$
then
$$
\int_I \int_{\R^d} |u(t,x)|^{\frac{2(d+2)}{d-2}}\, dx\, dt \leq C(E_*) < \infty.
$$
\end{theorem}

The key development that allows us to treat non-radial data is a proof that certain minimal kinetic energy blowup solutions
have finite mass; indeed, that they belong to $L_t^\infty \dot H^{-\eps}_x$ for some $\eps>0$.  This is done in Section~\ref{S:neg}.
Examination of the stationary solution $u(t,x)=W(x)$ shows that this is not true in dimensions three and four, indicating
a difficulty intrinsic to these dimensions.  At a more technical level, the main novelty presented here is a proof that
such minimal kinetic energy solutions exhibit additional decay in $L^p_x$-sense.  This is then harnessed in a double Duhamel trick
(cf. \cite{tao:attractor}) to prove the $L_x^2$-based properties mentioned above.

Combining Theorem~\ref{T:main} with the local theory gives the following more appealing formulation.
However, the slightly stronger statement in Theorem~\ref{T:main} is helpful for applications;
an example of this is Theorem~\ref{T:conc} below.

\begin{corollary}[Scattering]\label{C:main}
Let $d\geq 5$ and let $u$ be a solution to \eqref{nls} with maximal lifespan $I$.  Assume also that
$$
\sup_{t\in I} \|\nabla u(t)\|_2 < \|\nabla W\|_2.
$$
Then $I=\R$ and
$$
\int_\R \int_{\R^d} |u(t,x)|^{\frac{2(d+2)}{d-2}}\, dx\, dt < \infty.
$$
\end{corollary}

A more effective criterion for global well-posedness (depending directly on $u_0$) can be obtained from Corollary~\ref{C:main}
using an energy-trapping argument of Kenig and Merle \cite{kenig-merle}; see Corollary~\ref{trap}.

\begin{corollary}\label{C:inductE}
Let $d\geq 5$ and let $u_0\in \dot H^1_x(\R^d)$ be such that $\|\nabla u_0\|_2\leq \|\nabla W\|_2$ and $E(u_0)<E(W)$.
Then the corresponding solution $u$ to \eqref{nls} is global and moreover,
$$
\int_\R \int_{\R^d} |u(t,x)|^{\frac{2(d+2)}{d-2}}\, dx\, dt <\infty.
$$
\end{corollary}

We would like to note that the proof of Theorem~\ref{T:main} adapts without difficulty (indeed, with some simplifications)
to the defocusing case $F(u)=|u|^{\frac4{d-2}}u$; as such, it constitutes a new (more streamlined) derivation of the main results
in \cite{thesis:art, Monica:thesis}.

The result in Theorem~\ref{T:main} is sharp.  Indeed, the ground state $W$ is a solution to \eqref{nls} that blows up at infinity
(in the sense of Definition~\ref{D:blowup}) in both time directions.  Moreover, there exist solutions with kinetic energy
only slightly greater than that of $W$ which blow up in finite time.  More precisely, in Section~\ref{S:blowup} we prove

\begin{proposition}[Blowup]\label{P:blowup}
Let $d\geq 3$ and  $u_0\in \dot H_x^1(\R^d)$ with $E(u_0)<E(W)$ and $\|\nabla u_0\|_2 \geq \|\nabla W \|_2$.  Assume also that
either $xu_0\in L^2_x(\R^d)$ or $u_0\in H_x^1(\R^d)$ is radial.  Then the corresponding solution $u$ to \eqref{nls} blows up in
finite time.
\end{proposition}

Results of this type were previously obtained in the energy-subcritical case.  Ogawa and Tsutsumi \cite{OgawaTsutsumi} treated
the case of initial data with negative energy.  Holmer and Roudenko \cite{HolmerRoudenko} extended their result to include
certain positive-energy initial data.  Proposition~\ref{P:blowup} appears in \cite{kenig-merle}, where a proof is given in the case
$xu_0\in L^2_x(\R^d)$.  For a complete proof, see Section~\ref{S:blowup}.

In view of Corollary~\ref{C:inductE} and Proposition~\ref{P:blowup}, one may inquire about the case $E(u)=E(W)$.  For radial data
in dimensions $3$, $4$, and $5$, this is discussed in \cite{DuckyMerle}.

Our last result proves that blowup solutions with bounded kinetic energy must concentrate a fixed amount of
kinetic energy around the blowup time.  The argument given shows that this result follows from a positive answer
to Conjecture~\ref{conj}; the dimensional restriction below reflects only the fact that this conjecture
is currently open in three and four dimensions.

\begin{theorem}[Blowup solutions concentrate kinetic energy]\label{T:conc}
Fix $d\geq 5$ and let $u$ be a solution to \eqref{nls} that blows up at time $T^*\in[-\infty, \infty]$.  Assume also that
\begin{align}\label{type II}
\limsup_{t\to T^*}\|\nabla u(t)\|_2<\infty.
\end{align}
If $T^* $ is finite, then there exists a sequence $t_n\to T^*$ such that for any sequence $R_n\in (0,\infty)$ obeying
$|T^*-t_n|^{-\frac 12}R_n\to\infty$,
$$
\limsup_{n\to \infty} \sup_{x_0\in \R^d}\int_{|x-x_0|\leq R_n} |\nabla u(t_n,x)|^2\, dx \geq \|\nabla W\|_2^2.
$$
If $|T^*|=\infty$, then there exists a sequence $t_n\to T^*$ such that for any sequence $R_n\in (0,\infty)$ obeying
$|t_n|^{-\frac 12}R_n\to\infty$,
$$
\limsup_{n\to \infty} \sup_{x_0\in \R^d}\int_{|x-x_0|\leq R_n} |\nabla u(t_n,x)|^2\, dx \geq \|\nabla W\|_2^2.
$$
\end{theorem}

We prove Theorem~\ref{T:conc} in Section~\ref{S:Conc}.  The argument is inspired by that in \cite{KTV}, which employs some ideas from
\cite{bourg.2d}. The fact that the kinetic energy is not conserved introduces new subtleties.
A similar result in dimensions $d=3,4,5$ for spherically symmetric data was proposed by Kenig and Merle \cite{kenig-merle};
see Section~\ref{S:Conc} for a fuller discussion.

\subsection{Outline of the proof of Theorem~\ref{T:main}}

We argue by contradiction.  We show that if the theorem failed, it would imply the existence of a very special
type of counterexample.  Such counterexamples are then shown to have a wealth of properties not immediately
apparent from their construction, so many properties, in fact, that they cannot exist.

While we will make some further reductions later, the main property of the special counterexamples is
almost periodicity modulo symmetries:

\begin{definition}[Almost periodicity modulo symmetries]\label{D:ap}
Let $d \geq 3$.  A solution $u$ to \eqref{nls} with lifespan $I$ is said to be \emph{almost periodic modulo symmetries} if there
exist functions $N: I \to \R^+$, $x:I\to \R^d$, and $C: \R^+ \to \R^+$ such that for all
$t \in I$ and $\eta > 0$,
$$ \int_{|x-x(t)| \geq C(\eta)/N(t)} |\nabla u(t,x)|^2\, dx \leq \eta$$
and
$$ \int_{|\xi| \geq C(\eta) N(t)} |\xi|^2\, | \hat u(t,\xi)|^2\, d\xi \leq \eta.$$
We refer to the function $N$ as the \emph{frequency scale function} for the solution $u$, $x$ the \emph{spatial center function},
and to $C$ as the \emph{compactness modulus function}.
\end{definition}

\begin{remark}
The parameter $N(t)$ measures the frequency scale of the solution at time $t$, while $1/N(t)$ measures the spatial scale; see
\cite{tvz:cc} for further discussion.  It is possible to multiply $N(t)$ by any function of $t$ that is bounded both above
and below, provided that we also modify the compactness modulus function $C$ accordingly.
\end{remark}

\begin{remark}\label{R:pot energy}
By the Ascoli--Arzela Theorem, a family of functions is precompact in $\dot H^1_x(\R^d)$ if and only if it is norm-bounded and
there exists a compactness modulus function $C$ so that
$$
\int_{|x| \geq C(\eta)} |\nabla f(x)|^2\ dx + \int_{|\xi| \geq C(\eta)} |\xi|^2 \, |\hat f(\xi)|^2\ d\xi \leq \eta
$$
for all functions $f$ in the family.  Thus, an equivalent formulation of Definition~\ref{D:ap} is as follows: $u$ is almost
periodic modulo symmetries if and only if
$$
\{ u(t): t \in I \} \subseteq \{ \lambda^{\frac{d-2}2} f(\lambda (x+x_0)) : \, \lambda\in(0,\infty), \ x_0\in \R^d, \text{ and }f \in K \}
$$
for some compact subset $K$ of $\dot H^1_x(\R^d)$.  In particular, as every compact set in $\dot H^1_x(\R^d)$ is compact in
$L_x^{2d/(d-2)}(\R^d)$ (by Sobolev embedding), any solution $u:\ird\to \C$ to \eqref{nls} that is almost periodic modulo
symmetries must also satisfy
$$
 \int_{|x-x(t)| \geq C(\eta)/N(t)} |u(t,x)|^{\frac{2d}{d-2}}\, dx \leq \eta
$$
for all $t\in I$ and $\eta>0$.
\end{remark}

\begin{remark}\label{R:c small}
A further consequence of compactness modulo symmetries is the existence of a function $c: \R^+\to \R^+$ so that
$$
\int_{|x-x(t)| \leq c(\eta)/N(t)} |\nabla u(t,x)|^2\, dx
    + \int_{|\xi| \leq c(\eta) N(t)} |\xi|^2\, | \hat u(t,\xi)|^2\, d\xi \leq \eta
$$
for all $t \in I$ and $\eta > 0$.
\end{remark}

With these preliminaries out of the way, we can now describe the first major milestone in the proof of Theorem~\ref{T:main}.

\begin{theorem}[Reduction to almost periodic solutions]\label{T:reduct}
Suppose $d \geq 3$ is such that Conjecture~\ref{conj} failed.  Then there exists a maximal-lifespan solution
$u:\ird\to \C$ to \eqref{nls} such that $\sup_{t\in I} \|\nabla u(t)\|_2 < \|\nabla W\|_2$, $u$ is almost periodic modulo
symmetries, and $u$ blows up both forward and backward in time.  Moreover, $u$ has minimal kinetic energy among all blowup solutions,
that is,
$$
\sup_{t\in J} \|\nabla v(t)\|_2 \leq \sup_{t\in I} \|\nabla u(t)\|_2
$$
for all maximal-lifespan solutions $v:J\times\R^d \to \C$ that blowup in at least one time direction.
\end{theorem}

Most of the properties of $u$ described in this theorem stem directly from the fact that it is a minimal kinetic energy blowup
solution.  The innovative discovery that such minimal blowup solutions exist was made by Keraani \cite[Theorem 1.3]{keraani-l2}
in the context of the mass-critical NLS.  This was adapted to the energy-critical setting in \cite{kenig-merle}, which also
constitutes the first application of the existence of minimal blowup solutions to the well-posedness problem.

Following Bourgain~\cite{borg:scatter}, earlier works on the energy-critical NLS (e.g., \cite{ckstt:gwp,RV,thesis:art})
focussed their attention on \emph{almost}-minimal blowup solutions, which are then shown to have space and frequency localization
properties similar to (but slightly weaker than) those in Definition~\ref{D:ap}.  While these earlier methods are inherently
quantitative, they add significantly to the complexity of the argument.

One of the main ingredients in the proof of Theorem~\ref{T:reduct} is a linear profile decomposition of Keraani, which
is reproduced below as Lemma~\ref{L:cc}.  Using this result, Kenig and Merle proved a slight variant of Theorem~\ref{T:reduct}
in dimensions $3$, $4$, and~$5$; see \cite[Proposition~4.1]{kenig-merle}.  They further indicate that the argument may be
modified to give a proof of the full Theorem~\ref{T:reduct} in these dimensions (see the proof of Corollary~5.16 in
\cite{kenig-merle}).  

We will present a complete proof of Theorem~\ref{T:reduct} uniformly in $d\geq 3$.  In doing so, we uncovered
several difficulties in extracting a blowup solution with minimal kinetic energy that were not elaborated upon in \cite{kenig-merle}.
These difficulties are related to the fact that unlike the energy, the kinetic energy is not conserved.
Firstly, choosing a `bad profile' requires a certain amount of gymnastics; see, for instance, Lemma~\ref{L:bad profile}
and the discussion that follows it.  The difficulty arises from the possibility that the scattering size of several profiles
is large over short times, while their kinetic energy does not achieve the critical value until much later.
Secondly, having selected the `bad profile', one must then prove kinetic energy decoupling at the (potentially) \emph{later}
time when this profile achieves the critical kinetic energy; see Lemma~\ref{L:decouple ke}.
Related arguments (for the cubic NLS in three spatial dimensions) appear in \cite{kenig-merle:1/2}.

To prove Theorem~\ref{T:main}, we will need to demonstrate the existence of minimal kinetic energy blowup solutions
with more refined properties than those provided by Theorem~\ref{T:reduct}; in particular, we need to better constrain
the behaviour of the frequency scale function $N(t)$.  Theorem~1.16 in \cite{KTV} is the strongest result of this type
of which we are aware.  In Section~\ref{S:enemies}, we adapt the argument given there to obtain

\begin{theorem}[Three special scenarios for blowup]\label{T:enemies}
Fix $d\geq 3$ and suppose that Conjecture~\ref{conj} fails for this choice of $d$.
Then there exists a minimal kinetic energy, maximal-lifespan solution $u:I\times\R^d\to \C$,
which is almost periodic modulo symmetries, $S_I(u)=\infty$,
and obeys $\sup_{t\in I} \|\nabla u(t)\|_2 < \|\nabla W\|_2$.

We can also ensure that the lifespan $I$ and the frequency scale function $N:I\to\R^+$ match one of the following three scenarios:
\begin{itemize}
\item[I.] (Finite-time blowup) We have that either $|\inf I|<\infty$ or $\sup I<\infty$.
\item[II.] (Soliton-like solution) We have $I = \R$ and
\begin{equation*}
 N(t) = 1 \quad \text{for all} \quad t \in \R.
\end{equation*}
\item[III.] (Low-to-high frequency cascade) We have $I = \R$,
\begin{equation*}
\inf_{t \in \R} N(t) \geq 1, \quad \text{and} \quad \limsup_{t \to +\infty} N(t) = \infty.
\end{equation*}
\end{itemize}
\end{theorem}

Therefore, in order to prove Theorem~\ref{T:main} it suffices to preclude the existence of solutions that satisfy the criteria in
Theorem~\ref{T:enemies}.  The key step in all three scenarios above is to prove negative regularity, that is, the solution $u$
lies in $L_x^2$ or better.  In scenarios~II and~III, the proof that $u\in L^2_x$ requires $d\geq5$; note that in lower dimensions,
the ground-state solution $W$ does not belong to $L^2_x$.  Similar in spirit to \cite{KTV, KVZ}, negative regularity is deduced
from the minimality of the solution considered; recall that $u$ has minimal kinetic energy among all blowup solutions.  In this regard,
our approach differs from those in \cite{borg:scatter, ckstt:gwp, grillakis, RV, tao: gwp radial, thesis:art, Monica:thesis}
where various versions of the Morawetz inequality are used to prove negative regularity.

The fact that the solutions described in Theorem~\ref{T:enemies} belong to $L^2_x$ is a very peculiar property.  General $\dot H^1_x$
initial data do not decay this quickly at infinity.  That it holds for the solutions in Theorem~\ref{T:enemies} is closely tied
to the almost periodicity of these solutions, which is itself a very non-generic behaviour for a dispersive equation.  Indeed,
generic initial data lead to solutions containing waves that radiate to infinity as time progresses, while almost periodic solutions
remain concentrated in space.

The key property that leads to the solutions in Theorem~\ref{T:enemies} having finite mass (and being almost periodic) is their
selection as minimal kinetic energy blowup solutions.  Any low-frequency component that does not contribute directly to the blowup
behaviour would constitute a waste of kinetic energy.  A further manifestation of this minimality is the absence of a scattered wave
at the endpoints of the lifespan $I$; more formally, we have the following Duhamel formula, which plays an important role in
proving negative regularity.  For a proof, see \cite[Section~6]{tvz:cc}.

\begin{lemma}\label{L:duhamel}
Let $u$ be an almost periodic solution to \eqref{nls} on its maximal-lifespan $I$.  Then, for all $t\in I$,
\begin{equation}\label{Duhamel}
\begin{aligned}
u(t)&=\lim_{T\nearrow\,\sup I}i\int_t^T e^{i(t-t')\Delta} F(u(t'))\,dt'\\
&=-\lim_{T\searrow\,\inf I}i\int_T^t e^{i(t-t')\Delta} F(u(t'))\,dt',
\end{aligned}
\end{equation}
as weak limits in $\dot H^1_x$.
\end{lemma}

The finite-time blowup scenario is considered in Section~\ref{S:finite time}.  Arguing as in \cite{kenig-merle},
we prove that the $L_x^2$ norm of $u(t)$ converges to zero as $t$ approaches the finite endpoint.
Since mass is conserved, this implies that $u$ is identically zero.

For the remaining two cases, we prove negative regularity in Section~\ref{S:neg}.  This is the heart of the matter
and is achieved in two stages.  First, we prove that the solution belongs to $L^\infty_t L^p_x$ for certain values
of $p$ less than $2d/(d-2)$.  This demonstrates that the solution decays more quickly at infinity than a general
$u\in L^\infty_t \dot H^1_x$; recall that since $u$ has uniformly bounded kinetic energy, $u\in L^\infty_t L^{2d/(d-2)}_x$ by
virtue of Sobolev embedding.  The proof of this first step involves a bootstrap argument built off the Duhamel formulae
\eqref{Duhamel}.  In order to disentangle frequency interactions, we make use of an `acausal' Gronwall inequality,
Lemma~\ref{L:Gronwall}.  The need for such a result in this paper stems from the two ways in which the nonlinearity
can produce low frequencies: from combinations of higher frequencies in $u$ and from fractional powers of lower frequencies in $u$.

The second step in proving negative regularity is to upgrade the decay proved in the first step to $L^2_x$-based spaces.
To do this, we take advantage of the global existence together with a double Duhamel trick in the spirit of \cite{tao:attractor}.
In order to make the associated time integrals converge, we need both $d\geq 5$ and the decay proved in step one.

In Section~\ref{S:cascade}, we use the negative regularity proved in Section~\ref{S:neg} together with the conservation of mass
to preclude the low-to-high frequency cascade.

In Section~\ref{S:Soliton}, we preclude the soliton.  To achieve this, we first use the negative regularity proved
in Section~\ref{S:neg} to deduce compactness properties for $u$ in $L_x^2$.  Secondly, we argue as in \cite{DHR,kenig-merle:wave}
to show that a minimal kinetic energy blowup solution must have zero momentum.  Notice that in order to even define the
momentum we need $u(t) \in \dot H^{1/2}_x$, which is considerable negative regularity compared to $\dot H^1_x$.
Using the vanishing of the momentum, we will deduce that the spatial center function obeys $|x(t)|=o(t)$ as $t\to \infty$, rather than
merely $O(t)$.  This mimics similar arguments in \cite{DHR,kenig-merle:wave} and relies crucially
on the $L^2_x$-compactness properties proved in the first step. To preclude the soliton, we now use a truncated virial
inequality much as in \cite{DHR, kenig-merle}; negative regularity and the fact that $|x(t)|=o(t)$ are needed in this
last step.

Proposition~\ref{P:blowup} is proved in Section~\ref{S:blowup}.  In Section~\ref{S:Conc}, we prove Theorem~\ref{T:conc}.

\subsection*{Acknowledgements}
We would like to thank Xiaoyi Zhang for permission to incorporate portions of some earlier
joint work \cite{RadialHD} into this article.  This unpublished manuscript gave a proof of
Theorem~\ref{T:main} in the case of radial data.  We would also like to acknowledge comments on
that manuscript from C. Kenig, as well as other helpful correspondence.  We are grateful to Terry Tao
for comments on an earlier version of this manuscript.

R.~K. was supported by NSF grants DMS-0701085 and DMS-0401277 and by a Sloan Foundation Fellowship.
Both authors were supported by the Institute for Advanced Study through NSF grant DMS-0635607.

Any opinions, findings and conclusions or recommendations expressed are those of the authors and do not reflect the views of the
National Science Foundation.

%
%
%
%

\section{Notations and useful lemmas}\label{S:Not}

\subsection{Some notation}
We use $X \lesssim Y$ or $Y \gtrsim X$ whenever $X \leq CY$ for some constant $C>0$.  We use $O(Y)$ to denote any quantity $X$
such that $|X| \lesssim Y$.  We use the notation $X \sim Y$ whenever $X \lesssim Y \lesssim X$.  The fact that these constants
depend upon the dimension $d$ will be suppressed.  If $C$ depends upon some additional parameters, we will indicate this with subscripts;
for example, $X \lesssim_u Y$ denotes the assertion that $X \leq C_u Y$ for some $C_u$ depending on $u$;
similarly for $X \sim_u Y$, $X = O_u(Y)$, etc.  We denote by $X\pm$ any quantity of the form $X\pm\eps$ for any $\eps>0$.

For any spacetime slab $I\times \R^d$, we use $L_t^qL_x^r(I\times \R^d)$ to denote the Banach space of functions $u: I\times
\R^d\to \mathbb C$ whose norm is
$$
\|u\|_{L_t^qL_x^r(I\times\R^d)}:=\Bigl(\int_I\|u(t)\|_{L^r_x}^q \, dt\Bigr)^{\frac 1q}<\infty,
$$
with the usual modifications when $q$ or $r$ are equal to infinity.  When $q=r$ we abbreviate $L^q_t L^q_x$ as $L^q_{t,x}$.

We define the Fourier transform on $\R^d$ by
$$
\hat f(\xi):= (2\pi)^{-d/2} \int_{\R^d} e^{-ix\xi}f(x)\,dx.
$$
For $s\in \R$, we define the fractional differentiation/integral operator
$$
\widehat{|\nabla|^s f}(\xi):=|\xi|^s\hat f(\xi),
$$
which in turn defines the homogeneous Sobolev norm
$$
\|f\|_{\dot H_x^s(\R^d)}:=\||\nabla|^s f\|_{L_x^2(\R^d)}.
$$

\subsection{Basic harmonic analysis}\label{ss:basic}
Let $\varphi(\xi)$ be a radial bump function supported in the ball $\{ \xi \in \R^d: |\xi| \leq \tfrac {11}{10} \}$ and equal to
$1$ on the ball $\{ \xi \in \R^d: |\xi| \leq 1 \}$.  For each number $N > 0$, we define the Fourier multipliers
\begin{align*}
\widehat{P_{\leq N} f}(\xi) &:= \varphi(\xi/N) \hat f(\xi)\\
\widehat{P_{> N} f}(\xi) &:= (1 - \varphi(\xi/N)) \hat f(\xi)\\
\widehat{P_N f}(\xi) &:= \psi(\xi/N)\hat f(\xi) := (\varphi(\xi/N) - \varphi(2\xi/N)) \hat f(\xi)
\end{align*}
and similarly $P_{<N}$ and $P_{\geq N}$.  We also define
$$ P_{M < \cdot \leq N} := P_{\leq N} - P_{\leq M} = \sum_{M < N' \leq N} P_{N'}$$
whenever $M < N$.  We will usually use these multipliers when $M$ and $N$ are \emph{dyadic numbers} (that is, of the form $2^n$
for some integer $n$); in particular, all summations over $N$ or $M$ are understood to be over dyadic numbers.  Nevertheless, it
will occasionally be convenient to allow $M$ and $N$ to not be a power of $2$.

Like all Fourier multipliers, the Littlewood-Paley operators commute with the propagator $e^{it\Delta}$, as well as with
differential operators such as $i\partial_t + \Delta$. We will use basic properties of these operators many many times,
including

\begin{lemma}[Bernstein estimates]\label{Bernstein}
 For $1 \leq p \leq q \leq \infty$,
\begin{align*}
\bigl\| |\nabla|^{\pm s} P_N f\bigr\|_{L^p_x(\R^d)} &\sim N^{\pm s} \| P_N f \|_{L^p_x(\R^d)},\\
\|P_{\leq N} f\|_{L^q_x(\R^d)} &\lesssim N^{\frac{d}{p}-\frac{d}{q}} \|P_{\leq N} f\|_{L^p_x(\R^d)},\\
\|P_N f\|_{L^q_x(\R^d)} &\lesssim N^{\frac{d}{p}-\frac{d}{q}} \| P_N f\|_{L^p_x(\R^d)}.
\end{align*}
\end{lemma}

We will also need the following fractional chain rule \cite{ChW:fractional chain rule}.  For a textbook treatment, see
\cite[\S 2.4]{Taylor:book}.

\begin{lemma}[Fractional chain rule, \cite{ChW:fractional chain rule}]\label{F Lip}
Suppose $G\in C^1(\mathbb C)$, $s \in (0,1]$, and $1<p,p_1,p_2<\infty$ are such that $\frac 1p=\frac 1{p_1}+\frac 1{p_2}$.   Then,
$$
\||\nabla|^sG(u)\|_p\lesssim \|G'(u)\|_{p_1}\||\nabla|^s u\|_{p_2}.
$$
\end{lemma}

\subsection{Strichartz estimates}
Let $e^{it\Delta}$ be the free Schr\"odinger evolution.  From the explicit formula
$$ e^{it\Delta} f(x) = \frac{1}{(4\pi i t)^{d/2}} \int_{\R^d} e^{i|x-y|^2/4t} f(y)\, dy,$$
one easily obtains the standard dispersive inequality
\begin{equation}\label{dispersive}
\| e^{it\Delta} f \|_{L_x^\infty(\R^d)} \lesssim|t|^{-\frac d2} \| f \|_{L_x^1(\R^d)}
\end{equation}
for all $t\neq 0$.

\begin{definition}[Admissible pairs]
For $d\geq3$, we say that a pair of exponents $(q,r)$ is \emph{Schr\"odinger-admissible} if
\begin{equation}\label{admissible}
\frac 2q+\frac dr=\frac d2 \quad \text{and} \quad 2\le q, r \le\infty.
\end{equation}
For a fixed spacetime slab $I\times \R^d$, we define the \emph{Strichartz norms}
$$
\|u\|_{\dot S^0(I)}:=\sup_{(q,r)\ \text{admissible}} \|u\|_{L_t^qL_x^r(\ird)} \qquad \text{and} \qquad \|u\|_{\dot
S^1(I)}:=\|\nabla u\|_{\dot S^0(I)}.
$$
We write $\dot S^0(I)$ and $\dot S^1(I)$ for the closure of all test functions under these norms, respectively.
\end{definition}

A simple application of Sobolev embedding yields
\begin{align*}
\|\nabla u\|_{L_t^\infty L_x^2(\ird)} & + \|\nabla u\|_{L_{t,x}^{\frac{2(d+2)}d}(\ird)}
    + \|\nabla u\|_{L_t^2 L_x^{\frac {2d}{d-2}}(\ird)}\\
&+\|u\|_{L_t^\infty L_x^{\frac {2d}{d-2}}(\ird)} + \|u\|_{L_{t,x}^\sct(\ird)} \lesssim \|u\|_{\dot S^1(I)}
\end{align*}
for all $d\geq 3$.

As a consequence of the dispersive estimate \eqref{dispersive}, we have the following standard Strichartz estimate.

\begin{lemma}[Strichartz]\label{strichartz}
Let $k= 0,1$, let $I$ be a compact time interval, and let $u: I\times\R^d \to \mathbb C$ be a solution to the forced
Schr\"odinger equation
$$
iu_t+\Delta u=F.
$$
Then,
$$
\|u\|_{\dot S^k(I)}\lesssim \|u(t_0)\|_{\dot H_x^k}+\|F\|_{\dot N^k(I)}
$$
for any $t_0\in I$.
\end{lemma}

\begin{proof}
See, for example, \cite{gv:strichartz, strichartz}.  For the endpoint $(q,r)=\bigl(2,\frac{2d}{d-2}\bigr)$, see \cite{tao:keel}.
\end{proof}

The next result is Lemma~3.7 from \cite{keraani-h1} extended to all dimensions.  This will play an important role in the proof
of Lemma~\ref{L:bad profile}.  Below we offer a quantitative and more streamlined proof.

\begin{lemma}\label{L:Keraani3.7}
Given $\phi\in \dot H^1_x(\R^d)$,
$$
\| \nabla e^{it\Delta} \phi \|_{L^2_{t,x}([-T,T]\times\{|x|\leq R\})}^3 \lesssim
     T^{\frac2{d+2}} R^{\frac{3d+2}{2(d+2)}} \| e^{it\Delta} \phi \|_{L^{2(d+2)/(d-2)}_{t,x}}\| \nabla \phi \|_{L^2_x}^2.
$$
\end{lemma}

\begin{proof}
Given $N>0$, H\"older's and Bernstein's inequalities imply
\begin{align*}
\| \nabla e^{it\Delta} \phi_{< N} \|_{L^2_{t,x}([-T,T]\times\{|x|\leq R\})}
    &\lesssim T^{2/(d+2)} R^{2d/(d+2)} \| e^{it\Delta} \nabla \phi_{< N} \|_{L^{2(d+2)/(d-2)}_{t,x}} \\
&\lesssim T^{2/(d+2)} R^{2d/(d+2)} \, N\, \| e^{it\Delta} \phi \|_{L^{2(d+2)/(d-2)}_{t,x}}.
\end{align*}
On the other hand, the high frequencies can be estimated using local smoothing:
\begin{align*}
\| \nabla e^{it\Delta} \phi_{\geq N} \|_{L^2_{t,x}([-T,T]\times\{|x|\leq R\})}
    &\lesssim R^{1/2} \| |\nabla|^{1/2} \phi_{\geq N} \|_{L^2_x} \\
&\lesssim  N^{-1/2} R^{1/2} \| \nabla \phi \|_{L^2_x}.
\end{align*}
The lemma now follows by optimizing the choice of $N$.
\end{proof}

We will also make use of the following bilinear estimate:

\begin{lemma}[Bilinear Strichartz]\label{L:bilinear strichartz}
For any spacetime slab $I \times \R^d$ and any $M, N > 0$, we have
\begin{align*}
\|e^{it\Delta} \phi_N e^{it\Delta} \phi_M \|_{L^2_{t,x}(I \times \R^d)}
& \lesssim M^{\frac{d-4}2} N^{-1} \|\nabla \phi_M\|_{L_x^2} \|\nabla \phi_N\|_{L_x^2},
\end{align*}
for any function $\phi$.
\end{lemma}

\begin{proof}
See \cite[Lemma 2.5]{Monica:thesis}, which builds on earlier versions in \cite{borg:book, ckstt:gwp}.
\end{proof}

\subsection{Concentration compactness}\label{SS:cc}
In this subsection we record the linear profile decomposition statement due to Keraani \cite{keraani-h1},
which will lead to the reduction in Theorem~\ref{T:reduct}.
We first recall the symmetries of the equation \eqref{nls} which fix the initial surface $t=0$ and preserve
the energy.

\begin{definition}[Symmetry group]\label{D:sym}
For any phase $\theta \in \R/2\pi \Z$, position $x_0\in \R^d$, and scaling parameter $\lambda > 0$, we define a unitary transformation
$g_{\theta,x_0,\lambda}: \dot H^1_x(\R^d) \to \dot H^1_x(\R^d)$ by
$$
[g_{\theta,x_0, \lambda} f](x) :=  \lambda^{-\frac{d-2}2} e^{i\theta}  f\bigl( \lambda^{-1}(x-x_0) \bigr).
$$
Let $G$ denote the collection of such transformations.  For a function $u: I \times \R^d \to \C$, we define $T_{g_{\theta,x_0,
\lambda}} u: \lambda^2 I \times \R^d \to \C$ where $\lambda^2 I := \{ \lambda^2 t: t \in I \}$ by the formula
$$
[T_{g_{\theta,x_0, \lambda}} u](t,x) :=  \lambda^{-\frac{d-2}2} e^{i\theta} u\bigl( \lambda^{-2}t, \lambda^{-1}(x-x_0)\bigr).
$$
Note that if $u$ is a solution to \eqref{nls}, then $T_{g}u$ is a solution to \eqref{nls} with initial data $g u_0$.
\end{definition}

\begin{remark} It is easy to verify that $G$ is a group and that the map $g \mapsto T_g$ is a homomorphism.
The map $u \mapsto T_g u$ maps solutions to \eqref{nls} to solutions with the same energy and scattering size as $u$,
that is, $E(T_g u) = E(u)$ and $S(T_g u)=S(u)$.  Furthermore, $u$ is a maximal-lifespan solution if and only if $T_g u$
is a maximal-lifespan solution.
\end{remark}

We are now ready to state Keraani's linear profile decomposition.

\begin{lemma} [Linear profile decomposition, \cite{keraani-h1}]\label{L:cc}
Fix $d\ge 3$ and let $\{u_n\}_{n\geq 1}$ be a sequence of functions bounded in $\dot H_x^1(\R^d)$. Then,
after passing to a subsequence if necessary, there exist a sequence of functions $\{\phi^j\}_{j\geq 1}\subset \dot H_x^1(\R^d)$,
group elements $\gnj \in G$, and times $\tnj\in \R$ such that we have the decomposition
\begin{align}\label{decomp}
u_n = \sum_{j=1}^J \gnj e^{i\tnj\Delta}\phi^j + \wnJ
\end{align}
for all $J\geq 1$; here, $\wnJ \in \dot H^1_x(\R^d)$ obey
\begin{equation}\label{w scat}
\lim_{J\to \infty}\limsup_{n\to\infty} \bigl\| e^{it\Delta}\wnJ
\bigr\|_{L_{t,x}^{\frac{2(d+2)}{d-2}}(\R\times\R^d)}=0.
\end{equation}
Moreover, for any $j \neq j'$,
\begin{align}\label{crazy}
\frac{\lambda_n^j}{\lambda_n^{j'}} + \frac{\lambda_n^{j'}}{\lambda_n^{j}}
    + \frac{|x_n^j-x_n^{j'}|^2}{\lambda_n^j \lambda_n^{j'}}
    + \frac{\bigl|t_n^j(\lambda_n^j)^2- t_n^{j'}(\lambda_n^{j'})^2\bigr|}{\lambda_n^j\lambda_n^{j'}}\to\infty
    \quad \text{as } n\to \infty.
\end{align}
Furthermore, for any $J \geq 1$ we have the kinetic energy decoupling property
\begin{equation}\label{decouple}
\lim_{n \to \infty} \Bigl[ \|\nabla u_n\bigr\|_2^2 - \sum_{j=1}^J \|\nabla \phi^j\|_2^2 - \|\nabla \wnJ\|^2_2 \Bigr] = 0.
\end{equation}
\end{lemma}

Finally, we will need the following result, which shows that for all $J\geq 1$, the error term $\wnJ$ converges weakly to zero
in $\dot H^1_x(\R^d)$ modulo the symmetries associated to $\phi^j$ for $1\leq j\leq J$.
This property is actually built into the proof of Lemma~\ref{L:cc}; however, since it is not explicitly stated in \cite{keraani-h1}
and is easy to verify \emph{a posteriori}, we do that here.

\begin{lemma}[Strong decoupling]\label{L:strong decouple}
For all $J\geq 1$ and all $1\leq j\leq J$, the sequence $e^{-i\tnj\Delta}[(\gnj)^{-1}\wnJ]$ converges weakly to zero in $\dot
H^1_x(\R^d)$ as $n\to \infty$.  In particular, this implies the kinetic energy decoupling \eqref{decouple}.
\end{lemma}

\begin{proof}
Fix $J\geq 1$ and $1\leq j\leq J$.  By \eqref{decouple} and the fact that $\{u_n\}_{n\geq 1}$ is bounded in $\dot H^1_x(\R^d)$,
we deduce that $\{e^{-i\tnj\Delta}[(\gnj)^{-1}\wnJ]\}_{n\geq 1}$ is bounded in $\dot H^1_x(\R^d)$.  Using Alaoglu's Theorem (and
passing to a subsequence if necessary), we obtain that $e^{-i\tnj\Delta}[(\gnj)^{-1}\wnJ]$ converges weakly in $\dot
H^1_x(\R^d)$ to some $\psi\in \dot H_x^1(\R^d)$.  To prove the lemma, it thus suffices to show that $\psi\equiv 0$.

By weak convergence and \eqref{decomp},
\begin{equation}\label{decomp2}
\begin{aligned}
\|\psi\|_{\dot H^1_x(\R^d)}^2
&=\lim_{n\to \infty} \bigl\langle \nabla e^{-i\tnj\Delta}[(\gnj)^{-1}\wnJ], \nabla \psi \bigr\rangle\\
&= \sum_{l=J+1}^L \lim_{n\to \infty} \bigl\langle \nabla g_n^l e^{it_n^l\Delta}\phi^l, \nabla g_n^j e^{it_n^j\Delta}\psi \bigr\rangle\\
&\quad + \lim_{n\to \infty} \bigl\langle \nabla e^{-i\tnj\Delta}[(\gnj)^{-1} w_n^L], \nabla \psi \bigr\rangle,
\end{aligned}
\end{equation}
for all $L>J$.  The limits on the right-hand side are guaranteed to exist; indeed, using \eqref{crazy}, a change of variables
shows
$$
\lim_{n\to \infty} \bigl\langle \nabla g_n^l e^{it_n^l\Delta}\phi^l, \nabla g_n^j e^{it_n^j\Delta}\psi \bigr\rangle\to 0 \quad
\text{as} \quad n\to \infty
$$
for all $L\geq l\geq J+1>j$; see the proof of \cite[Lemma~2.7]{keraani-h1}.

On the other hand, combining the fact that the family $\{e^{-i\tnj\Delta}[(\gnj)^{-1} w_n^L]\}_{n,L\geq1}$ is bounded in $\dot
H^1_x(\R^d)$ with
$$
\lim_{L\to \infty}\limsup_{n\to \infty} S_{\R}\bigl( e^{it\Delta} e^{-i\tnj\Delta}[(\gnj)^{-1} w_n^L]\bigr) =\lim_{L\to
\infty}\limsup_{n\to \infty} S_{\R}\bigl( e^{it\Delta} w_n^L\bigr)=0,
$$
we deduce that $e^{-i\tnj\Delta}[(\gnj)^{-1} w_n^L]$ converges weakly to zero in $\dot H^1_x(\R^d)$ as $n,L\to \infty$
(cf. \cite[Lemma~3.63]{merle-vega}).
Thus, for $L$ sufficiently large,
$$
\limsup_{n\to \infty} \bigl|\bigl\langle \nabla e^{-i\tnj\Delta}[(\gnj)^{-1} w_n^L], \nabla \psi \bigr\rangle\bigr|
    \leq \tfrac 12 \|\psi\|_{\dot H^1_x(\R^d)}^2.
$$

Returning to \eqref{decomp2} and choosing $L$ large, we conclude $\psi\equiv 0$.  This finishes the proof of Lemma~\ref{L:strong
decouple}.
\end{proof}

\subsection{Additional harmonic analysis}
In this subsection we record tools that will be used to prove the concentration result Theorem~\ref{T:conc}.

The next lemma is an inverse Strichartz inequality.  It roughly states that the Strichartz norm cannot be large without there
being a bubble of concentration in spacetime.  Results of this type constitute an important precursor to the concentration
compactness technique.  The prototype is \cite[\S2--3]{bourg.2d}, which has been extended and elaborated in many subsequent
papers; see, for instance, \cite{BegoutVargas,borg:scatter,keraani-h1,keraani-l2,merle-vega,tao: gwp radial,tao:pseudo}.
While the lemma can be deduced a posteriori from these works, we give a self-contained argument using the ideas in them.

\begin{lemma}[Inverse Strichartz]\label{L:B conc}
Fix $d\geq 3$.  Let $\phi\in \dot H^1_x(\R^d)$ and $\eta>0$ such that
$$
\int_I \int_{\R^d} \bigl|e^{it\Delta} \phi \bigr|^{\frac{2(d+2)}{d-2}} \,dx\,dt \geq \eta
$$
for some interval $I\subseteq \R$.  Then there exist $C=C(\|\nabla \phi\|_2, \eta)$, $x_0\in \R^d$, and $J\subseteq I$ so that
\begin{equation*}
\int_{|x - x_0|\leq C |J|^{1/2}} \bigl|e^{it\Delta} \nabla \phi \bigr|^2 \,dx \geq C^{-1} \quad\text{for all} \quad t\in J.
\end{equation*}
Notice that $C$ does not depend on $I$ or $J$.
\end{lemma}

\begin{proof}
First we prove that
\begin{equation}\label{single M big}
\int_I \int_{\R^d} \bigl|e^{it\Delta} \phi_M \bigr|^{\frac{2(d+2)}{d-2}} \,dx\,dt \gtrsim \eta^c
\end{equation}
for some dyadic $M \gtrsim |I|^{-1/2}$ and some (dimension-dependent) $c>0$.  We begin with the argument for dimensions $d\geq 6$.

Using elementary Littlewood--Paley theory together with the Strichartz inequality and the bilinear Strichartz
inequality (Lemma~\ref{L:bilinear strichartz}), we argue as follows:
\begin{align}
\eta &\lesssim \Bigl\| \Bigl( \sum_M \bigl| e^{it\Delta} \phi_M\bigr|^2 \Bigr)^{\frac{d+2}{2(d-2)}}
        \Bigl( \sum_N \bigl| e^{it\Delta} \phi_N\bigr|^2 \Bigr)^{\frac{d+2}{2(d-2)}} \Bigr\|_{L^1_{t,x}} \notag \\
&\lesssim \Bigl\| \Bigl( \sum_M \bigl| e^{it\Delta} \phi_M\bigr|^{\frac{d+2}{d-2}} \Bigr)
        \Bigl( \sum_N \bigl| e^{it\Delta} \phi_N\bigr|^{\frac{d+2}{d-2}} \Bigr) \Bigr\|_{L^1_{t,x}} \notag \\
&\lesssim \sum_{M\leq N} \|e^{it\Delta} \phi_M e^{it\Delta} \phi_N\|_{L_{t,x}^2}^{\frac2{d-2}}
        \|e^{it\Delta} \phi_N\|_{L_{t,x}^{\frac{2(d+2)}d}}
        \|e^{it\Delta} \phi_N\|_{L_{t,x}^{\frac{2(d+2)}{d-2}}}^{\frac2{d-2}} \times \notag \\
&\qquad \qquad \qquad \qquad \qquad \times
    \|e^{it\Delta} \phi_M\|_{L_{t,x}^\infty}^{\frac{8}{(d-2)^2}}
    \|e^{it\Delta} \phi_M\|_{L_{t,x}^{\frac{2(d+2)}{d-2}}}^{\frac{(d+2)(d-4)}{(d-2)^2}} \notag \\
&\lesssim \sum_{M\leq N}\bigl(\tfrac{M}{N}\bigr)^{\frac d{d-2}} \|\nabla \phi_N\|_{L_x^2}^{\frac{d+2}{d-2}}
     \|\nabla \phi_M\|_{L_x^2}^{\frac{d-6}{d-2}}\|e^{it\Delta} \phi_M\|_{L_{t,x}^{\frac{2(d+2)}{d-2}}}^{\frac8{d-2}} \notag \\
&\lesssim \|\nabla \phi\|_{L_x^2}^2 \sup_{M} \|e^{it\Delta} \phi_M\|_{L_{t,x}^{\frac{2(d+2)}{d-2}}}^{\frac8{d-2}},
\label{big calc}
\end{align}
where all space-time norms are on $I\times\R^d$.

On the other hand, by Bernstein's inequality (Lemma~\ref{Bernstein}),
$$
\int_I \int_{\R^d} \bigl|e^{it\Delta} \phi_{M} \bigr|^{\frac{2(d+2)}{d-2}} \,dx\,dt
    \lesssim |I| M^2 \|\nabla \phi\|_{L^2_x}^{\frac{2(d+2)}{d-2}}.
$$
Combining this with the argument above, we see that there is $M\gtrsim_{\|\nabla \phi\|_2, \eta} |I|^{-1/2}$ so that
\eqref{single M big} holds with $c=(d-2)/8$.

To obtain \eqref{single M big} in dimensions $3\leq d <6$, we merely need to find a replacement for \eqref{big calc}.
We argue as follows using the same tools as before:
\begin{align*}
\eta
&\lesssim \|e^{it\Delta} \phi\|_{L_{t,x}^{\frac{2(d+2)}{d-2}}}^{\frac{2(6-d)}{d-2}}
        \Bigl\| \Bigl( \sum_M \bigl| e^{it\Delta} \phi_M\bigr|^2 \Bigr)
        \Bigl( \sum_N \bigl| e^{it\Delta} \phi_N\bigr|^2 \Bigr) \Bigr\|_{L^{\frac{d+2}{2(d-2)}}_{t,x}}\\
&\lesssim \|\nabla \phi\|_{L_x^2}^{\frac{2(6-d)}{d-2}} \sum_{M\leq N} \|e^{it\Delta} \phi_M\|_{L^{\frac{2(d+2)}{d-2}}_{t,x}}
        \|e^{it\Delta} \phi_N\|_{L^{\frac{2(d+2)}{d-2}}_{t,x}}
        \|e^{it\Delta} \phi_M e^{it\Delta} \phi_N\|_{L_{t,x}^2}^{\frac {2(d-2)}{d+2}}\times\\
&\qquad \qquad \qquad \qquad \qquad \times
        \|e^{it\Delta} \phi_M \|_{L_{t,x}^\infty}^{\frac {6-d}{d+2}} \|e^{it\Delta}\phi_N\|_{L_{t,x}^\infty}^{\frac {6-d}{d+2}}\\
&\lesssim \|\nabla \phi\|_{L_x^2}^{\frac{2(6-d)}{d-2}} \sum_{M\leq N} \bigl( \tfrac {M}{N}\bigr)^{\frac{(d-2)^2}{2(d+2)}}
        \|\nabla \phi_N \|_{L_x^2}^2  \|\nabla \phi_M \|_{L_x^2} \|e^{it\Delta} \phi_M\|_{L^{\frac{2(d+2)}{d-2}}_{t,x}}\\
&\lesssim \|\nabla \phi\|_{L_x^2}^{\frac{d+6}{d-2}} \sup_M \|e^{it\Delta} \phi_M\|_{L^{\frac{2(d+2)}{d-2}}_{t,x}}.
\end{align*}

Having proved \eqref{single M big}, we continue as follows:  Using Bernstein combined with the Strichartz inequality,
we obtain the upper bound
$$
\|e^{it\Delta} \phi_{M}\|_{L_{t,x}^{\frac{2(d+2)}d}(I\times\R^d)}\lesssim M^{-1}\|\nabla \phi\|_{L_x^2};
$$
this combined with \eqref{single M big} and H\"older's inequality yields
$$
\|e^{it\Delta} \phi_{M}\|_{L_{t,x}^\infty(I\times\R^d)}\gtrsim_{\|\nabla \phi\|_2, \eta} M^{\frac{d-2}2}.
$$
Thus, there exist $t_0\in I$ and $x_0\in\R^d$ so that
$$
\bigl| [e^{it_0\Delta}\phi_M](x_0) \bigr| \gtrsim_{\|\nabla \phi\|_2, \eta} M^{\frac{d-2}2}.
$$
Using basic properties of the kernel of $e^{it\Delta} P_M$, we may deduce
$$
\int_{|x-x_0|\lesssim M^{-1}} \bigl| e^{it\Delta} \nabla \phi (x) \bigr|^2\,dx \gtrsim_{\|\nabla \phi\|_2,\eta}  1
$$
for all $|t-t_0|\lesssim M^{-2}$.  Let $J:=\{t\in I:\, |t-t_0|\lesssim M^{-2}\}$. To obtain the claim,
we simply note that because of our lower bound on $M$, the length of $J$ obeys $|J|\gtrsim_{\|\nabla \phi\|_2,\eta} M^{-2}$.
\end{proof}

Next, we recall \cite[Lemma~10.2]{KTV}.  While \cite[Lemma~10.2]{KTV} is stated
and proved in dimension $d=2$, the proof extends without difficulty to higher dimensions.

\begin{lemma}[Tightness of profiles]\label{L:comp prof}
Let $d\geq 3$ and let $\psi\in \dot H^1_x(\R^d)$.  Assume that
$$
\int_{|x-x_k|\leq r_k} \bigl|e^{it_k\Delta} \nabla \psi \bigr|^2 \,dx \geq \eps
$$
for some $\eps>0$ and sequences $t_k\in\R$, $x_k \in \R^d$, and $r_k>0$.  Then for any sequence $a_k\to\infty$,
\begin{equation*}
\int_{|x|\leq a_k r_k} \bigl|e^{it_k\Delta} \nabla \psi \bigr|^2 \,dx \to \|\nabla \psi\|_2^2.
\end{equation*}
\end{lemma}

As the kinetic energy is not conserved, we need to upgrade this lemma as follows:

\begin{proposition}[Tightness of trajectories]\label{P:comp prof}
Let $\psi:I\times\R^d\to\C$ be a solution to \eqref{nls} with $S_I(\psi)<\infty$.  Suppose
$$
\int_{|x-x_k|\leq r_k} \bigl|e^{it_k\Delta} \nabla \psi(\tau_k) \bigr|^2 \,dx \geq \eps
$$
for some $\eps>0$ and sequences $t_k\in\R$, $x_k \in \R^d$, $\tau_k\in I$, and $r_k>0$.  Then
\begin{equation*}
\Bigl|\|\nabla \psi(\tau_k)\|_2^2- \int_{|x|\leq a_k r_k} \bigl|e^{it_k\Delta} \nabla \psi(\tau_k) \bigr|^2 \,dx \Bigr|\to 0
\end{equation*}
for any sequence $a_k\to\infty$.
\end{proposition}

\begin{proof}
It suffices to treat the case where the sequence $\tau_k$ converges (possibly to $\pm\infty$).  By Theorem~\ref{T:local}, we may
assume that $I$ is closed.

If $\tau_k$ converges to a finite point (in $I$), then the claim follows from Lemma~\ref{L:comp prof} and the $\dot
H^1_x$-continuity of the flow.

Next we treat the case $\tau_k\to\infty$; a similar argument settles the case $\tau_k\to-\infty$.  In particular, $\sup I
=\infty$. As $\psi$ has finite scattering size on $I$, Theorem~\ref{T:local} implies the existence of $\psi_+\in \dot H^1_x$ so
that
$$
\bigl\| \psi(\tau_k) - e^{i\tau_k\Delta} \psi_+ \bigr\|_{\dot H^1_x} \to 0.
$$
We may now apply Lemma~\ref{L:comp prof} to complete the argument.
\end{proof}

\subsection{A Gronwall inequality}

Our last technical tool is the most elementary.  It is a form of Gronwall's inequality that involves both the past
and the future, `acausal' in the terminology of \cite{tao:book}.  It will be used in Section~\ref{S:neg}.

\begin{lemma}\label{L:Gronwall}
Given $\gamma>0$, $0<\eta<\tfrac12(1-2^{-\gamma})$, and $\{b_k\}\in\ell^\infty(\Z^+)$,
let $x_k\in\ell^\infty(\Z^+)$ be a non-negative sequence obeying
\begin{align}\label{Gron rec}
x_k \leq b_k + \eta \sum_{l=0}^\infty 2^{-\gamma|k-l|} x_l \qquad \text{for all $k\geq 0$.}
\end{align}
Then
\begin{align}\label{Gron bound}
x_k \lesssim \sum_{l=0}^{k} r^{|k-l|} b_l  \qquad \text{for all $k\geq 0$}
\end{align}
for some $r=r(\eta)\in (2^{-\gamma},1)$.  Moreover, $r\downarrow 2^{-\gamma}$ as $\eta\downarrow 0$.
\end{lemma}

\begin{proof}
Our proof follows a well-travelled path.
By decreasing entries in $b_k$ we can achieve equality in \eqref{Gron rec}; since this also reduces the righthand side
of \eqref{Gron bound}, it suffices to prove the lemma in this case.  Note that since $x_k\in\ell^\infty$, $b_k$ will remain
a bounded sequence.

Let $A$ denote the doubly infinite matrix with entries $A_{k,l}=2^{-\gamma|k-l|}$ and let $P$ denote the natural projection
from $\ell^2(\Z)$ onto $\ell^2(\Z^+)$.  Our goal is to show that \eqref{Gron bound} holds for any solution of
\begin{equation}\label{groneq}
(1-\eta PAP^*)x =b.
\end{equation}
First we observe that since
$$
\|A\|=\sum_{k\in\Z} 2^{-\gamma|k|} = \frac{1+2^{-\gamma}}{1-2^{-\gamma}},
$$
$\eta A$ is a contraction on $\ell^\infty$.  Thus we may write
$$
x = \sum_{p=0}^\infty (\eta PAP^*)^p b \leq \sum_{p=0}^\infty P (\eta A)^p P^* b = P (1-\eta A)^{-1} P^* b,
$$
where the inequality is meant entry-wise.  The justification for this inequality is simply that the matrix $A$
has non-negative entries.  We will complete the proof of \eqref{Gron bound} by computing the entries of $(1-\eta A)^{-1}$.
This is easily done via Fourier methods:  Let
$$
a(z) := \sum_{k\in\Z} 2^{-\gamma|k|} z^k = 1 + \frac{2^{-\gamma}z}{1-2^{-\gamma}z} + \frac{2^{-\gamma}z^{-1}}{1-2^{-\gamma}z^{-1}}
$$
and
\begin{align*}
f(z):= \frac{1}{1-\eta a(z)} &= \frac{(z-2^\gamma)(z-2^{-\gamma})}{z^2-(2^{-\gamma}+2^{\gamma}-\eta2^{\gamma}+\eta2^{-\gamma})z+1} \\
&=1 + \frac{(1-r2^{-\gamma})(r2^{\gamma}-1)}{(1-r^2)} \Bigl[ 1 + \frac{rz}{1-rz} + \frac{rz^{-1}}{1-rz^{-1}}\Bigr],
\end{align*}
where $r\in(0,1)$ and $1/r$ are the roots of $z^2-(2^{-\gamma}+2^{\gamma}-\eta2^{\gamma}+\eta2^{-\gamma})z+1=0$.  From this formula,
we can immediately read off the Fourier coefficients of $f$, which give us the matrix elements of $(1-\eta A)^{-1}$.  In particular,
they are $O(r^{|k-l|})$.
\end{proof}

%
%
%
%

\section{Reduction to almost periodic solutions}\label{S:Reduct}

The goal of this section is to prove Theorem~\ref{T:reduct}.  In order to achieve this, we will first prove a Palais-Smale
condition modulo symmetries.

For any $0\leq E_0\le \|\nabla W\|_2^2$, we define
$$
L(E_0):=\sup \{S(u):\, u:\ird\to \C \text{ such that } \sup_{t\in I}\|\nabla u(t)\|_2^2 \leq  E_0\},
$$
where the supremum  is taken over all solutions $u:\ird\to \C$ to \eqref{nls} obeying $\sup_{t\in I}\|\nabla u(t)\|_2^2 \leq E_0$.
Thus, $L:\bigl[0, \|\nabla W\|_2^2\bigr] \to [0, \infty]$ is a non-decreasing function with
$L\bigl(\|\nabla W\|_2^2\bigr)=\infty$. Moreover, from Theorem~\ref{T:local},
\begin{align*}
L(E_0)\lesssim_d E_0^{\frac{d+2}{d-2}} \quad \text{for} \quad  E_0\leq \eta_0,
\end{align*}
where $\eta_0=\eta_0(d)$ is the threshold from the small data theory.

From Lemma~\ref{L:stab}, we see that $L$ is continuous.  Therefore, there must exist a unique \emph{critical kinetic energy}
$E_c$ such that $L(E_0)<\infty$ for $E_0<E_c$ and $L(E_0)=\infty$ for $E_0\geq E_c$.  In particular, if $u:\ird\to \C$ is a
maximal-lifespan solution to \eqref{nls} such that $\sup_{t\in I}\|\nabla u(t)\|_2^2 < E_c$, then $u$ is global and
$$
S_\R(u)\leq L\bigl(\sup_{t\in I}\|\nabla u(t)\|_2^2\bigr).
$$

Failure of Conjecture~\ref{conj} is equivalent to $0 < E_c < \|\nabla W\|_2^2$.

\subsection{The key convergence result}
In this subsection we prove the folowing

\begin{proposition}[Palais-Smale condition modulo symmetries]\label{P:palais-smale}
Fix $d\geq 3$.  Let $u_n:I_n\times\R^d\mapsto \C$ be a sequence of solutions to \eqref{nls} such that
\begin{align}\label{max ke}
\limsup_{n\to \infty} \sup_{t\in I_n}\|\nabla u_n(t)\|_2^2 =E_c
\end{align}
and let $t_n\in I_n$ be a sequence of times such that
\begin{equation*}
\lim_{n\to \infty} S_{\ge t_n}(u_n) = \lim_{n\to \infty} S_{\le t_n}(u_n) = \infty.
\end{equation*}
Then the sequence $u_n(t_n)$ has a subsequence which converges in $\dot H^1_x(\R^d)$ modulo symmetries.
\end{proposition}

\begin{proof}
Using the time-translation symmetry of \eqref{nls}, we may set $t_n=0$ for all $n\geq 1$.  Thus,
\begin{equation}\label{blow up in two}
\lim_{n\to \infty} S_{\ge 0} (u_n) = \lim_{n\to \infty} S_{\le 0}(u_n) = \infty.
\end{equation}

Applying Lemma~\ref{L:cc} to the sequence $u_n(0)$ (which is bounded in $\dot H^1_x(\R^d)$ by \eqref{max ke}) and passing to a
subsequence if necessary, we obtain the decomposition
$$
u_n(0) = \sum_{j=1}^J g_n^j e^{i\tnj\Delta} \phi^j + w_n^J
$$
as in Lemma~\ref{L:cc}.

Refining the subsequence once for every $j$ and using a diagonal argument, we may assume that for each $j$, the sequence
$\{\tnj\}_{n\geq 1}$ converges to some $t^j\in [-\infty, \infty]$.  If $t^j\in (-\infty, \infty)$, then by replacing $\phi^j$ by
$e^{i t^j \Delta}\phi^j$, we may assume that $t^j=0$; moreover, absorbing the error $e^{i \tnj \Delta}\phi^j - \phi^j$ into the
error term $\wnJ$, we may assume that $\tnj\equiv 0$. Thus, either $\tnj\equiv 0$ or $\tnj\to\pm \infty$.

We now define the nonlinear profiles $v^j:I^j\times\R^d \to \C$ associated to $\phi^j$ and $\tnj$ as follows:

\begin{CI}
\item If $\tnj\equiv 0$, then $v^j$ is the maximal-lifespan solution to \eqref{nls} with initial data $v^j(0)=\phi^j$.
\item If $\tnj\to \infty$, then $v^j$ is the maximal-lifespan solution to \eqref{nls} that scatters forward in time to $e^{it\Delta}\phi^j$.
\item If $\tnj\to -\infty$, then $v^j$ is the maximal-lifespan solution to \eqref{nls} that scatters backward in time to $e^{it\Delta}\phi^j$.
\end{CI}

For each $j,n\geq 1$, we introduce $\vnj:I_n^j\times\R^d\to \C$ defined by
$$
\vnj(t):= T_{\gnj}\bigl[ v^j(\cdot + \tnj)\bigr](t),
$$
where $I_n^j:=\{t\in \R:\, (\lambda_n^j)^{-2} t + \tnj \in I^j\}$. Each $\vnj$ is a solution to \eqref{nls} with initial data at
time $t=0$ given by $\vnj(0)=\gnj v^j(\tnj)$ and maximal lifespan $I_n^j= (-T^-_{n,j}, T^+_{n,j})$, where $-\infty\leq
-T^-_{n,j}<0<T^+_{n,j}\leq \infty$.

By \eqref{decouple}, there exists $J_0\geq 1$ such that
$$
\|\nabla \phi^j \|_2\leq \eta_0 \quad \text{for all} \quad j\geq J_0,
$$
where $\eta_0=\eta_0(d)$ is the threshold for the small data theory.  Hence, by Theorem~\ref{T:local}, for all $n\geq 1$ and all
$j\geq J_0$ the solutions $\vnj$ are global and moreover,
\begin{align}\label{tail}
\sup_{t\in \R} \|\nabla \vnj(t)\|_2^2 + S_\R(\vnj)\lesssim \|\nabla \phi^j \|_2^2.
\end{align}

\begin{lemma}[At least one bad profile]\label{L:bad profile}
There exists $1\leq j_0<J_0$ such that
$$
\limsup_{n\to \infty} S_{[0, T^+_{n,j_0})}(v_n^{j_0})=\infty.
$$
\end{lemma}

\begin{proof}
Assume for a contradiction that for all $1\leq j<J_0$,
\begin{align}\label{S ass}
\limsup_{n\to \infty} S_{[0, T^+_{n,j})}(v_n^j)< \infty.
\end{align}
In particular, this implies $T^+_{n,j}=\infty$ for all $1\leq j<J_0$ and all sufficiently large $n$.  Moreover, subdividing
$[0,\infty)$ into intervals where the scattering size of $\vnj$ is small, applying the Strichartz inequality on each such
interval, and then summing, we obtain
\begin{align}\label{S1 ass}
\limsup_{n\to \infty} \|\vnj\|_{\dot S^1([0,\infty))}< \infty
\end{align}
for all $1\leq j<J_0$.

Combining \eqref{S ass} with \eqref{tail}, and then using \eqref{decouple} and \eqref{max ke},
\begin{equation}\label{S decouple}
\sum_{j\geq 1} S_{[0,\infty)}(\vnj)\lesssim 1 + \sum_{j\geq J_0}\|\nabla \phi^j \|_2^2 \lesssim 1 + E_c
\end{equation}
for all $n$ sufficiently large.

From these assumptions, we will deduce a bound on the scattering size of $u_n$ forward in time (for $n$ sufficiently large),
thus contradicting \eqref{blow up in two}.  In order to achieve this, we will use Lemma~\ref{L:stab}.  To this end, we define the
approximation
$$
\unJ(t):=\sum_{j=1}^J \vnj(t) + e^{it\Delta}\wnJ.
$$
Note that
\begin{align*}
\|\unJ(0)-u_n(0)\|_{\dot H^1_x(\R^d)}
&\lesssim \bigl\| \sum_{j=1}^J \bigl( \gnj v^j(\tnj) - \gnj e^{i\tnj\Delta} \phi^j\bigr) \bigr\|_{\dot H^1_x(\R^d)}\\
&\lesssim \sum_{j=1}^J \bigl\| v^j(\tnj) - e^{i\tnj\Delta} \phi^j\bigr\|_{\dot H^1_x(\R^d)},
\end{align*}
and hence, by our choice of $v^j$,
$$
\limsup_{n\to \infty}\|u_n(0)-\unJ(0)\|_{\dot H^1_x(\R^d)} =0.
$$

We now show that $\unJ$ does not blowup forward in time.  Indeed, by \eqref{crazy} and the fact that $\vnj$ does not blow up
forward in time for any $j\geq1$ and all $n$ sufficiently large,
\begin{align*}
\limsup_{n\to \infty} S_{[0,\infty)}\bigl(|\vnj|^{1-\theta} |v_n^{j'}|^\theta \bigl) = 0
\end{align*}
for any $0<\theta<1$ and $j\neq j'$; see \cite{keraani-h1}.  Thus, by \eqref{w scat} and \eqref{S decouple},
\begin{equation}\label{S unj}
\begin{aligned}
\lim_{J\to \infty}\limsup_{n\to \infty} S_{[0,\infty)}(\unJ) &\lesssim \lim_{J\to \infty} \limsup_{n\to \infty} \Bigl(
S_{[0,\infty)}\bigl( \sum_{j=1}^J \vnj\bigr)
        + S_{[0,\infty)}\bigl(e^{it\Delta}\wnJ\bigr)\Bigr)\\
&\lesssim \lim_{J\to \infty} \limsup_{n\to \infty} \sum_{j=1}^J S_{[0,\infty)}(\vnj) \lesssim 1 + E_c.
\end{aligned}
\end{equation}
By the same argument as that used to derive \eqref{S1 ass} from \eqref{S ass}, we obtain
\begin{align}\label{S1 unj}
\lim_{J\to \infty}\limsup_{n\to \infty} \|\unJ\|_{\dot S^1([0, \infty))} \leq C(E_c)<\infty.
\end{align}

To apply Lemma~\ref{L:stab}, it suffices to show that $\unJ$ asymptotically solves \eqref{nls} in the sense that
\begin{align*}
\lim_{J\to\infty} \limsup_{n \to \infty}
    \bigl\| \nabla \bigl[(i \partial_t + \Delta) \unJ - F(\unJ) \bigr]\bigr\|_{L_{t,x}^{\frac{2(d+2)}{d+4}}([0,\infty)\times\R^d)} =0,
\end{align*}
which by the triangle inequality reduces to proving
\begin{align}\label{eq good approx 1}
\lim_{J\to\infty} \limsup_{n \to \infty}
    \Bigl\| \nabla \Bigl[\sum_{j=1}^J F(\vnj) - F\bigl(\sum_{j=1}^J \vnj\bigr) \Bigr]\Bigr\|_{L_{t,x}^{\frac{2(d+2)}{d+4}}([0,\infty)\times\R^d)} =0
\end{align}
and
\begin{align}\label{eq good approx 2}
\lim_{J\to\infty} \limsup_{n \to \infty}
    \bigl\| \nabla \bigl[F\bigl( \unJ-e^{it\Delta} \wnJ\bigr) - F(\unJ)\bigr] \bigr\|_{L_{t,x}^{\frac{2(d+2)}{d+4}}([0,\infty)\times\R^d)} =0.
\end{align}
The arguments we use to prove \eqref{eq good approx 1} and \eqref{eq good approx 2} owe much to the proof of
\cite[Proposition 3.4]{keraani-h1}, particularly to Keraani's treatment of the most delicate point, \eqref{to prove 2}.
We are grateful to C.~Kenig for drawing our attention to this aspect of Keraani's work.

We first address \eqref{eq good approx 1}.  Note that we can write
$$
\Bigl|\nabla \Bigl(\sum_{j=1}^J F(f_j) - F\bigl(\sum_{j=1}^J f_j\bigr) \Bigr) \Bigr| \lesssim_J \sum_{j\neq j'} |\nabla f_j|
|f_{j'}|^{\frac{4}{d-2}}.
$$
Next, recall that by \eqref{tail} and \eqref{S1 ass}, $\vnj\in \dot S^1([0,\infty))$ for all $j\geq 1$ and all $n$ sufficiently
large; invoking \eqref{crazy}, a simple computation shows
$$
\limsup_{n\to \infty} \bigl\| |v_n^{j'}|^{\frac{4}{d-2}} \nabla \vnj \bigl\|_{L_{t,x}^{\frac{2(d+2)}{d+4}}([0,\infty)\times\R^d)} =0
$$
for any $j\neq j'$; see \cite{keraani-h1}.  Thus,
\begin{align*}
\limsup_{n \to \infty} \Bigl\|\nabla \Bigl[ \sum_{j=1}^J F(\vnj) & - F\bigl(\sum_{j=1}^J \vnj\bigr)\Bigr] \Bigr\|_{L_{t,x}^{\frac{2(d+2)}{d+4}}([0,\infty)\times\R^d)}\\
&\lesssim_J \limsup_{n \to \infty}\sum_{j\neq j'} \bigl\|\nabla \vnj |v_n^{j'}|^{\frac{4}{d-2}} \bigl\|_{L_{t,x}^{\frac{2(d+2)}{d+4}}([0,\infty)\times\R^d)} =0
\end{align*}
and \eqref{eq good approx 1} follows.

We now consider \eqref{eq good approx 2}.  In what follows, all spacetime norms are taken on the slab $[0, \infty)\times\R^d$, unless noted
otherwise.  In dimensions $d\geq 6$,
\begin{align*}
\bigl\| \nabla \bigl[F\bigl( \unJ-e^{it\Delta} \wnJ\bigr)  -  F(\unJ) \bigr]\bigr\|_{L_{t,x}^{\frac{2(d+2)}{d+4}}}
&\lesssim \|\nabla e^{it\Delta} \wnJ\|_{L_{t,x}^{\frac{2(d+2)}d}}
        \|e^{it\Delta} \wnJ\|_{L_{t,x}^{\frac{2(d+2)}{d-2}}}^{\frac4{d-2}}\\
&\quad +\|\nabla \unJ\|_{L_{t,x}^{\frac{2(d+2)}d}}
        \|e^{it\Delta} \wnJ\|_{L_{t,x}^{\frac{2(d+2)}{d-2}}}^{\frac4{d-2}}\\
&\quad + \bigl\| |\unJ|^{\frac4{d-2}} \nabla e^{it\Delta} \wnJ \bigr\|_{L_{t,x}^{\frac{2(d+2)}{d+4}}}
\end{align*}
by H\"older's inequality.  In dimensions $d=3,4,5$, one must add the term
$$
\|\nabla \unJ\|_{L_{t,x}^{\frac{2(d+2)}d}} \|e^{it\Delta} \wnJ\|_{L_{t,x}^{\frac{2(d+2)}{d-2}}}
\|\unJ\|_{L_{t,x}^{\frac{2(d+2)}{d-2}}}^{\frac{6-d}{d-2}}
$$
to the right-hand side above.
Using \eqref{w scat}, \eqref{S unj}, \eqref{S1 unj}, and the Strichartz inequality combined with the fact that $\wnJ$
is bounded in $\dot H^1_x$, we see that the claim \eqref{eq good approx 2} follows once we establish
\begin{align}\label{to prove}
\lim_{J\to\infty} \limsup_{n \to \infty}
        \bigl\| |\unJ|^{\frac4{d-2}} \nabla e^{it\Delta} \wnJ \bigr\|_{L_{t,x}^{\frac{2(d+2)}{d+4}}([0, \infty)\times\R^d)}=0.
\end{align}
By H\"older, \eqref{S unj}, and the Strichartz inequality,
\begin{align*}
\bigl\|  |\unJ|^{\frac4{d-2}}& \nabla e^{it\Delta} \wnJ \bigr\|_{L_{t,x}^{\frac{2(d+2)}{d+4}}}\\
&\lesssim \| \unJ \|_{L_{t,x}^{\frac{2(d+2)}{d-2}}}^{\frac3{d-2}}
    \|\nabla e^{it\Delta} \wnJ \bigr\|_{L_{t,x}^{\frac{2(d+2)}d}}^{\frac{d-3}{d-2}}
    \| \unJ \nabla e^{it\Delta} \wnJ \bigr\|_{L_{t,x}^{\frac{d+2}{d-1}}}^{\frac1{d-2}}\\
&\lesssim \bigl\| \bigl(\sum_{j=1}^J \vnj\bigr) \nabla e^{it\Delta} \wnJ  \bigr\|_{L_{t,x}^{\frac{d+2}{d-1}}}^{\frac1{d-2}}
        + \|e^{it\Delta} \wnJ\|_{L_{t,x}^{\frac{2(d+2)}{d-2}}}^{\frac1{d-2}}
            \|\nabla e^{it\Delta} \wnJ\|_{L_{t,x}^{\frac{2(d+2)}d} }^{\frac1{d-2}}\\
&\lesssim \bigl\| \bigl(\sum_{j=1}^J \vnj \bigr) \nabla e^{it\Delta} \wnJ \bigr\|_{L_{t,x}^{\frac{d+2}{d-1}}}^{\frac1{d-2}}
        + \|e^{it\Delta} \wnJ\|_{L_{t,x}^{\frac{2(d+2)}{d-2}}}^{\frac1{d-2}}.
\end{align*}
Invoking \eqref{w scat}, proving \eqref{to prove} reduces to proving
\begin{align}\label{to prove 1}
\lim_{J\to\infty} \limsup_{n \to \infty}
        \bigl\| \bigl(\sum_{j=1}^J \vnj\bigr) \nabla e^{it\Delta} \wnJ  \bigr\|_{L_{t,x}^{\frac{d+2}{d-1}}([0, \infty)\times\R^d)}=0.
\end{align}
Let $\eta>0$.  By \eqref{S decouple}, we see that there exists $J'=J'(\eta)\geq 1$ such that
$$
\sum_{j\geq J'} S_{[0, \infty)}(\vnj) \leq \eta.
$$
Thus, using H\"older's inequality and arguing as for \eqref{S unj},
\begin{align*}
\limsup_{n\to \infty} \bigl\| \bigl(\sum_{j=J'}^J \vnj\bigr) \nabla e^{it\Delta} \wnJ  \bigr\|_{L_{t,x}^{\frac{d+2}{d-1}}}^{\frac{2(d+2)}{d-2}}
&\lesssim \limsup_{n\to \infty} \Bigl( \sum_{j\geq J'} S_{[0, \infty)}(\vnj) \Bigr) \|\nabla e^{it\Delta} \wnJ\|_{L_{t,x}^{\frac{2(d+2)}d} }^{\frac{2(d+2)}{d-2}} \\
&\lesssim \eta.
\end{align*}
As $\eta>0$ is arbitrary, proving \eqref{to prove 1} reduces to showing
\begin{align}\label{to prove 2}
\lim_{J\to\infty} \limsup_{n \to \infty}\| \vnj \nabla e^{it\Delta} \wnJ  \|_{L_{t,x}^{\frac{d+2}{d-1}}([0, \infty)\times\R^d)}=0
    \quad \text{for}\quad 1\leq j\leq J'.
\end{align}
Fix $1\leq j\leq J'$.  By a change of variables,
\begin{align*}
\| \vnj \nabla e^{it\Delta} \wnJ  \|_{L_{t,x}^{\frac{d+2}{d-1}}}
=\bigl\| v^j  \nabla {\twnJ} \bigr\|_{L_{t,x}^{\frac{d+2}{d-1}}},
\end{align*}
where $\twnJ:=\bigl[T_{(\gnj)^{-1}} \bigl(e^{it\Delta} \wnJ\bigr)\bigr](\cdot-\tnj) $.  Note that
\begin{align}\label{tilde w}
S_{\R}(\twnJ) = S_{\R}(e^{it\Delta} \wnJ) \quad \text{and} \quad
        \|\nabla\twnJ\|_{L_{t,x}^{\frac{2(d+2)}d}} = \|\nabla e^{it\Delta} \wnJ\|_{L_{t,x}^{\frac{2(d+2)}d}}.
\end{align}
By density, we may assume $v^j\in C^\infty_c(\R\times\R^d)$.  Invoking H\"older's inequality, it thus suffices to show
\begin{align}\label{to prove 2a}
\lim_{J\to\infty} \limsup_{n \to \infty}\| \nabla \twnJ \|_{L_{t,x}^2(K)}=0
\end{align}
for any compact $K\in \R\times\R^d$.  This is a consequence of Lemma~\ref{L:Keraani3.7}, \eqref{tilde w}, and \eqref{w scat}.
Tracing back through the argument we see that we have verified \eqref{to prove 1} and hence \eqref{eq good approx 2}.

We are now in a position to apply Lemma~\ref{L:stab}; invoking \eqref{S unj}, we conclude that for $n$ sufficiently large,
\begin{align}\label{S un}
S_{[0,\infty)}(u_n) \lesssim 1 + E_c,
\end{align}
thus contradicting \eqref{blow up in two}.  This finishes the proof of Lemma~\ref{L:bad profile}.
\end{proof}

Returning to the proof of Proposition~\ref{P:palais-smale} and rearranging the indices, we may assume that there exists $1\leq
J_1 <J_0$ such that
$$
\limsup_{n\to \infty}S_{[0, T^+_{n,j})}(\vnj)=\infty \text{ for } 1\leq j\leq J_1 \ \ \text{and} \ \ %
        \limsup_{n\to \infty} S_{[0,\infty)}(\vnj)<\infty \text{ for } j> J_1.
$$
Passing to a subsequence in $n$, we can guarantee that $S_{[0, T^+_{n,1})}(v_n^1)\to\infty$.

For each $m,n\geq 1$ let us define an integer $j(m,n)\in \{1, \ldots, J_1\}$ and an interval $K^m_n$ of the form $[0, \tau]$ by
\begin{equation}\label{Kdefn}
\sup_{1\leq j \leq J_1} S_{K^m_n}(\vnj) = S_{K^m_n}(v_n^{j(m,n)}) = m.
\end{equation}
By the pigeonhole principle, there is a $1\leq j_1\leq J_1$ so that for infinitely many $m$ one has $j(m,n)=j_1$ for infinitely
many $n$. Note that the infinite set of $n$ for which this holds may be $m$-dependent.  By reordering the indices, we may assume
that $j_1=1$. Then, by the definition of the critical kinetic energy, we obtain
\begin{align}\label{baddie}
\limsup_{m \to \infty} \; \limsup_{n\to \infty}  \; \sup_{t\in K_n^m} \|\nabla v_n^1(t)\|_2^2 \geq E_c.
\end{align}

On the other hand, by virtue of \eqref{Kdefn}, all $\vnj$ have finite scattering size on $K^m_n$ for each $m\geq 1$. Thus, by
the same argument used in Lemma~\ref{L:bad profile}, we see that for $n$ and $J$ sufficiently large, $\unJ$ is a good
approximation to $u_n$ on each $K_n^m$.  More precisely,
\begin{align}\label{good approx ke}
\lim_{J\to \infty} \limsup_{n\to \infty} \| \unJ - u_n\|_{L_t^\infty \dot H^1_x (K_n^m\times\R^d)}=0
\end{align}
for each $m\geq 1$.

Our next result proves asymptotic kinetic energy decoupling for $\unJ$ at all times of existence.

\begin{lemma}[Kinetic energy decoupling for $\unJ$]\label{L:decouple ke}
For all $J\geq 1$ and $m\geq1$,
\begin{align*}
\limsup_{n\to \infty} \sup_{t\in K^m_n} \Bigl|\|\nabla \unJ(t)\|^2_2 -\sum_{j=1}^J \|\nabla \vnj(t)\|_2^2 - \|\nabla \wnJ\|_2^2
\Bigr|=0.
\end{align*}
\end{lemma}

\begin{proof}
Fix $J\geq 1$ and $m\geq1$. Then, for all $t\in K^m_n$,
\begin{align*}
\|\nabla \unJ(t)\|_2^2
&= \langle \nabla \unJ(t), \nabla \unJ(t)\rangle\\
&= \sum_{j=1}^J \|\nabla \vnj(t)\|_2^2 + \|\nabla \wnJ\|_2^2 + \sum_{j\neq j'} \langle \nabla \vnj(t), \nabla v_n^{j'} (t)\rangle\\
&\quad + \sum_{j=1}^J \bigl( \bigl\langle \nabla e^{it\Delta}\wnJ , \nabla \vnj(t) \bigr\rangle
        + \bigl\langle \nabla \vnj(t), \nabla e^{it\Delta}\wnJ \bigr\rangle \bigr).
\end{align*}

To prove Lemma~\ref{L:decouple ke}, it thus suffices to show that for all sequences $t_n\in K_n^m$,
\begin{align}\label{orthog'}
 \langle \nabla \vnj(t_n), \nabla v_n^{j'} (t_n)\rangle \to 0 \quad \text{as } n\to \infty
\end{align}
and
\begin{align}\label{orthog}
\bigl\langle \nabla e^{it_n\Delta}\wnJ , \nabla \vnj(t_n) \bigr\rangle \to 0 \quad \text{as } n\to \infty
\end{align}
for all $1\leq j,j'\leq J$ with $j\neq j'$.  We will only demonstrate the latter, which requires Lemma~\ref{L:strong decouple};
the former can be deduced in much the same manner using \eqref{crazy}.

By a change of variables,
\begin{align}\label{rescale}
\bigl\langle \nabla e^{it_n\Delta}\wnJ , \nabla \vnj(t_n) \bigr\rangle = \bigl\langle \nabla e^{i t_n (\lambda_n^j)^{-2}\Delta}[
(\gnj)^{-1}\wnJ] , \nabla v^j\bigl(\tfrac{t_n}{(\lambda_n^j)^2} +\tnj \bigr) \bigr\rangle.
\end{align}
As $t_n\in K_n^m \subset [0, T_{n,j}^+)$ for all $1\leq j\leq J_1$, we have $t_n(\lambda_n^j)^{-2} +\tnj \in I^j$ for all $j\geq
1$. Recall that $I^j$ is the maximal lifespan of $v^j$; for $j>J_1$ this is $\R$.  By refining the sequence once for every $j$
and using the standard diagonalisation argument, we may assume $t_n(\lambda_n^j)^{-2} +\tnj$ converges for every $j$.

Fix $1\leq j\leq J$.  If $t_n(\lambda_n^j)^{-2} +\tnj$ converges to some point $\tau^j$ in the interior of $I^j$, then by the
continuity of the flow, $v^j\bigl(t_n(\lambda_n^j)^{-2} +\tnj \bigr)$ converges to $v^j(\tau^j)$ in $\dot H^1_x(\R^d)$. On the
other hand, by \eqref{decouple},
\begin{align}\label{bounded}
\limsup_{n\to \infty} \bigl\| e^{it_n(\lambda_n^j)^{-2}\Delta}[ (\gnj)^{-1}\wnJ]\bigr\|_{\dot H^1_x(\R^d)} =\limsup_{n\to
\infty} \|\wnJ\|_{\dot H^1_x(\R^d)} \lesssim E_c.
\end{align}
Combining this with \eqref{rescale}, we obtain
\begin{align*}
\lim_{n \to \infty}\bigl\langle \nabla e^{it_n\Delta}\wnJ , \nabla \vnj(t_n) \bigr\rangle
&= \lim_{n \to \infty}\bigl\langle \nabla e^{it_n(\lambda_n^j)^{-2}\Delta}[ (\gnj)^{-1}\wnJ] , \nabla v^j (\tau^j) \bigr\rangle\\
&= \lim_{n \to \infty}\bigl\langle \nabla e^{-i\tnj\Delta}[ (\gnj)^{-1}\wnJ] , \nabla e^{-i \tau^j\Delta} v^j (\tau^j)
\bigr\rangle.
\end{align*}
Invoking Lemma~\ref{L:strong decouple}, we deduce \eqref{orthog}.

Consider now the case when $t_n(\lambda_n^j)^{-2} +\tnj$ converges to $\sup I^j$.  Then we must have $\sup I^j=\infty$ and $v^j$
scatters forward in time.  This is clearly true if $\tnj\to \infty$ as $n\to \infty$; in the other cases, failure would imply
$$
\limsup_{n\to \infty} S_{[0,t_n]}(\vnj) = \limsup_{n\to \infty} S_{\bigl[\tnj,t_n(\lambda_n^j)^{-2} + \tnj\bigr]}(v^j) =\infty,
$$
which contradicts $t_n\in K_n^m$.  Therefore, there exists $\psi^j\in \dot H^1_x(\R^d)$ such that
$$
\lim_{n\to \infty} \Bigl\| v^j \bigl(t_n(\lambda_n^j)^{-2} +\tnj\bigr)
            - e^{i\bigl(t_n(\lambda_n^j)^{-2} +\tnj\bigr)\Delta}\psi^j \Bigr\|_{\dot H^1_x(\R^d)} = 0.
$$
Together with \eqref{rescale}, this yields
$$
\lim_{n \to \infty}\bigl\langle \nabla e^{it_n\Delta}\wnJ , \nabla \vnj(t_n) \bigr\rangle = \lim_{n \to \infty}\bigl\langle
\nabla e^{-i\tnj\Delta}[ (\gnj)^{-1}\wnJ] , \nabla \psi^j \bigr\rangle,
$$
which by Lemma~\ref{L:strong decouple} implies \eqref{orthog}.

Finally, we consider the case when $t_n(\lambda_n^j)^{-2} +\tnj$ converges to $\inf I^j$.  Since $t_n(\lambda_n^j)^{-2}\geq 0$
and $\inf I^j<\infty$ for all $j\geq 1$ we see that $\tnj$ does not converge to $+\infty$.  Moreover, if $\tnj\equiv 0$, then
$\inf I^j<0$; as $t_n(\lambda_n^j)^{-2}\geq 0$, we see that $\tnj$ cannot be identically zero.  This leaves $\tnj\to -\infty$ as
$n\to \infty$. Thus $\inf I^j=-\infty$ and $v^j$ scatters backward in time to $e^{it\Delta}\phi^j$.  We obtain
$$
\lim_{n\to \infty} \Bigl\| v^j \bigl(t_n(\lambda_n^j)^{-2} +\tnj\bigr)
            - e^{i\bigl(t_n(\lambda_n^j)^{-2} +\tnj\bigr)\Delta}\phi^j \Bigr\|_{\dot H^1_x(\R^d)} = 0,
$$
which by \eqref{rescale} implies
$$
\lim_{n \to \infty}\bigl\langle \nabla e^{it_n\Delta}\wnJ , \nabla \vnj(t_n) \bigr\rangle = \lim_{n \to \infty}\bigl\langle
\nabla e^{-i\tnj\Delta}[ (\gnj)^{-1}\wnJ] , \nabla \phi^j \bigr\rangle.
$$
Invoking Lemma~\ref{L:strong decouple} once again, we derive \eqref{orthog}.

This finishes the proof of Lemma~\ref{L:decouple ke}.
\end{proof}

We now return to the proof of Proposition~\ref{P:palais-smale}.  By \eqref{max ke}, \eqref{good approx ke}, and
Lemma~\ref{L:decouple ke},
\begin{align*}
E_c\geq \limsup_{n\to\infty} \sup_{t\in K^m_n} \|\nabla u_n(t)\|_2^2
    = \lim_{J\to\infty} \limsup_{n\to\infty} \, \Bigl\{ \|\nabla \wnJ\|_2^2 + \sup_{t\in K^m_n} \sum_{j=1}^J \|\nabla \vnj(t)\|_2^2 \Bigr\}.
\end{align*}
Invoking \eqref{baddie}, this implies $J_1=1$, $\vnj\equiv0$ for all $j\geq 2$, and $w_n:=w_n^1$ converges to zero strongly in
$\dot H^1_x(\R^d)$.  In other words,
\begin{align}\label{just one}
u_n(0)=g_n e^{i\tau_n\Delta}\phi + w_n
\end{align}
for some $g_n\in G$, $\tau_n\in \R$, and some functions $\phi, w_n \in \dot H^1_x(\R^d)$ with $w_n\to 0$ strongly in $\dot H^1_x(\R^d)$.
Moreover, the sequence $\tau_n$ obeys $\tau_n\equiv 0$ or $\tau_n\to \pm \infty$.

If $\tau_n\equiv 0$, \eqref{just one} immediately implies that $u_n(0)$ converges modulo symmetries to $\phi$, which proves
Proposition~\ref{P:palais-smale} in this case.

Finally, we will show that this is the only possible case, that is, $\tau_n$ cannot converge to either $\infty$ or $-\infty$. We
argue by contradiction.  Assume that $\tau_n$ converges to $\infty$; the proof in the negative time direction is essentially the
same. By the Strichartz inequality, $S_\R(\propagate\phi)<\infty$; thus we have
$$
\lim_{n\to \infty} S_{\ge 0}\bigl(\propagate e^{i\tau_n\Delta}\phi\bigr)=0.
$$
Since the action of $G$ preserves linear solutions and the scattering size, this implies
$$
\lim_{n\to \infty}S_{\ge 0}\bigl(\propagate g_n e^{i\tau_n\Delta}\phi\bigr)=0.
$$
Combining this with \eqref{just one} and $w_n\to0$ in $\dot H^1_x$, we conclude
$$
\lim_{n\to \infty}S_{\ge 0}\bigl(\propagate u_n(0)\bigr)=0.
$$
An application on Lemma~\ref{L:stab} yields
$$
\lim_{n\to \infty}S_{\ge 0}(u_n)=0,
$$
which contradicts \eqref{blow up in two}.

This completes the proof of Proposition~\ref{P:palais-smale}.
\end{proof}

\subsection{Proof of Theorem~\ref{T:reduct}}
Suppose $d\geq 3$ is such that Conjecture~\ref{conj} failed.  Then the critical kinetic energy $E_c$ must obey $E_c<\|\nabla
W\|^2_2$. By the definition of the critical kinetic energy, we can find a sequence $u_n:I_n \times \R^d \to \C$ of solutions
to \eqref{nls} with $I_n$ compact,
\begin{align}\label{hyp}
\sup_{n\geq 1} \sup_{t\in I_n} \|\nabla u_n(t)\|_2^2 = E_c, \quad \text{and} \quad \lim_{n\to \infty}S_{I_n}(u_n)=\infty.
\end{align}
Let $t_n\in I_n$ be such that $S_{\geq t_n}(u_n)= S_{\leq t_n}(u_n)$.  Then,
\begin{align}\label{hyp 2}
\lim_{n\to \infty} S_{\ge t_n}(u_n)=\lim_{n\to \infty}S_{\le t_n}(u_n)=\infty.
\end{align}
Using the time-translation symmetry of \eqref{nls}, we may take all $t_n=0$.

Applying Proposition~\ref{P:palais-smale} and passing to a subsequence if necessary, we can find $g_n\in G$ and a
function $u_0\in \dot H_x^1$ such that $g_n u_n(0)$ converges to $u_0$ strongly in $\dot H_x^1$.  By
applying the group action $T_{g_n}$ to the solution $u_n$, we may take $g_n$ to all be the identity.  Thus $u_n(0)$ converges
strongly to $u_0$ in $\dot H_x^1$.

Let $u:\ird \to \C$ be the maximal-lifespan solution to \eqref{nls} with initial data $u(0)=u_0$.  As $u_n(0)\to u_0$ in $\dot H_x^1$,
Lemma~\ref{L:stab} shows that $I\subseteq \liminf I_n$ and
$$
\lim_{n\to \infty} \|u_n -u\|_{L_t^\infty \dot H^1_x(K\times\R^d)}=0,\quad\text{for all compact $K\subset I$}.
$$
Thus by \eqref{hyp},
\begin{align}\label{ke<crit}
\sup_{t\in I} \|\nabla u(t)\|_2^2\leq E_c.
\end{align}

Next we prove that $u$ blows up both forward and backward in time.  Indeed, if $u$ does not blow up forward in time, then
$[0,\infty)\subset I$ and $S_{\ge 0}(u)< \infty$.  By Lemma~\ref{L:stab}, this implies $S_{\ge 0}(u_n)<\infty$ for sufficiently
large $n$, which contradicts \eqref{hyp 2}.  A similar argument proves that $u$ blows up backward in time.

Therefore, by our definition of $E_c$,
$$
\sup_{t\in I} \|\nabla u(t)\|_2^2\geq E_c.
$$
Combining this with \eqref{ke<crit}, we obtain
$$
\sup_{t\in I} \|\nabla u(t)\|_2^2 = E_c.
$$

It remains to show that $u$ is almost periodic modulo symmetries.  Consider an arbitrary sequence $\tau_n \in I$. As $u$ blows up
in both time directions
$$
S_{\ge \tau_n}(u)=S_{\le \tau_n}(u)=\infty.
$$
Applying Proposition~\ref{P:palais-smale}, we conclude that $u(\tau_n)$ admits a convergent subsequence in $\dot H^1_x(\R^d)$
modulo symmetries.  Thus the orbit $\{G u(t): \, t\in I\}$ is precompact in $G\backslash \dot H_x^1$. This concludes the
proof of Theorem~\ref{T:reduct}.\qed

%
%
%
%

\section{The enemies}\label{S:enemies}

In this section, we prove Theorem~\ref{T:enemies}.  The argument owes much to \cite[\S4]{KTV};
indeed, readers seeking a fuller treatment of certain details may consult that paper.

Let $v:J\times\R^d\to\C$ denote a minimal kinetic energy blowup solution whose existence (under the hypotheses of
Theorem~\ref{T:enemies}) is guaranteed by Theorem~\ref{T:reduct}.  We denote the symmetry parameters of $v$ by
$N_v(t)$ and $x_v(t)$. We will construct our solution $u$ by taking a subsequential limit of various normalizations of $v$:

\begin{definition}
Given $t_0\in J$, we define the \emph{normalisation} of $v$ at $t_0$ by
\begin{equation}\label{untn}
v^{[t_0]} := T_{g_{0, -x_v(t_0)N_v(t_0),N_v(t_0)}}\bigr( v( \cdot + t_0) \bigr).
\end{equation}
This solution is almost periodic and has symmetry parameters
$$
N_{v^{[t_0]}}(t) = \frac{N_v(t_0+tN_v(t_0)^{-2})}{N_v(t_0)}
\text{ and }
x_{v^{[t_0]}}(t) = N_v(t_0)[x_v(t_0+tN_v(t_0)^{-2})-x_v(t_0)].
$$
\end{definition}

Note that by the definition of almost periodicity, any sequence of $t_n\in J$ admits a subsequence so that
$v^{[t_n]}(0)$ converges in $\dot H^1_x$.  Furthermore, if $u_0$ denotes this limit and $u:I\times\R^d\to\C$ denotes the
maximal-lifespan solution with $u(0)=u_0$, then $u$ is almost periodic modulo symmetries with the same compactness
modulus function as $v$.  Lastly, $v^{[t_n]}\to u$ in $\dot S^1$ (along the subsequence) uniformly on any compact subset of $I$.

Our first goal is to find a soliton from among the normalizations of $v$ if this is at all possible.  To this end,
for any $T \geq 0$, we define the quantity
\begin{equation}\label{cdef}
\osc(T) := \inf_{t_0 \in J} \,\frac{\sup\, \{ N_v(t) : t \in J \text{ and } |t-t_0| \leq T N_v(t_0)^{-2} \}}
    {\inf\, \{ N_v(t) : t \in J \text{ and } |t-t_0| \leq T N_v(t_0)^{-2} \}},
\end{equation}
which measures the least possible oscillation that one can find in $N_v(t)$ on time intervals of normalised duration $T$.

\medskip

{\bf Case 1:} $\lim_{T \to\infty} \osc(T) < \infty$.  Under this hypothesis, we will be able to extract a soliton-like solution.

Choose $t_n$ so that
$$
\limsup_{n\to\infty} \frac{\sup\, \{ N_v(t) : t \in J \text{ and } |t-t_n| \leq n N_v(t_n)^{-2} \}}
    {\inf\, \{ N_v(t) : t \in J \text{ and } |t-t_n| \leq n N_v(t_n)^{-2} \}} <\infty.
$$
Then a few computations reveal that any subsequential limit $u$ of $v^{[t_n]}$ fulfils the requirements to be classed
as a soliton in the sense of Theorem~\ref{T:enemies}.  In particular, $u$ is global because an almost periodic (modulo
symmetries) solution cannot blow up in finite time without its frequency scale function converging to infinity.

\medskip

When $\osc(T)$ is unbounded, we must seek a solution belonging to one of the remaining two scenarios.
To aid in distinguishing between them, we introduce the quantity
$$
a(t_0) := \frac{ N_v(t_0) }{ \sup\,\{ N_v(t) : t \in J \text{ and } t \leq t_0\}  }
    + \frac{ N_v(t_0) }{ \sup\,\{ N_v(t) : t \in J \text{ and } t \geq t_0\}  }
$$
associated to each $t_0 \in J$.  First we treat the case where $a(t_0)$ can be arbitrarily small.  As we will see,
this may lead to either a finite-time blowup solution or to a cascade.

\medskip

{\bf Case 2:} $\lim_{T\to \infty} \osc(T) = \infty$ and $\inf_{t_0 \in J} a(t_0) = 0$.  From the behaviour of $a(t_0)$
we may choose sequences $t_n^{-}<t_n<t_n^{+}$ from $J$ so that $a(t_n)\to 0$, $N_v(t_n^{-})/N_v(t_n)\to\infty$,
and $N_v(t_n^{+})/N_v(t_n)\to\infty$.  Next we choose times $t_n'\in(t_n^{-},t_n^{+})$ so that
\begin{align}\label{N from below}
N_v(t_n') \leq 2 \inf\, \{ N(t) : t\in [t_n^{-},t_n^{+}] \}.
\end{align}
In particular, $N(t_n') \leq 2N(t_n)$, which allows us to deduce that
\begin{align}\label{N to infty}
\frac{ N_v(t_n^{-}) }{ N_v(t_n') }\to\infty \quad\text{and}\quad \frac{ N_v(t_n^{+}) }{ N_v(t_n') } \to\infty.
\end{align}

Let $u$ denote a subsequential limit of $v^{[t_n']}$ and let $I$ denote its maximal lifespan.  If $I$ has a finite endpoint,
then $u$ is a finite-time blowup solution in the sense of Theorem~\ref{T:enemies} and we are done.  Thus we are left to
consider the case $I=\R$.

Let $s_n^\pm := (t_n^\pm - t_n') N_v(t_n')^2$.  From \eqref{N to infty} we see that $N_u(s_n^\pm)\to\infty$ and so
deduce $s_n^\pm\to\pm\infty$ from the fact that $u$ is a global solution.  Combining this with \eqref{N from below} we
find that $N_u(t)$ is bounded from below uniformly for $t\in\R$.  Rescaling $u$ slightly, we may ensure that $N_u(t)\geq 1$
for all $t\in\R$.

From the fact that $\osc(T)\to\infty$, we see that $N_v(t)$ must show significant oscillation in neighbourhoods of $t_n'$.
Transferring this information to $u$ and using the lower bound on $N_u(t)$ we may conclude that
$\limsup_{|t|\to\infty} N_u(t) =\infty$.  Using time-reversal symmetry, if necessary, we obtain a low-to-high cascade
in the sense of Theorem~\ref{T:enemies}.

\medskip

Finally, we treat the case where $a(t_0)$ is strictly positive; we will construct a finite-time blowup solution.

\medskip

{\bf Case 3:} $\lim_{T\to \infty} \osc(T) = \infty$ and $\inf_{t_0 \in J} a(t_0) = 2\eps >0$.  Let us call a $t_0\in J$
\emph{future-spreading} if $N(t)\leq\eps^{-1} N(t_0)$ for all $t\geq t_0$; we call $t_0$ \emph{past-spreading} if
$N(t)\leq\eps^{-1} N(t_0)$ for all $t\leq t_0$.  Note that by hypothesis, every $t_0\in J$ is future-spreading, past-spreading,
or possibly both.

The fact that even a single time is future- or past-spreading guarantees that $J$ must be infinite in the forward or
reverse time direction, respectively; recall that finite-time blowup is accompanied by $N_v(t)\to\infty$ as $t$ approaches
the blowup time.  Next we argue
that either all sufficiently late times are future-spreading or all sufficiently early times are past-spreading.  If this
were not the case, one would be able to find arbitrarily long time intervals beginning with a future-spreading time and
ending with a past-spreading time.  The existence of such intervals would contradict the divergence of $\osc(T)$.
By appealing to time-reversal symmetry, we restrict our attention to the case where all $t\geq t_0$ are future-spreading.

Choose $T$ so that $\osc(T) > 2\eps^{-1}$.
We will now recursively construct an increasing sequence of times $\{t_n\}_{n=0}^\infty$ so that
\begin{align}\label{t_n props}
0 \leq t_{n+1} - t_n \leq 8T N(t_n)^{-2} \quad\text{ and }\quad
N(t_{n+1}) \leq \tfrac12 N(t_n).
\end{align}
Given $t_n$, set $t_n':=t_n + 4T N(t_n)^{-2}$.  If $N(t_n')\leq \frac12 N(t_n)$ we choose $t_{n+1}=t_n'$ and the properties
set out above follow immediately.  If $N_v(t_n') > \frac12 N_v(t_n)$, then
$$
J_n:=[t_n' - T N_v(t_n')^{-2},t_n' + T N_v(t_n')^{-2}] \subseteq [t_n,t_n+8T N(t_n)^{-2}].
$$
As $t_n$ is future-spreading, this allows us to conclude that $N(t)\leq \eps^{-1} N(t_n)$ on $J_n$, but then by the
way $T$ is chosen, we may find $t_{n+1}\in J_n$ so that $N(t_{n+1})\leq \frac12 N(t_n)$.

Having obtained a sequence of times obeying \eqref{t_n props}, we may conclude that any subsequential limit $u$ of $v^{[t_n]}$
is a finite-time blowup solution.  To elaborate, set $s_n := (t_0-t_n)N(t_n)^2$ and note that
$N_{v^{[t_n]}}(s_n)\geq 2^n$.  However, $s_n$ is a bounded sequence; indeed,
\begin{align*}
|s_n| = N(t_n)^2 \sum^{n-1}_{k=0} \bigl[ t_{k+1} - t_k \bigr] \leq 8 T \sum^{n-1}_{k=0} \frac{N(t_n)^2}{N(t_k)^2}
\leq 8 T \sum^{n-1}_{k=0} 2^{-(n-k)} \leq 8T.
\end{align*}
In this way, we see that the solution $u$ must blow up at some time $-8T\leq t < 0$.

This completes the proof of Theorem~\ref{T:enemies}.
\qed

%
%
%
%

\section{Finite-time blowup}\label{S:finite time}

In this section we preclude scenario I from Theorem~\ref{T:enemies}.  In this particular case, we do not need to restrict
to dimensions $d\geq 5$.  The argument is essentially taken from \cite{kenig-merle}.

\begin{theorem}[No finite-time blowup]\label{T:no ftb}
Let $d\geq 3$.  Then there are no maximal-lifespan solutions $u:I\times\R^d \to \C$ to \eqref{nls} that are almost periodic
modulo symmetries, obey
\begin{align}\label{infinite norm}
S_I(u)=\infty,
\end{align}
and
\begin{equation}\label{hype}
\sup_{t\in I} \|\nabla u(t)\|_2 < \|\nabla W\|_2,
\end{equation}
and are such that either $|\inf I|<\infty$ or $\sup I <\infty$.
\end{theorem}

\begin{proof}
Suppose for a contradiction that there existed such a solution $u$.  Without loss of generality, we may assume $\sup I<\infty$.
We first argue that
\begin{align}\label{volcano}
\liminf_{t\nearrow \sup I} N(t) =\infty.
\end{align}

Assume for contradiction that $\liminf_{t\nearrow \sup I} N(t) <\infty$.  Let $t_n \in  I$ such that $t_n\nearrow \sup I$, and
define the rescaled functions $v_n: I_n\times\R^d\to \C$ by
$$
v_n(t,x):= u^{[t_n]}(t,x) = N(t_n)^{-\frac{d-2}2} u\bigl( t_n + t N(t_n)^{-2}, x(t_n) + x N(t_n)^{-1} \bigr),
$$
where $0\in I_n:=\{ t_n + N(t_n)^{-2}t: \, t\in I\}$.  Then each $v_n$ is a solution to \eqref{nls} and $\{v_n(0)\}_n$ is
precompact in $\dot H^1_x(\R^d)$.  Thus, after passing to a subsequence if necessary, we may assume that $v_n(0)$ converges
strongly in $\dot H^1_x(\R^d)$ to some function $v_0$.  As $\|\nabla v_n(0)\|_2 = \|\nabla u(t_n)\|_2$
and, by assumption, $u$ is not identically zero, we conclude (using Sobolev embedding and the conservation of energy) that $v_0$
is not identically zero.

Let $v$ be the solution to \eqref{nls} with initial data $v_0$ at time $t=0$ and maximal lifespan $(-T_-, T_+)$ with $-\infty
\leq -T_- < 0 < T_+ \leq \infty$.  From the local theory for \eqref{nls} (see, for example, Lemma~\ref{L:stab}), $v_n$ is
well-posed and has finite scattering size on any compact interval $J\in (-T_-, T_+)$.  In particular, $u$ is well-posed with
finite scattering size on $\{ t_n + N(t_n)^{-2}t: \, t\in J\}$.  However, as $t_n\nearrow \sup I$ and
$\liminf_{n\to \infty} N(t_n)<\infty$, this means that $u$ has finite scattering size beyond $\sup I$, which
contradicts the fact that, by assumption, $u$ blows up forward in time on $I$.  Thus \eqref{volcano} must hold.

We now show that \eqref{volcano} implies
\begin{align}\label{mass leaves balls}
\limsup_{t\nearrow \,\sup I} \int_{|x|\leq R} |u(t,x)|^2\, dx = 0 \quad \text{for all $R>0$.}
\end{align}
Indeed, let $0<\eta<1$ and $t\in I$.  By H\"older's inequality, Sobolev embedding, and \eqref{hype},
\begin{align*}
\int_{|x|\leq R} |u(t,x)|^2\, dx
&\leq  \int_{|x-x(t)|\leq \eta R} |u(t,x)|^2\, dx + \mathop{\int_{ |x|\leq R}}_{|x-x(t)|>\eta R} |u(t,x)|^2\, dx\\
&\lesssim \eta^2 R^2 \|u(t)\|_{\frac{2d}{d-2}}^2 + R^2 \Bigl(  \int_{|x-x(t)|>\eta R} |u(t,x)|^{\frac{2d}{d-2}} \, dx \Bigr)^{\frac{d-2}d} \\
&\lesssim \eta^2 R^2 \|\nabla W\|_2^2 + R^2 \Bigl(  \int_{|x-x(t)|>\eta R} |u(t,x)|^{\frac{2d}{d-2}} \, dx \Bigr)^{\frac{d-2}d}.
\end{align*}
Letting $\eta\to 0$, we can make the first term on the right-hand side of the inequality above as small as we wish.
On the other hand, by \eqref{volcano}, almost periodicity modulo symmetries, and Remark~\ref{R:pot energy}, we see that
$$
\limsup_{t\nearrow \sup I} \int_{|x-x(t)|>\eta R} |u(t,x)|^{\frac{2d}{d-2}} \, dx =0.
$$
This proves \eqref{mass leaves balls}.

The next step is to prove that \eqref{mass leaves balls} implies the solution $u$ is identically zero, thus contradicting \eqref{infinite norm}.
For $t\in I$ define
$$
M_R(t):=\int_{\R^d} \phi\bigl(\tfrac{|x|}R\bigr) |u(x,t)|^2\,dx,
$$
where $\phi$ is a smooth, radial function, such that
\begin{align*}
\phi(r)=\begin{cases}
1 & \text{for } r\leq 1\\
0 & \text{for } r\geq 2.
\end{cases}
\end{align*}
By \eqref{mass leaves balls},
\begin{align}\label{M_R -> 0}
\limsup_{t\nearrow \sup I} M_R(t)=0 \quad \text{for all $R>0$.}
\end{align}
On the other hand, a simple computation involving Hardy's inequality and \eqref{hype} shows
$$
|\partial_t M_R(t)|\lesssim \|\nabla u(t)\|_2 \Bigl\|\frac{u(t)}{|x|}\Bigr\|_2\lesssim \|\nabla u(t)\|_2^2 \lesssim \|\nabla
W\|_2^2.
$$
Thus, by the Fundamental Theorem of Calculus,
\begin{align*}
M_R(t_1)= M_R(t_2) + \int_{t_2}^{t_1} \partial_t M_R(t)\, dt \lesssim M_R(t_2) + |t_1-t_2| \|\nabla W\|_2^2
\end{align*}
for all $t_1, t_2 \in I$ and $R>0$.  Letting $t_2\nearrow \sup I$ and invoking \eqref{M_R -> 0}, we deduce
$$
M_R(t_1) \lesssim |\sup I-t_1| \|\nabla W\|_2^2.
$$
Now letting $R\to \infty$ and using the conservation of mass, we obtain $u_0\in L_x^2(\R^d)$.  Finally, letting $t_1\nearrow
\sup I$ we conclude $u_0=0$.  By the uniqueness statement in Theorem~\ref{T:local}, this implies that the solution $u$ is
identically zero, contradicting \eqref{infinite norm}.

This concludes the proof of Theorem~\ref{T:no ftb}.
\end{proof}

%
%
%
%

\section{Negative regularity}\label{S:neg}

In this section we prove

\begin{theorem}[Negative regularity in the global case]\label{T:-reg}
Let $d\geq 5$ and let $u$ be a global solution to \eqref{nls} that is almost periodic modulo symmetries.  Suppose also that
\begin{align}\label{ke bounded}
\sup_{t\in \R} \|\nabla u(t)\|_{L_x^2} <\infty
\end{align}
and
\begin{align}\label{inf bounded}
\inf_{t\in \R} N(t)\geq 1.
\end{align}
Then $u\in L_t^\infty \dot H_x^{-\eps}(\R\times\R^d)$ for some $\eps=\eps(d)>0$.  In particular, $u\in L^\infty_t L^2_x$.
\end{theorem}

The proof of Theorem~\ref{T:-reg} is achieved in two steps: First, we `break' scaling in a Lebesque space; more precisely,
we prove that our solution lives in $L^\infty_t L_x^p$ for some $2<p<\tfrac{2d}{d-2}$.  Next, we use a double Duhamel trick to upgrade this
to $u\in L_t^\infty \dot H_x^{1-s}$ for some $s=s(p,d)>0$.  Iterating the second step finitely many times, we derive Theorem~\ref{T:-reg}.

We learned the double Duhamel trick from \cite{tao:attractor} where it is used for a similar purpose; however, in that paper, the
breach of scaling comes directly from the subcritical nature of the nonlinearity.

Let $u$ be a solution to \eqref{nls} that obeys the hypotheses of Theorem~\ref{T:-reg}.  Let $\eta>0$ be a small constant to be chosen
later.  Then by Remark~\ref{R:c small} combined with \eqref{inf bounded}, there exists $N_0=N_0(\eta)$ such that
\begin{align}\label{ke small}
\|\nabla u_{\leq N_0}\|_{L_t^\infty L_x^2 (\R\times\R^d)}\leq \eta.
\end{align}

We turn now to our first step, that is, breaking scaling in a Lebesgue space.  To this end, we define
\begin{equation}
A(N):=
\begin{cases}
N^{-\frac2{d-2}} \sup_{t\in \R} \|u_N(t)\|_{L_x^{\frac{2(d-2)}{d-4}}} & \quad \text{for}\quad d\geq 6\\
N^{-\frac 12} \sup_{t\in \R} \|u_N(t)\|_{L_x^5} & \quad \text{for}\quad d=5.
\end{cases}
\end{equation}
for frequencies $N\leq 10N_0$.  Note that by Bernstein's inequality combined with Sobolev embedding and \eqref{ke bounded},
$$
A(N)\lesssim \|u_N\|_{L_t^\infty L_x^{\frac{2d}{d-2}}} \lesssim \|\nabla u\|_{L_t^\infty L_x^2} <\infty.
$$

We next prove a recurrence formula for $A(N$).

\begin{lemma}[Recurrence]\label{L:recurrence}
For all $N\leq 10 N_0$,
$$
A(N)\lesssim_u \bigl(\tfrac{N}{N_0}\bigr)^{\alpha}
  + \eta^{\frac4{d-2}} \sum_{\frac{N}{10}\leq N_1\leq N_0} \bigl(\tfrac{N}{N_1}\bigr)^{\alpha}A(N_1)
  + \eta^{\frac4{d-2}} \sum_{N_1<\frac{N}{10}} \bigl(\tfrac{N_1}{N}\bigr)^{\alpha}A(N_1),
$$
where $\alpha:=\min\{\tfrac2{d-2}, \tfrac 12\}$.
\end{lemma}

\begin{proof}
We first give the proof in dimensions $d\geq 6$.  Once this is completed, we will explain the changes necessary to treat $d=5$.

Fix $N\leq 10 N_0$.  By time-translation symmetry, it suffices to prove
\begin{align}\label{rec goal}
N^{-\frac2{d-2}} \| u_N(0)\|_{L_x^{\frac{2(d-2)}{d-4}}}
&\lesssim_u \bigl(\tfrac{N}{N_0}\bigr)^{\frac2{d-2}}
  + \eta^{\frac4{d-2}} \sum_{\frac{N}{10}\leq N_1\leq N_0} \bigl(\tfrac{N}{N_1}\bigr)^{\frac2{d-2}}A(N_1) \notag\\
&\qquad + \eta^{\frac4{d-2}} \sum_{N_1<\frac{N}{10}} \bigl(\tfrac{N_1}{N}\bigr)^{\frac2{d-2}}A(N_1).
\end{align}

Using the Duhamel formula \eqref{Duhamel} into the future followed by the triangle inequality, Bernstein, and
the dispersive inequality, we estimate
\begin{align}\label{to estimate}
N^{-\frac2{d-2}} \| u_N(0)\|_{L_x^{\frac{2(d-2)}{d-4}}}
&\leq N^{-\frac2{d-2}} \Bigl\| \int_0^{N^{-2}} e^{-it\Delta} P_N F(u(t))\, dt \Bigr\|_{L_x^{\frac{2(d-2)}{d-4}}} \notag\\
&\quad + N^{-\frac2{d-2}}  \int_{N^{-2}}^\infty \bigl\|e^{-it\Delta}  P_N F(u(t))\, dt \bigr\|_{L_x^{\frac{2(d-2)}{d-4}}} \notag \\
&\lesssim N \Bigl\| \int_0^{N^{-2}} e^{-it\Delta}  P_N F(u(t)) \, dt \Bigr\|_{L_x^2}\notag \\
&\quad + N^{-\frac2{d-2}} \| P_N F(u)\|_{L_t^\infty L_x^{\frac{2(d-2)}d}} \int_{N^{-2}}^\infty t^{-\frac d{d-2}}\, dt  \notag \\
&\lesssim N^{-1} \| P_N F(u)\|_{L_t^\infty L_x^2} + N^{\frac2{d-2}} \|  P_N  F(u)\|_{L_t^\infty L_x^{\frac{2(d-2)}d}}\notag\\
&\lesssim N^{\frac2{d-2}} \|  P_N  F(u)\|_{L_t^\infty L_x^{\frac{2(d-2)}d}}.
\end{align}

Using the Fundamental Theorem of Calculus, we decompose
\begin{align}\label{decomposition}
F(u)&= O(|u_{>N_0}| |u_{\leq N_0}|^{\frac4{d-2}}) + O(|u_{>N_0}|^{\frac{d+2}{d-2}}) + F(u_{\frac N{10}\leq \cdot\leq N_0})\notag \\
&\quad + u_{<\frac N{10}} \int_0^1 F_z \bigl( u_{\frac N{10}\leq \cdot\leq N_0} + \theta u_{<\frac N{10}}\bigr)\, d\theta  \\
&\quad + \overline{u_{<\frac N{10}}} \int_0^1 F_{\bar z} \bigl( u_{\frac N{10}\leq \cdot\leq N_0} + \theta u_{<\frac N{10}}\bigr)\, d\theta. \notag
\end{align}

The contribution to the right-hand side of \eqref{to estimate} coming from terms that contain at least one copy of $u_{> N_0}$
can be estimated in the following manner:  Using H\"older, Bernstein, and \eqref{ke bounded},
\begin{align}\label{1}
N^{\frac2{d-2}} \| P_N O(|u_{>N_0}| |u|^{\frac4{d-2}}) \bigr\|_{L_t^\infty L_x^{\frac{2(d-2)}d}}
&\lesssim N^{\frac2{d-2}} \|u_{>N_0}\|_{L_t^\infty L_x^{\frac{2d(d-2)}{d^2-4d+8}}} \|u\|_{L_t^\infty L_x^{\frac{2d}{d-2}}}^{\frac4{d-2}}\notag\\
&\lesssim_u N^{\frac2{d-2}} N_0^{-\frac2{d-2}}.
\end{align}
Thus, this contribution is acceptable.

Next we turn to the contribution to the right-hand side of \eqref{to estimate} coming from the last two terms in \eqref{decomposition};
it suffices to consider the first of them since similar arguments can be used to deal with the second.

First we note that as $\nabla u \in L_t^\infty L_x^2$, we have $F_z(u) \in \dot \Lambda^{\frac{d-2}2,\infty}_{\frac4{d-2}}$.
Furthermore, as $P_{>\frac N{10}}F_z(u)$ is restricted to high frequencies, the Besov characterization of the
homogeneous H\"older continuous functions (see \cite[\S VI.7.8]{stein:large}) yields
$$
\bigl\|P_{>\frac N{10}}F_z(u)\bigr\|_{L_t^\infty L_x^{\frac{d-2}2}}
\lesssim N^{-\frac4{d-2}} \|\nabla u\|_{L_t^\infty L_x^2}^{\frac{4}{d-2}}.
$$
Thus, by H\"older's inequality and \eqref{ke small},
\begin{align}\label{2}
N^{\frac2{d-2}}&\Bigl\| P_N\Bigl( u_{<\frac N{10}}  \int_0^1 F_z \bigl( u_{\frac N{10}\leq \cdot\leq N_0}
            + \theta u_{<\frac N{10}}\bigr)\, d\theta\Bigr) \Bigr\|_{L_t^\infty L_x^{\frac{2(d-2)}d}}\notag\\
&\lesssim N^{\frac2{d-2}}\|u_{<\frac N{10}}\|_{L_t^\infty L_x^{\frac{2(d-2)}{d-4}}}
        \Bigl\| P_{>\frac N{10}}\Bigl(\int_0^1 F_z \bigl( u_{\frac N{10}\leq \cdot\leq N_0}+ \theta u_{<\frac N{10}}\bigl)\, d\theta\Bigr)\Bigr\|_{L_t^\infty L_x^{\frac{d-2}2}}\notag\\
&\lesssim N^{-\frac2{d-2}} \|u_{<\frac N{10}}\|_{L_t^\infty L_x^{\frac{2(d-2)}{d-4}}} \|\nabla u_{\leq N_0}\|_{L_t^\infty L_x^2}^{\frac{4}{d-2}}\notag\\
&\lesssim \eta^{\frac4{d-2}}  \sum_{N_1<\frac{N}{10}} \bigl(\tfrac{N_1}{N}\bigr)^{\frac2{d-2}}A(N_1).
\end{align}
Hence, the contribution coming from the last two terms in \eqref{decomposition} is acceptable.

We are left to estimate the contribution of $F(u_{\frac N{10}\leq \cdot\leq N_0})$ to the right-hand side of \eqref{to estimate}.
We need only show
\begin{align}\label{last goal}
\|F(u_{\frac N{10}\leq \cdot\leq N_0})\|_{L_t^\infty L_x^{\frac{2(d-2)}d}}
\lesssim \eta^{\frac4{d-2}} \sum_{\frac{N}{10}\leq N_1\leq N_0} N_1^{-\frac2{d-2}}A(N_1).
\end{align}
As $d\geq 6$, we have $\tfrac4{d-2}\leq 1$.  Using the triangle inequality, Bernstein,
\eqref{ke small}, and H\"older, we estimate
\begin{align*}
\|F(& u_{\frac N{10}\leq \cdot\leq N_0})\|_{L_t^\infty L_x^{\frac{2(d-2)}d}}\\
&\lesssim \sum_{\frac N{10}\leq N_1 \leq N_0} \bigl\|u_{N_1} |u_{\frac N{10}\leq \cdot\leq N_0}|^{\frac4{d-2}}\bigr\|_{L_t^\infty L_x^{\frac{2(d-2)}d}}\\
&\lesssim \sum_{\frac N{10}\leq N_1, N_2 \leq N_0} \bigl\|u_{N_1} |u_{N_2}|^{\frac4{d-2}}\bigr\|_{L_t^\infty L_x^{\frac{2(d-2)}d}}\\
&\lesssim \sum_{\frac N{10}\leq N_1\leq N_2 \leq N_0} \|u_{N_1}\|_{L_t^\infty L_x^{\frac{2(d-2)}{d-4}}} \|u_{N_2}\|_{L_t^\infty L_x^2}^{\frac4{d-2}}\\
&\qquad + \sum_{\frac N{10}\leq N_2\leq N_1 \leq N_0} \|u_{N_1}\|_{L_t^\infty L_x^2}^{\frac4{d-2}}
            \|u_{N_1}\|_{L_t^\infty L_x^{\frac{2(d-2)}{d-4}}}^{\frac{d-6}{d-2}} \|u_{N_2}\|_{L_t^\infty L_x^{\frac{2(d-2)}{d-4}}}^{\frac4{d-2}}\\
&\lesssim \sum_{\frac N{10}\leq N_1\leq N_2 \leq N_0} \|u_{N_1}\|_{L_t^\infty L_x^{\frac{2(d-2)}{d-4}}} \eta^{\frac4{d-2}}N_2^{-\frac4{d-2}}\\
&\qquad + \sum_{\frac N{10}\leq N_2\leq N_1 \leq N_0}\eta^{\frac4{d-2}}N_1^{-\frac4{d-2}}
            \|u_{N_1}\|_{L_t^\infty L_x^{\frac{2(d-2)}{d-4}}}^{\frac{d-6}{d-2}} \|u_{N_2}\|_{L_t^\infty L_x^{\frac{2(d-2)}{d-4}}}^{\frac4{d-2}}\\
&\lesssim \eta^{\frac4{d-2}} \sum_{\frac N{10}\leq N_1\leq N_0} N_1^{-\frac2{d-2}}A(N_1)\\
&\qquad + \eta^{\frac4{d-2}} \sum_{\frac N{10}\leq N_2\leq N_1 \leq N_0} \bigl(\tfrac{N_2}{N_1} \bigr)^{\frac{16}{(d-2)^2}}
            \bigl(N_1^{-\frac2{d-2}}A(N_1)\bigr)^{\frac{d-6}{d-2}} \bigl(N_2^{-\frac2{d-2}}A(N_2)\bigr)^{\frac{4}{d-2}}\\
&\lesssim \eta^{\frac4{d-2}} \sum_{\frac N{10}\leq N_1\leq N_0} N_1^{-\frac2{d-2}}A(N_1).
\end{align*}
This proves \eqref{last goal} and so completes the proof of the lemma in dimensions $d\geq 6$.

Consider now $d=5$.  Arguing as for \eqref{to estimate}, we have
$$
N^{-\frac 12} \|u_N(0)\|_{L_x^5}\lesssim N^{\frac 12} \|P_N F(u)\|_{L_t^\infty L_x^{\frac54}},
$$
which we estimate by decomposing the nonlinearity as in \eqref{decomposition}.  The analogue of \eqref{1} in this case is
\begin{align*}
N^{\frac12} \| P_N O(|u_{>N_0}| |u|^{\frac4{d-2}}) \bigr\|_{L_t^\infty L_x^{\frac54}}
&\lesssim N^{\frac12} \|u_{>N_0}\|_{L_t^\infty L_x^{\frac52}} \|u\|_{L_t^\infty L_x^{\frac{10}{3}}}^{\frac43}
\lesssim_u N^{\frac12} N_0^{-\frac12}.
\end{align*}
Using Bernstein and Lemma~\ref{F Lip} together with \eqref{ke small}, we replace \eqref{2} by
\begin{align*}
N^{\frac12}&\Bigl\| P_N\Bigl( u_{<\frac N{10}}  \int_0^1 F_z \bigl( u_{\frac N{10}\leq \cdot\leq N_0}
            + \theta u_{<\frac N{10}}\bigr)\, d\theta\Bigr) \Bigr\|_{L_t^\infty L_x^{\frac54}}\\
&\lesssim N^{\frac12}\|u_{<\frac N{10}}\|_{L_t^\infty L_x^5}
        \Bigl\| P_{>\frac N{10}}\Bigl(\int_0^1 F_z \bigl( u_{\frac N{10}\leq \cdot\leq N_0}+ \theta u_{<\frac N{10}}\bigl)\, d\theta\Bigr)\Bigr\|_{L_t^\infty L_x^{\frac53}}\\
&\lesssim N^{-\frac12} \|u_{<\frac N{10}}\|_{L_t^\infty L_x^5} \|\nabla u_{\leq N_0}\|_{L_t^\infty L_x^2} \|u_{\leq N_0}\|_{L_t^\infty L_x^{\frac{10}3}}^{\frac13}\\
&\lesssim \eta^{\frac43}  \sum_{N_1<\frac{N}{10}} \bigl(\tfrac{N_1}{N}\bigr)^{\frac12}A(N_1).
\end{align*}
Finally, arguing as for \eqref{last goal}, we estimate
\begin{align*}
\|F&( u_{\frac N{10}\leq \cdot\leq N_0})\|_{L_t^\infty L_x^{\frac54}}\\
&\lesssim \sum_{\frac N{10}\leq N_1, N_2 \leq N_0} \bigl\|u_{N_1} u_{N_2}|u_{\frac N{10}\leq \cdot\leq N_0}|^{\frac13}\bigr\|_{L_t^\infty L_x^{\frac54}}\\
&\lesssim \sum_{\frac N{10}\leq N_1\leq N_2, N_3 \leq N_0} \|u_{N_1}\|_{L_t^\infty L_x^5}\|u_{N_2}\|_{L_t^\infty L_x^{\frac{20}9}}\|u_{N_3}\|_{L_t^\infty L_x^{\frac{20}9}}^{\frac13}\\
&\qquad + \sum_{\frac N{10}\leq N_3\leq N_1\leq N_2 \leq N_0} \|u_{N_1}\|_{L_t^\infty L_x^5}^{\frac23}\|u_{N_1}\|_{L_t^\infty L_x^{\frac{20}9}}^{\frac13}
        \|u_{N_2}\|_{L_t^\infty L_x^{\frac{20}9}}\|u_{N_3}\|_{L_t^\infty L_x^5}^{\frac13}\\
&\lesssim \sum_{\frac N{10}\leq N_1\leq N_2, N_3 \leq N_0} \|u_{N_1}\|_{L_t^\infty L_x^5} \eta N_2^{-\frac 34} \eta^{\frac13} N_3^{-\frac14}\\
&\qquad + \sum_{\frac N{10}\leq N_3\leq N_1\leq N_2 \leq N_0} \|u_{N_1}\|_{L_t^\infty L_x^5}^{\frac23} \eta^{\frac13} N_1^{-\frac14}
        \eta N_2^{-\frac34} \|u_{N_3}\|_{L_t^\infty L_x^5}^{\frac13}\\
&\lesssim \eta^{\frac43} \sum_{\frac N{10}\leq N_1 \leq N_0} N_1^{-\frac12} A(N_1)\\
&\qquad + \eta^{\frac43} \sum_{\frac N{10}\leq N_3\leq N_1 \leq N_0} \bigl( \tfrac{N_3}{N_1}\bigr)^{\frac13}
            \bigl(N_1^{-\frac12}A(N_1)\bigr)^{\frac23} \bigl(N_3^{-\frac12}A(N_3)\bigr)^{\frac13}\\
&\lesssim \eta^{\frac43} \sum_{\frac N{10}\leq N_1 \leq N_0} N_1^{-\frac12} A(N_1).
\end{align*}
Putting everything together completes the proof of the lemma in the case $d=5$.
\end{proof}

This lemma leads very quickly to our first goal:

\begin{proposition}[$L^p$ breach of scaling]\label{P:L^p breach}
Let $u$ be as in Theorem~\ref{T:-reg}.  Then
\begin{align}\label{step 1}
u\in L_t^\infty L_x^p \quad \text{for} \quad \tfrac{2(d+1)}{d-1}\leq p<\tfrac{2d}{d-2}.
\end{align}
In particular, by H\"older's inequality,
\begin{align}\label{breach 2}
\nabla F(u) \in L_t^\infty L_x^r \quad \text{for} \quad \tfrac{2(d-2)(d+1)}{d^2+3d-6} \leq r<\tfrac{2d}{d+4}.
\end{align}
\end{proposition}

\begin{proof}
We only present the details for $d\geq 6$.  The treatment of $d=5$ is completely analogous.

Combining Lemma~\ref{L:recurrence} with Lemma~\ref{L:Gronwall}, we deduce
\begin{align}\label{breach 1}
\|u_N\|_{L_t^\infty L_x^{\frac{2(d-2)}{d-4}}}\lesssim_u N^{\frac4{d-2}-} \quad \text{for all} \quad N\leq 10 N_0.
\end{align}
In applying Lemma~\ref{L:Gronwall}, we set $N=10 \cdot 2^{-k} N_0$, $x_k=A(10 \cdot 2^{-k} N_0)$, and take $\eta$ sufficiently small.

By interpolation followed by \eqref{breach 1}, Bernstein, and \eqref{ke bounded},
\begin{align*}
\|u_N\|_{L_t^\infty L_x^p}
&\leq \|u_N\|_{L_t^\infty L_x^{\frac{2(d-2)}{d-4}}}^{(d-2)(\frac12-\frac1p)} \|u_N\|_{L_t^\infty L_x^2}^{\frac{d-2}p-\frac{d-4}2}\\
&\lesssim_u N^{\frac{2(p-2)}p-} N^{\frac{d-4}2-\frac{d-2}p}\\
&\lesssim_u N^{\frac1{d+1}-}
\end{align*}
for all $N\leq 10 N_0$.  Thus, using Bernstein together with \eqref{ke bounded}, we obtain
\begin{align*}
\|u\|_{L_t^\infty L_x^p}
\leq \|u_{\leq N_0}\|_{L_t^\infty L_x^p} + \|u_{> N_0}\|_{L_t^\infty L_x^p}
\lesssim_u \sum_{N\leq N_0} N^{\frac1{d+1}-} + \sum_{N>N_0} N^{\frac{d-2}2-\frac dp}
\lesssim_{u} 1,
\end{align*}
which completes the proof of the proposition.
\end{proof}

\begin{remark}
With a few modifications, the argument used in dimension five can be adapted to dimensions three and four.
However, $u(t,x)=W(x)$ provides an explicit counterexample to Theorem~\ref{T:-reg} in  these dimensions.
At a technical level, the obstruction is that the strongest dispersive estimate available is $|t|^{-d/2}$,
which is insufficient to perform both integrals in the double Duhamel trick for $d\leq 4$.
\end{remark}

Our second step is to use the double Duhamel trick to upgrade \eqref{step 1} to `honest' negative regularity (i.e. in Sobolev
sense).  We start with

\begin{proposition}[Some negative regularity]\label{P:some -reg}
Let $d\geq 5$ and let $u$ be as in Theorem~\ref{T:-reg}.  Assume further that $|\nabla|^s F(u) \in L_t^\infty L_x^r$ for some
$\tfrac{2(d-2)(d+1)}{d^2+3d-6}\leq r<\tfrac{2d}{d+4}$ and some $0\leq s\leq 1$.  Then there exists $s_0=s_0(r,d)>0$ such that
$u \in L_t^\infty \dot H^{s-s_0+}_x$.
\end{proposition}

\begin{proof}
The proposition will follow once we establish
\begin{align}\label{neg reg 0}
\bigl\| |\nabla|^s u_N \bigr\|_{L_t^\infty L_x^2} \lesssim_u N^{s_0} \quad \text{for all} \quad N>0 \quad \text{and} \quad
s_0:=\tfrac{d}r-\tfrac{d+4}2>0.
\end{align}
Indeed, by Bernstein combined with \eqref{ke bounded},
\begin{align*}
\bigl\| |\nabla|^{s-s_0+} u \bigr\|_{L_t^\infty L_x^2}
&\leq \bigl\| |\nabla|^{s-s_0+} u_{\leq 1} \bigr\|_{L_t^\infty L_x^2} + \bigl\| |\nabla|^{s-s_0+} u_{>1} \bigr\|_{L_t^\infty L_x^2}\\
&\lesssim_u \sum_{N\leq 1} N^{0+} + \sum_{N>1} N^{(s-s_0+)-1}\\
&\lesssim_u 1.
\end{align*}

Thus, we are left to prove \eqref{neg reg 0}.  By time-translation symmetry, it suffices to prove
\begin{align}\label{neg reg 1}
\bigl\| |\nabla|^s u_N(0) \bigr\|_{L_x^2} \lesssim_u N^{s_0} \quad \text{for all} \quad N>0 \quad \text{and} \quad
s_0:=\tfrac{d}r-\tfrac{d+4}2>0.
\end{align}
Using the Duhamel formula \eqref{Duhamel} both in the future and in the past, we write
\begin{align*}
\bigl\| |\nabla|^s & u_N(0) \bigr\|_{L_x^2}^2 \\
&= \lim_{T\to\infty} \lim_{T'\to-\infty} \bigl \langle i\int_0^T e^{-it\Delta}P_N |\nabla|^s F(u(t))\, dt,
    -i\int_{T'}^0 e^{-i\tau\Delta}P_N |\nabla|^s F(u(\tau))\, d\tau  \bigr\rangle\\
&\leq  \int_0^\infty \int_{-\infty}^0 \Bigl| \bigl\langle P_N |\nabla|^s F(u(t)) ,
    e^{i(t-\tau)\Delta}P_N |\nabla|^s F(u(\tau))  \bigr\rangle \Bigr| \,dt\, d\tau.
\end{align*}
We estimate the term inside the integrals in two ways.  On one hand, using H\"older and the dispersive estimate,
\begin{align*}
\Bigl|\bigl\langle P_N |\nabla|^s F(u(t)) ,  e^{i(t-\tau)\Delta} & P_N |\nabla|^s F(u(\tau))  \bigr\rangle\Bigr|\\
&\lesssim \bigl\| P_N |\nabla|^s F(u(t))\bigr\|_{L_x^r} \bigl\| e^{i(t-\tau)\Delta} P_N |\nabla|^s F(u(\tau))\bigr\|_{L_x^{r'}}\\
&\lesssim |t-\tau|^{\frac d2-\frac dr} \bigl\| |\nabla|^s F(u)\bigr\|_{L_t^\infty L_x^r}^2.
\end{align*}
On the other hand, using Bernstein,
\begin{align*}
\Bigl|\bigl\langle P_N |\nabla|^s F(u(t)) ,  e^{i(t-\tau)\Delta} & P_N |\nabla|^s F(u(\tau))  \bigr\rangle\Bigr|\\
&\lesssim \bigl\| P_N |\nabla|^s F(u(t))\bigr\|_{L_x^2} \bigl\| e^{i(t-\tau)\Delta} P_N |\nabla|^s F(u(\tau))\bigr\|_{L_x^2}\\
&\lesssim N^{2(\frac d2-\frac dr)} \bigl\| |\nabla|^s F(u)\bigr\|_{L_t^\infty L_x^r}^2.
\end{align*}
Thus,
\begin{align*}
\bigl\| |\nabla|^s u_N(0) \bigr\|_{L_x^2}^2
&\lesssim \bigl\| |\nabla|^s F(u)\bigr\|_{L_t^\infty L_x^r}^2 \int_0^\infty \int_{-\infty}^0 \min\{ |t-\tau|^{-1} , N^2 \}^{\frac dr -\frac d2} \, dt\, d\tau\\
&\lesssim N^{2s_0} \bigl\| |\nabla|^s F(u)\bigr\|_{L_t^\infty L_x^r}^2.
\end{align*}
To obtain the last inequality we used the fact that $\tfrac dr -\tfrac d2 >2$ since $r<\tfrac{2d}{d+4}$.
Thus \eqref{neg reg 1} holds; this finishes the proof of the proposition.
\end{proof}

The proof of Theorem~\ref{T:-reg} will follow from iterating Proposition~\ref{P:some -reg} finitely many times.

\begin{proof}[Proof of Theorem~\ref{T:-reg}]
Proposition~\ref{P:L^p breach} allows us to apply Proposition~\ref{P:some -reg} with $s=1$.
We conclude that $u \in L_t^\infty \dot H^{1-s_0+}_x$ for some $s_0=s_0(r,d)>0$.
Combining this with the fractional chain rule Lemma~\ref{F Lip} and \eqref{step 1}, we deduce that $|\nabla|^{1-s_0+} F(u) \in L_t^\infty L_x^r$
for some $\tfrac{2(d-2)(d+1)}{d^2+3d-6}\leq r<\tfrac{2d}{d+4}$.  We are thus in the position to apply Proposition~\ref{P:some -reg} again
and obtain $u \in L_t^\infty \dot H^{1-2s_0+}_x$.  Iterating this procedure finitely many times, we derive
$u\in L_t^\infty \dot H_x^{-\eps}$ for any $0<\eps<s_0$.

This completes the proof of Theorem~\ref{T:-reg}.
\end{proof}

%
%
%
%

\section{The low-to-high frequency cascade}\label{S:cascade}

In this section, we use the negative regularity provided by Theorem~\ref{T:-reg} to preclude low-to-high frequency cascade solutions.

\begin{theorem}[Absence of cascades]\label{T:no cascade}
Let $d\geq 5$.  There are no global solutions to \eqref{nls} that are low-to-high frequency cascades
in the sense of Theorem~\ref{T:enemies}.
\end{theorem}

\begin{proof}
Suppose for a contradiction that there existed such a solution $u$.  Then, by Theorem~\ref{T:-reg}, $u\in L_t^\infty L_x^2$;
thus, by the conservation of mass,
$$
0\leq M(u)=M(u(t)) = \int_{\R^d} |u(t,x)|^2\, dx <\infty \quad \text{for all} \quad t\in \R.
$$

Fix $t\in \R$ and let $\eta>0$ be a small constant.  By compactness (see Remark~\ref{R:c small}),
$$
\int_{|\xi|\leq c(\eta)N(t)} |\xi|^2 |\hat u(t,\xi)|^2\, d\xi\leq \eta.
$$
On the other hand, as $u\in L_t^\infty \dot H^{-\eps}_x$ for some $\eps>0$,
$$
\int_{|\xi|\leq c(\eta)N(t)} |\xi|^{-2\eps} |\hat u(t,\xi)|^2\, d\xi\lesssim_u 1.
$$
Hence, by H\"older's inequality,
\begin{align}\label{close}
\int_{|\xi|\leq c(\eta)N(t)} |\hat u(t,\xi)|^2\, d\xi\lesssim_u \eta^{\frac \eps{1+\eps}}.
\end{align}

Meanwhile, by elementary considerations and \eqref{ke bounded},
\begin{align}\label{far}
\int_{|\xi|\geq c(\eta)N(t)} |\hat u(t,\xi)|^2\, d\xi
&\leq [c(\eta)N(t)]^{-2} \int_{\R^d} |\xi|^2 |\hat u(t,\xi)|^2\, d\xi \notag \\
&\leq [c(\eta)N(t)]^{-2} \|\nabla u(t)\|^2_2 \notag\\
&\leq [c(\eta)N(t)]^{-2} \|\nabla W\|^2_2.
\end{align}

Collecting \eqref{close} and \eqref{far} and using Plancherel's theorem, we obtain
$$
0\leq M(u)\lesssim_u c(\eta)^{-2} N(t)^{-2} + \eta^{\frac \eps{1+\eps}}
$$
for all $t\in \R$.  As $u$ is a low-to-high cascade, there is a sequence of times $t_n\to\infty$ so that
$N(t_n)\to \infty$.  As $\eta>0$ is arbitrary, we may conclude $M(u)=0$ and hence $u$ is identically zero.
This contradicts the fact that $S_I(u)=\infty$, thus settling Theorem~\ref{T:no cascade}.
\end{proof}

%
%
%
%

\section{The soliton}\label{S:Soliton}
In this case the contradiction will follow from a virial-type argument.  In order to successfully use the virial
inequality, we need to control the motion of $x(t)$.  As we now know that soliton solutions have finite mass
(see Theorem~\ref{T:-reg}), we will be able to use an argument from \cite{DHR}
to prove that $|x(t)|=o(t)$ as $t\to \infty$.  The first step is to note that a minimal kinetic energy
blowup solution with finite mass must have zero momentum.  Let us quickly remark that simple arguments show
$|x(t)|=O(t)$ without the knowledge that the mass is finite; however, the passage from $O(t)$ to $o(t)$ is essential
for the virial argument.

\begin{proposition}[Zero momentum]\label{P:p=0}
Let $u$ be a minimal kinetic energy blowup solution to \eqref{nls} that obeys $u\in L_t^\infty H^1_x$.
Then its total momentum, which is a conserved quantity, vanishes:
$$P(u):=2\, \Im \int_{\R^d} \overline{u(t,x)}\nabla u(t,x)\, dx \equiv 0.$$
\end{proposition}

\begin{proof}
Let $u:I\times\R^d\to \C$ be as in Proposition~\ref{P:p=0}.  Then the momentum $P(u)$ and the mass $M(u)$ are finite and conserved.
Moreover, $M(u)\neq 0$ since otherwise $u$ would be identically zero and hence not a blowup solution.

Let $\tilde u$ be the Galilei boost of $u$ by $\xi_0:= - [2M(u)]^{-1} P(u)$:
$$
\tilde u(t,x):= e^{ix\xi_0}e^{-it|\xi_0|^2} u(t,x-2\xi_0t).
$$
Simple computations then show that
\begin{equation}\label{ECofM}
\|\nabla \tilde u(t)\|_2^2
= \|\nabla u(t)\|_2^2 + |\xi_0|^2 M(u) + \xi_0 P(u)
= \|\nabla u(t)\|_2^2 - [4M(u)]^{-1}P(u)^2.
\end{equation}
Equivalently, we may write $E(u) = E(\tilde u) + [4M(u)]^{-1}P(u)^2$, which expresses the well-know physical fact that the
total energy can be decomposed as the energy viewed in the center of mass frame plus the energy arising from the motion of
the center of mass (cf \cite[\S 8]{ll}).

As $\tilde u$ is also a blowup solution of \eqref{nls}, indeed $S_I(\tilde u)=S_I(u)=\infty$, we see that $P(u)=0$;
for otherwise, $\tilde u$ would have less kinetic energy than $u$.
\end{proof}

A second ingredient needed to control the motion of $x(t)$ is a compactness property of the orbit $\{u(t)\}$ in $L_x^2$.
This requires the full force of Theorem~\ref{T:-reg}.

\begin{lemma}[Compactness in $L_x^2$]\label{L:mass compact}
Let $d\geq 5$ and let $u$ be a soliton in the sense of Theorem~\ref{T:enemies}.
Then for every $\eta>0$ there exists $C(\eta)>0$ such that
$$
\sup_{t\in \R} \int_{|x-x(t)|\geq C(\eta)}|u(t,x)|^2\, dx \lesssim_u \eta.
$$
\end{lemma}

\begin{proof}
The entire argument takes place at a fixed $t$; in particular, we may assume $x(t)=0$.

First we control the contribution from the low frequencies: by Theorem~\ref{T:-reg},
$$
\bigl\|u_{< N}(t)\bigl\|_{L_x^2(|x|\geq R)}
\leq \bigl\|u_{< N}(t)\bigl\|_{L_x^2} \lesssim N^{\eps} \bigl\| |\nabla|^{-\eps} u\bigl\|_{L_t^\infty L_x^2}\lesssim_u N^{\eps}.
$$
This can be made smaller than $\eta$ by choosing $N=N(\eta)$ small enough.

We now turn to the contribution from the high frequencies.
A simple application of Schur's test reveals the following: For any $m\geq 0$,
$$
\bigl\| \chi_{|x|\geq 2R} \Delta^{-1} \nabla P_{\geq N} \chi_{|x|\leq R} \bigr\|_{L^2\to L^2}
    \lesssim N^{-1} \langle RN \rangle^{-m}
$$
uniformly in $R,N>0$.  On the other hand, by Bernstein,
$$
\bigl\| \chi_{|x|\geq 2R} \Delta^{-1} \nabla P_{\geq N} \chi_{|x|\geq R} \bigr\|_{L^2\to L^2} \lesssim N^{-1}.
$$
Together, these lead quickly to
\begin{align*}
\int_{|x|\geq 2R} |u_{\geq N}(t,x)|^2 \,dx \lesssim N^{-2} \langle RN \rangle^{-2} \|\nabla u(t)\|_{L^2_x}^2
    + N^{-2} \int_{|x|\geq R} |\nabla u(t,x)|^2 \,dx.
\end{align*}
By choosing $R$ large enough, we can render the first term smaller than $\eta$; the same is true of the
second summand by virtue of $\dot H^1$-compactness:
\begin{align*}
\sup_{t\in\R} \int_{|x-x(t)|\geq C(\eta_1)}|\nabla u(t,x)|^2\, dx \leq \eta_1.
\end{align*}

The lemma follows by combining our estimates for $u_{< N}$ and $u_{\geq N}$.
\end{proof}

Following the argument in \cite{DHR}, we can now prove

\begin{lemma}[Control over $x(t)$]\label{L:x(t) control}
Fix $d\geq 5$ and let $u$ be a minimal kinetic energy soliton in the sense of Theorem~\ref{T:enemies}.  Then
\begin{align*}
|x(t)| =o(t) \quad \text{as} \quad t\to \infty.
\end{align*}
\end{lemma}

\begin{proof}
We argue by contradiction.  Suppose there exist $\delta>0$ and
a sequence $t_n\to \infty$ such that
\begin{align}\label{assume not}
|x(t_n)|> \delta t_n \quad \text{for all} \quad n\geq 1.
\end{align}
By spatial-translation symmetry, we may assume $x(0)=0$.

Let $\eta>0$ be a small constant to be chosen later.   By compactness and Lemma~\ref{L:mass compact},
\begin{align}\label{compact}
\sup_{t\in \R} \int_{|x-x(t)|>C(\eta)}\bigl( |\nabla u(t,x)|^2 + |u(t,x)|^2 \bigr)\, dx\leq \eta
\end{align}
Define
\begin{align}\label{define}
T_n:=\inf_{t\in[0, t_n]} \{|x(t)|=|x(t_n)|\}\leq t_n \quad \text{and} \quad R_n:=C(\eta) + \sup_{t\in[0,T_n]}|x(t)|.
\end{align}

Now let $\phi$ be a smooth, radial function such that
\begin{align*}
\phi(r)=\begin{cases}
1 & \text{for } r\leq 1\\
0 & \text{for } r\geq 2,
\end{cases}
\end{align*}
and define the truncated `position'
$$
X_R(t): =\int_{\R^d} x\phi\bigl(\tfrac{|x|}R\bigr) |u(t,x)|^2\,dx.
$$
By Theorem~\ref{T:-reg}, $u\in L_t^\infty L_x^2$; together with \eqref{compact} this implies
\begin{align*}
|X_{R_n}(0)|
&\leq \bigl|\int_{|x|\leq C(\eta)}x\phi\bigl(\tfrac{|x|}R\bigr) |u(t,x)|^2\,dx\Bigr|
 + \bigl|\int_{|x|\geq C(\eta)} x\phi\bigl(\tfrac{|x|}R\bigr) |u(t,x)|^2\,dx\Bigr|\\
&\leq C(\eta) M(u)  + 2 \eta R_n.
\end{align*}
On the other hand, by the triangle inequality combined with \eqref{compact} and \eqref{define},
\begin{align*}
|X_{R_n}(T_n)|
&\geq |x(T_n)|M(u) - |x(T_n)| \Bigl| \int_{\R^d}\Bigl[1-\phi\bigl(\tfrac{|x|}R\bigr)\bigr] |u(T_n,x)|^2\,dx\Bigr|\\
& \quad- \Bigl| \int_{|x-x(T_n)|\leq C(\eta)}\bigl[x-x(T_n)\bigr]\phi\bigl(\tfrac{|x|}R\bigr) |u(T_n,x)|^2\,dx\Bigr|\\
& \quad - \Bigl| \int_{|x-x(T_n)|\geq C(\eta)}\bigl[x-x(T_n)\bigr]\phi\bigl(\tfrac{|x|}R\bigr) |u(T_n,x)|^2\,dx\Bigr|\\
& \geq |x(T_n)| [M(u)-\eta] - C(\eta)M(u) - \eta [R_n + |x(T_n)|]\\
& \geq |x(T_n)| [M(u)-3\eta] - 2C(\eta)M(u).
\end{align*}
Thus, taking $\eta>0$ sufficiently small (depending on $M(u)$),
\begin{align*}
\bigl| X_{R_n}(T_n) - X_{R_n}(0) \bigr|\gtrsim_{M(u)}|x(T_n)| - C(\eta).
\end{align*}

A simple computation establishes
\begin{align*}
\partial_t X_R(t)
&= 2\Im \int_{\R^d} \phi\bigl(\tfrac{|x|}R\bigr)\nabla u(t,x) \overline{u(t,x)}\, dx\\
&\quad + 2\Im \int_{\R^d} \frac{x}{|x|R}\phi'\bigl(\tfrac{|x|}R\bigr)\,x\cdot \nabla u(t,x) \overline{u(t,x)}\, dx.
\end{align*}
By Proposition~\ref{P:p=0}, $P(u)=0$; together with Cauchy-Schwarz and \eqref{compact} this yields
\begin{align*}
\bigl|\partial_t X_{R_n}(t)\bigr|
&\leq \Bigl| 2\Im \int_{\R^d}\Bigl[1- \phi\bigl(\tfrac{|x|}R\bigr)\Bigr]\nabla u(t,x) \overline{u(t,x)}\, dx \Bigr|\\
&\quad + \Bigl|2\Im \int_{\R^d} \frac{x}{|x|R}\phi'\bigl(\tfrac{|x|}R\bigr) \, x\cdot \nabla u(t,x) \overline{u(t,x)}\, dx\Bigr|\\
&\leq 6\eta
\end{align*}
for all $t\in [0, T_n]$.

Thus, by the Fundamental Theorem of Calculus,
$$
|x(T_n)| - C(\eta)\lesssim _{M(u)} \eta T_n.
$$
Recalling that $|x(T_n)|=|x(t_n)|>\delta t_n\geq \delta T_n$ and letting $n\to \infty$ we derive a contradiction.
\end{proof}

We are finally in a position to preclude the soliton-like enemy by using a truncated virial identity.  When $x(t)\equiv0$,
as in the radial case, the necessary argument can be found in \cite{kenig-merle}.
As the reader will see, it is the finiteness of the $L^2_x$ norm that allows us to extend the argument to the case $|x(t)|=o(t)$.

\begin{theorem}[No soliton]\label{T:no soliton}
Let $d\geq 5$.  A minimal kinetic energy blowup solution of \eqref{nls} cannot be a soliton in the sense of Theorem~\ref{T:enemies}.
\end{theorem}

\begin{proof}
Suppose for a contradiction that there existed such a solution $u$.

Let $\eta>0$ be a small constant to be specified later.  Then, by Definition~\ref{D:ap} and Remark~\ref{R:pot energy},
\begin{equation}\label{decay u}
\sup_{t\in \R} \int_{|x-x(t)|>C(\eta)}\bigl( |\nabla u(t,x)|^2+|u(t,x)|^{\tdt} \bigr)\, dx\leq \eta.
\end{equation}
Moreover, by Lemma~\ref{L:x(t) control}, $|x(t)|=o(t)$ as $t \to \infty$.  Thus, there exists $T_0=T_0(\eta)\in \R$ such that
\begin{align}\label{control}
|x(t)|\leq \eta t \quad \text{for all} \quad t\geq T_0.
\end{align}

Now let $\phi$ be a smooth, radial function such that
\begin{align*}
\phi(r)=\begin{cases}
r & \text{for } r\leq 1\\
0 & \text{for } r\geq 2,
\end{cases}
\end{align*}
and define
$$
V_R(t): =\int_{\R^d} \psi(x) |u(t,x)|^2\,dx,
$$
where $\psi(x):=R^2\phi\bigl(\tfrac{|x|^2}{R^2}\bigr)$ for some $R>0$.

Differentiating $V_R$ with respect to the time variable, we find
$$
\partial_t V_R(t) = 4\Im \int_{\R^d} \phi'\bigl(\tfrac{|x|^2}{R^2}\bigr) \overline{u(t,x)} \; x\cdot \nabla u(t,x) \,dx.
$$
By Theorem~\ref{T:-reg}, $u\in L_t^\infty L_x^2$ and so
\begin{align}\label{one deriv}
|\partial_t V_R(t)|\lesssim R \|\nabla u(t)\|_2 \|u(t)\|_2
                    \lesssim_u R
\end{align}
for all $t\in I$ and $R>0$.

Further computations establish
\begin{align*}
\partial_{tt}V_R(t) &= 4 \Re \int_{\R^d}\psi_{ij}(x)u_{i}(t,x)\bar u_j(t,x)\, dx
        - \tfrac 4d \int_{\R^d} \bigl(\Delta \psi\bigr)(x) |u(t,x)|^{\frac{2d}{d-2}}\, dx\\
&\quad  - \int_{\R^d} \bigl(\Delta \Delta \psi\bigr)(x) |u(t,x)|^2\, dx,
\end{align*}
where subscripts denote spatial derivatives and repeated indices are summed. Substituting our choice of $\psi$ and using
H\"older's inequality on the last term,
\begin{align*}
\partial_{tt}V_R(t)
&=8\int_{\R^d} \bigl(|\nabla u(t,x)|^2-|u(t,x)|^{\tdt}\bigr)\, dx\\
&\quad + O\Bigl(\int_{|x|\geq R} \bigl(|\nabla u(t,x)|^2 + |u(t,x)|^{\tdt}\bigr)\, dx\Bigr)\\
&\quad + O\Bigl(\int_{R\leq |x|\leq 2R} |u(t,x)|^{\tdt}\, dx\Bigr)^{\frac{d-2}{d}}.
\end{align*}
From \eqref{hype} and Lemma~\ref{equal},
\begin{align*}
\int_{\R^d} \bigl(|\nabla u(t,x)|^2-|u(t,x)|^{\tdt}\bigr)\, dx \gtrsim \|\nabla u_0\|_2^2.
\end{align*}
Thus, choosing $\eta>0$ sufficiently small and $R:=C(\eta) +\sup_{T_0\leq t\leq T_1} |x(t)|$ and invoking \eqref{decay u},
\begin{align}\label{two deriv}
\partial_{tt}V_R(t) \gtrsim \|\nabla u_0\|_2^2.
\end{align}

Using the Fundamental Theorem of Calculus on the interval $[T_0, T_1]$ together with \eqref{one deriv} and \eqref{two deriv}, we obtain
$$
(T_1-T_0) \|\nabla u_0\|_2^2\lesssim_u R \lesssim_u C(\eta) +\sup_{T_0\leq t\leq T_1} |x(t)|
$$
for all $T_1\geq T_0$.  Invoking \eqref{control} and taking $\eta$ sufficiently small and $T_1$ sufficiently large,
we derive a contradiction unless $u_0\equiv 0$.  But $u_0\equiv 0$ is not consistent with the fact that for a soliton, $S_\R(u)=\infty$.
\end{proof}

%
%
%
%

\section{Blowup}\label{S:blowup}

In this section we prove Proposition~\ref{P:blowup}.  To this end, let $u_0\in \dot H^1_x(\R^d)$ and $\delta_0>0$ be such that
$$
\|\nabla u_0\|_2 \geq \|\nabla W\|_2 \quad \text{and} \quad E(u_0)\leq (1-\delta_0) E(W).
$$
Let $u:\ird\to \C$ be the maximal-lifespan solution to \eqref{nls} with initial data $u_0$ at time $t=t_0\in I$. By
Corollary~\ref{trap}, there exist $\delta_2, \delta_3>0$ such that for all $t\in I$
\begin{align}
\|\nabla u(t)\|_2^2 &\geq (1+\delta_2) \|\nabla W\|_2^2 \label{high ke}\\
\int_{\R^d}\bigl(|\nabla u(t,x)|^2-|u(t,x)|^{\tdt}\bigr)\,dx &\leq -\delta_3. \label{neg viriel}
\end{align}

To prove that the solution $u$ blows up in finite time (in either of the two cases described in Proposition~\ref{P:blowup}), we
will use the convexity method \cite{glassey, zack}.

Let us first treat the case when $xu_0\in L_x^2(\R^d)$; see also \cite{kenig-merle}.  In this case, the second moment
$$
V(t):=\int_{\R^d} |x|^2 |u(t,x)|^2 \, dx
$$
is well-defined and moreover, $V\in C^2(I)$; see, for example, \cite{cazenave:book}.  As $u$ is not identically zero (by
\eqref{high ke}), $V(t)>0$ for all $t\in I$. On the other hand, a quick computation together with \eqref{neg viriel} shows
$$
\partial_{tt} V(t) = 8\int_{\R^d}\bigl(|\nabla u(t,x)|^2-|u(t,x)|^{\tdt}\bigr)\,dx \leq -8\delta_3.
$$
Thus, the graph of $V$ lies under an inverted parabola, so the solution $u$ blows up in both time directions.

We consider next the case when $u_0\in H^1_x(\R^d)$ is radial.  Blowup for the energy-subcritical problem for this type of
initial data and negative energy was addressed by Ogawa and Tsutsumi \cite{OgawaTsutsumi}; see also \cite{cazenave:book}.
While our exposition is a little different, the argument is very close to theirs.

As the second moment is no longer finite for this initial data, we define the truncated virial quantity
$$
V_R(t):=\int_{\R^d}  \psi(x)|u(t,x)|^2\,dx,
$$
where $\psi(x):=R^2 \phi\bigl(\frac{|x|^2}{R^2}\bigr)$ with $R>0$ and $\phi$ a smooth, concave function on $[0,\infty)$ such
that
\begin{equation*}
\phi (r) = \begin{cases}
r & \text{for } r\leq 1\\
2  & \text{for } r\geq 3
\end{cases}
\end{equation*}
and $\phi''(r)$ is non-increasing on $r\leq 2$ and non-decreasing on $r\geq 2$.

A computation establishes
\begin{align*}
\partial_{tt}V_R(t) &= 4 \Re \int_{\R^d}\psi_{ij}(x)u_{i}(t,x)\bar u_j(t,x)\, dx
        - \tfrac 4d \int_{\R^d} \bigl(\Delta \psi\bigr)(x) |u(t,x)|^{\frac{2d}{d-2}}\, dx\\
&\quad  - \int_{\R^d} \bigl(\Delta \Delta \psi\bigr)(x) |u(t,x)|^2\, dx.
\end{align*}
Substituting our choice of $\psi$ in the formula above and recalling that $u$ is radial,
\begin{align*}
\partial_{tt}V_R(t)
&= 8\int_{\R^d}\bigl( |\nabla u(t,x)|^2 - |u(t,x)|^{\frac{2d}{d-2}}\bigr)\, dx + \tfrac1{R^2} O\Bigl( \int_{|x|\sim R} \! |u(t,x)|^2 \, dx\Bigr)\\
&\quad + 8\int_{\R^d}\! \Bigl( \phi'\bigl(\tfrac{|x|^2}{R^2}\bigr) -1 + \tfrac{2|x|^2}{R^2}
\phi''\bigl(\tfrac{|x|^2}{R^2}\bigr)\Bigr)
        \bigl(|\nabla u(t,x)|^2- |u(t,x)|^{\frac{2d}{d-2}}\bigr)\,dx\\
&\quad + \tfrac{16(d-1)}{d} \int_{\R^d} \tfrac{|x|^2}{R^2} \phi''\bigl(\tfrac{|x|^2}{R^2}\bigr)|u(t,x)|^{\frac{2d}{d-2}}\,dx.
\end{align*}
By our choice of $\phi$, we have $\phi''\leq 0$.  Moreover, as $u\in L_x^2(\R^d)$, one can choose $R$ sufficiently large
(depending on the mass of $u$) so that the contribution of the second term on the right-hand side of the equality above is less
than half that of the first term. Thus, invoking \eqref{neg viriel},
\begin{align}\label{deriv V}
\partial_{tt}V_R(t) \leq -4 \delta_3
    - 8\int_{\R^d}\! \omega(x)\bigl(|\nabla u(t,x)|^2- |u(t,x)|^{\frac{2d}{d-2}}\bigr) \, dx,
\end{align}
where $\omega(x):=1-\phi'\bigl(\tfrac{|x|^2}{R^2}\bigr) - \tfrac{2|x|^2}{R^2} \phi''\bigl(\tfrac{|x|^2}{R^2}\bigr)$. Note that
$0\leq \omega\lesssim 1$ is radial, $\supp(\omega)\subseteq \{ |x|\geq R\}$, and $\omega(x)\lesssim \omega(y)$ uniformly
for all $|x|\leq |y|$.

As in the first case, finite time blowup for $u$ will follow once we establish $\partial_{tt}V_R<0$. To achieve this we will
need the following

\begin{lemma}[Weighted radial Sobolev embedding]\label{wait}
Let $\omega$ be as above and let $f$ be a radial function on $\R^d$.  Then
\begin{align*}
\bigl\| |x|^{\frac{d-1}2} \omega^{\frac 14} f \bigr\|_{L_x^\infty(\R^d)}^2 \lesssim \| f \bigr\|_{L_x^2(\R^d)} \|\omega^{\frac
12} \nabla f \bigr\|_{L_x^2(\R^d)}.
\end{align*}
\end{lemma}

\begin{proof}
It suffices to establish the claim for radial Schwartz functions $f$.  Let $r\geq 0$.  By the Cauchy--Schwarz inequality,
\begin{align*}
r^{d-1}\omega(r)^{\frac12} |f(r)|^2
&= 2 r^{d-1}\omega(r)^{\frac12} \Re \int_r^\infty \bar f(\rho) f'(\rho) \, d\rho\\
&\lesssim \int_r^\infty \rho^{d-1}\omega(\rho)^{\frac12} |f(\rho)|\,|f'(\rho)| \, d\rho\\
&\lesssim \Bigl(\int_r^\infty \rho^{d-1} |f(\rho)|^2 \, d\rho \Bigr)^{\frac12}
        \Bigl(\int_r^\infty \rho^{d-1} \omega(\rho) |f'(\rho)|^2 \, d\rho \Bigr)^{\frac12}\\
&\lesssim \| f \bigr\|_{L^2(\rho^{d-1}d\rho)}\|\omega^{\frac 12} f' \bigr\|_{L^2(\rho^{d-1}d\rho)}.
\end{align*}
The claim follows.
\end{proof}

Returning to the proof of Proposition~\ref{P:blowup}, we use Lemma~\ref{wait} together with the fact that $u_0\in L_x^2(\R^d)$,
the conservation of mass, and the properties of $\omega$ described earlier to estimate
\begin{align*}
\int_{\R^d}\! \omega(x) |u(t,x)|^{\frac{2d}{d-2}}\, dx
&\lesssim \|\omega^{\frac14}u(t)\|_{L_x^\infty}^{\frac4{d-2}} \int_{\R^d} |u(t,x)|^2\, dx\\
&\lesssim R^{-\frac{2(d-1)}{d-2}}\bigl\| |x|^{\frac{d-1}2} \omega^{\frac 14} u(t) \bigr\|_{L_x^\infty}^{\frac4{d-2}}\|u_0\|_{L_x^2}^2\\
&\lesssim R^{-\frac{2(d-1)}{d-2}}\|\omega^{\frac12} \nabla u(t)\|_{L_x^2}^{\frac2{d-2}}\|u_0\|_{L_x^2}^{\frac{2(d-1)}{d-2}}\\
&\lesssim \bigl(R^{-1}\|u_0\|_{L_x^2}\bigr)^{\frac{2(d-1)}{d-2}} \bigl(\|\omega^{\frac12} \nabla u(t)\|_{L_x^2}^2 +1\bigr).
\end{align*}
Thus, taking $R$ sufficiently large depending on the mass of $u$ and recalling that $\omega$ is positive, \eqref{deriv V} yields
$\partial_{tt}V_R<0$.  This finishes the proof of Proposition~\ref{P:blowup}. \qed

%
%
%
%

\section{Concentration at blowup}\label{S:Conc}

The paper of Kenig and Merle contains a sketch of an argument to prove Theorem~\ref{T:conc} in the low dimensional
spherically symmetric case treated in that paper; see \cite[Corollary~5.18]{kenig-merle}.
As far as we understand it, that paper does not satisfactorily address the problem of quadratic oscillation,
namely, that there are radial functions $\phi_n$ obeying
$$
\|\phi_n\|_{L^2_x(\R^d)} =1,\quad \|e^{it\Delta} \phi_n\|_{L^{\frac{2(d+2)}{d}}_{t,x}([0,1]\times\R^d)} \gtrsim 1,
    \quad \text{but}\quad \int_{|x|\leq R} |\phi_n(x)|^2\,dx \to 0
$$
for all $R>0$.  This difficulty is described, for instance, in the works of Merle--Vega \cite{merle-vega} and Keraani
\cite{keraani-l2} on the mass-critical equation.

The approach we take here is inspired by \cite{bourg.2d,KTV}.  These papers considered the mass-critical equation; however,
unlike mass, kinetic energy is not conserved and this leads to several additional difficulties.  One example is that we
need to appeal to the full strength of Theorem~\ref{T:main}; Corollary~\ref{C:main} is not sufficient.

\begin{proof}[Proof of Theorem~\ref{T:conc}]
Without loss of generality, we may assume that the solution $u$ blows up forward in time at $0<T^*\leq \infty$. We will further
assume that $T^*$ is finite; the proof in the case when $T^*=\infty$ requires only a few minor changes.

Let $t_n\nearrow T^*$ and define $u_n(0):=u(t_n)$; then each $u_n$ is a solution to \eqref{nls} on $[0, T^*-t_n)$. Invoking
\eqref{type II}, we apply Lemma~\ref{L:cc} to decompose
$$
u_n(0)= \sum_{j=1}^J \gnj e^{i\tnj \Delta} \phi^j +\wnJ
$$
and define $\vnj:I_n^j\times\R^d\to \C$ to be the maximal-lifespan solution to \eqref{nls} with initial data
\begin{align*}
\vnj(0):=\gnj e^{i\tnj \Delta} \phi^j.
\end{align*}

By \eqref{decouple}, there exists $J_0\geq 1$ such that
$$
\|\nabla \phi^j\|_2\leq \eta_0 \quad \text{for all} \quad j\geq J_0,
$$
where $\eta_0=\eta_0(d)$ is the threshold from the small data theory.  By Theorem~\ref{T:local}, for all $n\geq 1$ and $j\geq
J_0$, the solutions $\vnj$ are global and moreover,
$$
\sup_{t\in \R}\|\nabla \vnj(t)\|_2^2 + S_{\R}(\vnj)\lesssim \|\nabla \phi^j\|_2^2 \quad \text{for all } n\geq 1 \text{ and }
j\geq J_0.
$$

\begin{lemma}[A bad profile]\label{L:bad profile II}
There exists $1\leq j_0<J_0$ such that
\begin{align}\label{E:bad}
\limsup_{n\to \infty}S_{[0, T^*-t_n)}(v_n^{j_0})=\infty.
\end{align}
\end{lemma}

\begin{proof}
The proof is similar to that of Lemma~\ref{L:bad profile}.  Arguing by contradiction, one approximates $u_n$ by
$\unJ:=\sum_{j=1}^J \vnj + e^{it\Delta}\wnJ$ and concludes that if \eqref{E:bad} were to fail, then $u$ would not blow up at
$T^*$.
\end{proof}

Reordering the indices, we may assume that there exists $1\leq J_1< J_0$ such that
\begin{align}\label{J1 def}
\limsup_{n\to \infty}S_{[0, T^*-t_n)}(\vnj)=\infty \quad \text{for all} \quad j\leq J_1
\end{align}
and
$$
\limsup_{n\to \infty}S_{[0, T^*-t_n)}(\vnj)<\infty \quad \text{for all} \quad j> J_1.
$$
To continue, we pass to a subsequence in $n$ so that $S_{[0,T^*-t_n)}(v_n^1) \to \infty$.

For each $m,n\geq 1$, there exist $1\leq j(m,n) \leq J_1$ and $0<T_n^m<T^*-t_n$ such that
\begin{equation}\label{S is m}
\sup_{1\leq j \leq J_1} S_{[0, T_n^m]}(\vnj) = S_{[0, T_n^m]}(v_n^{j(m,n)}) = m.
\end{equation}
By the pigeonhole principle, there exists $1\leq j_1 \leq J_1$ so that for infinitely many $m$ we have $j(m,n)=j_1$ for
infinitely many $n$.  Reordering the indices if necessary, we may assume $j_1=1$. Now by Theorem~\ref{T:main}, there exists $0\leq
\tau_n^m \leq T_n^m$ such that
$$
\limsup_{m\to \infty} \limsup_{n\to\infty} \|\nabla v_n^1(\tau_n^m)\|_2 \geq \|\nabla W\|_2.
$$

Given $\eps>0$, we may set $m_0=m_0(\eps)$ so that
\begin{align*}
\|\nabla v_n^1(\tau_n^{m_0})\|_2\geq \|\nabla W\|_2 -\eps \quad \text{for infinitely many } n.
\end{align*}
In what follows we will drop the superscript $m_0$ and denote $\tau_n^{m_0}$, $T_n^{m_0}$ by $\tau_n$, $T_n$, respectively.  We
will also pass to an $\eps$-dependent subsequence in $n$ so that
\begin{align}\label{big ke}
 \|\nabla v_n^1(\tau_n)\|_2\geq \|\nabla W\|_2 -\eps \ \text{ for all $n$ and }
    \ \lim_{n\to\infty} \|\nabla v_n^1(\tau_n)\|_2 \ \text{ exists.}
\end{align}

Let $\eta>0$ be a small constant to be chosen later.  Fix $n$.  As $S_{[0, T_n]}(v_n^1)=m_0$, there exists an interval
$[\tau_n^-, \tau_n^+]\subset [0, T_n]$ containing $\tau_n$ such that
\begin{equation}\label{eta S}
S_{[\tau_n^-, \tau_n^+]}(v_n^1)=\eta.
\end{equation}
Using \cite[Lemma~5.1]{KVZ-Q} as in that paper, one may deduce
$$
S_{[\tau_n^- -\tau_n, \tau_n^+ -\tau_n]}\bigl( e^{it\Delta} v_n^1(\tau_n) \bigr) \gtrsim \eta^{C(d)},
$$
where $C(d)$ is a dimension-dependent constant. Thus, by Lemma~\ref{L:B conc} there exist $x_n\in \R^d$ and $\tau_n^- -\tau_n\leq s_n\leq
\tau_n^+ -\tau_n$ such that
\begin{equation}\label{something is there}
\int_{|x - x_n|\lesssim  |T^*-t_n'|^{\frac12}} |e^{is_n\Delta} \nabla v_n^1(\tau_n)|^2 \, dx \gtrsim_\eta 1,
\end{equation}
where $t_n':= t_n + s_n +\tau_n$.  Note that we choose $s_n:=\inf J$ where $J$ is the interval from Lemma~\ref{L:B conc}. In the
case $T^*=\infty$, we choose $s_n:=\sup J$ and thus the diameter of the bubble is no more than $|t_n'|^{\frac12}$.

Passing to a subsequence in $n$, we may assume $t_n^1$ converges (possibly to $\pm\infty$).  Arguing as in the beginning of the
proof of Proposition~\ref{P:palais-smale}, we may assume that $t_n^1\equiv 0$ or $t_n^1\to\pm\infty$.  Continuing as there, we
define $v^1$ to be the maximal-lifespan solution that matches $\phi^1$ at $t=0$ (in the case $t_n^1\equiv0$) or scatters to
$e^{it\Delta} \phi^1$ (in the case $t_n^1\to\pm\infty$).  In truth, $t_n^1\to \infty$ is incompatible with \eqref{J1 def},
though we will not make use of this fact.

Tracing back through the definitions, we obtain
\begin{equation}\label{friends}
\bigl\| v_n^1 - T_{g_n^1} [ v^1(\cdot + t_n^1) ] \bigr\|_{\dot S^1([0,T_n])} \to 0.
\end{equation}
Note that in the case $t_n^1\equiv0$, the left-hand side is actually identically zero.  Combining this with \eqref{S is m} we
may deduce
\begin{equation}\label{can choose I}
S_{[t_n^1, t_n^1 + T_n(\lambda_n^1)^{-2}]}(v^1) \to S_{[0,T_n]}(v_n^1) = m_0.
\end{equation}
Moreover, combining \eqref{friends} with \eqref{something is there} and accounting for scaling and space translation,
\begin{equation*}
\int_{|\lambda_n^1 y + x_n^1 - x_n|\lesssim  |T^*-t_n'|^{\frac12}} |e^{is_n(\lambda^1_n)^{-2}\Delta} \nabla
v^1((\lambda^1_n)^{-2}\tau_n+t_n^1,y)|^2 \, dy \gtrsim_\eta 1.
\end{equation*}
We now apply Proposition~\ref{P:comp prof}, noting that \eqref{can choose I} implies that $v^1$ has finite scattering size on
the relevant interval.  Rescaling and invoking \eqref{friends}, we find
\begin{align}\label{bubble}
 \Bigl| \|\nabla v_n^1(\tau_n)\|_2^2 - \int_{|x- x_n^1|\leq R_n} |e^{is_n\Delta} \nabla v_n^1(\tau_n)|^2 \, dx \Bigr| \to 0
\end{align}
for any sequence $R_n\in (0, \infty)$ such that $(T^*-t_n')^{-\frac12}R_n\to \infty$ as $n\to \infty$.

It remains to show that a similar bubble can be found inside $u(t_n')$.  In view of \eqref{S is m}, we see that $\unJ$ is a good
approximation to $u_n$ on $[0, T_n]$ for $n$ and $J$ sufficiently large.  In particular,
\begin{align}\label{good approx ke II}
\lim_{J\to \infty}\limsup_{n\to \infty} \|\unJ(s_n + \tau_n) -u(t_n')\|_{\dot H^1_x(\R^d)}=0.
\end{align}
On the other hand, using \eqref{crazy} and Lemma~\ref{L:strong decouple}, one may deduce (arguing as in the proof of
Lemma~\ref{L:decouple ke}) that
\begin{align}\label{orthog II}
\limsup_{n\to \infty} \bigl|\bigl\langle \nabla \unJ(s_n +\tau_n), \nabla v_n^1(s_n+\tau_n) \bigr\rangle\bigr| = \limsup_{n\to
\infty}\|\nabla v_n^1(s_n +\tau_n)\|_2^2
\end{align}
for all $J\geq 1$.  Combining \eqref{good approx ke II} and \eqref{orthog II}, we derive
\begin{align*}
\limsup_{n\to \infty} \bigl|\bigl\langle \nabla u_n(t_n'), \nabla v_n^1(s_n+\tau_n) \bigr\rangle\bigr| = \limsup_{n\to
\infty}\|\nabla v_n^1(s_n +\tau_n)\|_2^2.
\end{align*}
Invoking \eqref{eta S} and using the Strichartz inequality, we see that
$$
\bigl\| v_n^1(s_n+\tau_n)- e^{is_n\Delta} \nabla v_n^1(\tau_n) \bigr\|_{\dot H^1_x} \lesssim \eta^{\frac2{d+2}}.
$$
Applying this on both sides of the equality above leads to
\begin{align}\label{some bubble}
\limsup_{n\to \infty} \bigl|\bigl\langle \nabla u_n(t_n'), e^{is_n\Delta} \nabla v_n^1(\tau_n) \bigr\rangle\bigr|
    \geq \lim_{n\to \infty}\|\nabla v_n^1(\tau_n)\|_2^2 -c(\eta),
\end{align}
provided $\eta$ is chosen sufficiently small.  Here $c(\eta)$ denotes a small power of $\eta$ (depending upon the dimension
$d$). Using the Cauchy--Schwarz inequality together with \eqref{some bubble} and \eqref{bubble}, we obtain
\begin{align*}
\limsup_{n\to \infty} \int_{|x- x_n^1|\leq R_n} |\nabla u(t_n',x)|^2 \, dx
        \geq \frac{ \bigl[\lim_{n\to \infty}\|\nabla v_n^1(\tau_n)\|_2^2-c(\eta)\bigr]^2}{\lim_{n\to \infty}\|\nabla v_n^1(\tau_n)\|_2^2}.
\end{align*}
Invoking \eqref{big ke} and recalling that $\eps$ and $\eta$ may be taken arbitrarily small completes the proof of
Theorem~\ref{T:conc}.
\end{proof}

\begin{remark}
One may wonder whether it is possible to show that kinetic energy concentrates along every sequence
approaching the blowup time.
In general, we are not able to prove results of this nature; in particular, we are unable to verify such a
claim in \cite[Corollary~5.18]{kenig-merle}.  The obstruction is as follows: we cannot preclude the possibility
that the solution rapidly alternates between being spread out and being concentrated as one approaches the blowup time.
One exception is when $\sup_{t\in I} \|\nabla u(t)\|_{L_x^2}^2 < 2 \|\nabla W\|_{L_x^2}^2$.  In this case, one may apply
Keraani's argument from the proof of \cite[Theorem~1.6]{keraani-l2}.
\end{remark}

%
%
%
%

\appendix

\section{Properties of $W$}\label{A: W}

In this paper,
$$
W(x):=\frac 1{(1+\frac{|x|^2}{d(d-2)})^{\dtt}},
$$
which solves the nonlinear elliptic equation
\begin{align}\label{W eq}
\Delta W + |W|^{\frac{4}{d-2}} W = 0.
\end{align}
Thus $W(t,x)=W(x)$ is a stationary solution to \eqref{nls}.

From the work of Aubin~\cite{aubin} and Talenti~\cite{talenti} (see also \cite[\S 8.3]{LiebLoss}) we know that $W$ is a
maximizer in the sharp Sobolev embedding inequality:

\begin{theorem}[Sharp Sobolev Embedding]\label{T:sharp embed}
\begin{equation}\label{embedding}
\|f\|_{\tdt}\le C_d\|\nabla f\|_2,
\end{equation}
with equality if and only if
\begin{equation}\label{embedding equality}
f(x) = c \lambda^{-\frac{d-2}2} W\bigl(\tfrac{x-x_0}\lambda\bigr)
\end{equation}
for some $c\in\C$, $x_0\in\R^d$, and $\lambda >0$.
\end{theorem}

It is not hard to see that $W\in \dot H_x^1(\R^d)$, but $W\notin L_x^2(\R^d)$ for $d=3,4$. Moreover, multiplying \eqref{W
eq} by $W$ and integrating by parts, we obtain
\begin{equation*}
\|\nabla W\|_2^2=\|W\|_{\tdt}^{\tdt}.
\end{equation*}
Combining this with Theorem~\ref{T:sharp embed}, we get
\begin{align}
\|\nabla W\|_2^2 &=\|W\|_{\tdt}^{\tdt}=C_d^{-d} \label{W ke}\\
E(W) &=\tfrac 12 \|\nabla W\|_2^2-\tfrac{d-2}{2d}\|W\|_{\tdt}^{\tdt}=d^{-1} C_d^{-d}. \label{W e}
\end{align}

The next three lemmas reproduce observations of Kenig and Merle \cite{kenig-merle}; they are reminiscent of Weinstein's
work in the mass-critical setting \cite{weinstein}.  We include details for the sake of completeness.

\begin{lemma}[Coercivity I, \cite{kenig-merle}]\label{coercive}
Assume $E(u_0)\le (1-\delta_0)E(W)$ for some $\delta_0>0$.  Then there exist two positive constants $\delta_1$ and $\delta_2$
(depending on $\delta_0$) such that

a) If $\|\nabla u_0\|_2\le \|\nabla W\|_2$, then
\begin{align*}
\|\nabla u_0\|_2^2\le (1-\delta_1)\|\nabla W\|_2^2.
\end{align*}

b) If $\|\nabla u_0\|_2\ge \|\nabla W\|_2$, then
\begin{align*}
&1)\ \|\nabla u_0\|_2^2\ge (1+\delta_2)\|\nabla W\|_2^2 \\
&2)\ \int_{\R^d} \bigl(|\nabla u_0(x)|^2-|u_0(x)|^{\tdt}\bigr)\,dx\le -\tfrac {2\delta_2}{(d-2)C_d^d}.
\end{align*}
\end{lemma}

\begin{proof}
We first consider the case when $\|\nabla u_0\|_2\le \|\nabla W\|_2$.  We define
\begin{equation*}
f(y):=\tfrac 12 y-\tfrac {d-2}{2d}C_d^{\tdt}y^{\frac d{d-2}}.
\end{equation*}
By \eqref{embedding},
\begin{equation}\label{a1}
f\bigl(\|\nabla u_0\|_2^2\bigr)\le E(u_0),
\end{equation}
while by \eqref{W ke} and \eqref{W e},
\begin{align*}
f\bigl(\|\nabla W\|_2^2\bigr) = E(W) = d^{-1} C_d^{-d}.
\end{align*}
Computing the first and second derivatives of $f$, we learn that $f$ is concave and attains its maximum value (that is,
$d^{-1}C_d^{-d}$) only at $C_d^{-d}=\|\nabla W\|_2^2$.  Moreover, $f$ is strictly increasing on $[0,C_d^{-d}]$ and strictly
decreasing on $[C_d^{-d}, \infty)$.  In particular, the inverse function (which we denote by $f^{-1}$) is well-defined and
strictly increasing on $[0, d^{-1}C_d^{-d}]$; hence, there exists $\delta_1> 0$ such that
$$
f^{-1}\bigl((1-\delta_0)E(W)\bigr)=(1-\delta_1)\|\nabla W\|_2^2.
$$
By hypothesis and \eqref{a1},
$$
\|\nabla u_0\|_2^2\le f^{-1}(E(u_0))\le f^{-1}((1-\delta_0)E(W))=(1-\delta_1)\|\nabla W\|_2^2.
$$
This settles the first claim.

We now consider the case $\|\nabla u_0\|_2\ge \|\nabla W\|_2$.  From the analysis above, we know that $f$ is strictly decreasing
on $[C_d^{-d}, \infty)$; thus, the inverse function is well-defined and strictly decreasing on $(-\infty, d^{-1}C_d^{-d}]$.
Therefore there exists $\delta_2>0$ such that
$$
f^{-1}\bigl((1-\delta_0)E(W)\bigr)=(1+\delta_2)\|\nabla W\|_2^2.
$$
Invoking the hypothesis and \eqref{a1},
$$
\|\nabla u_0\|_2^2\ge f^{-1} (E(u_0)) \geq f^{-1}[(1-\delta_0)E(W)]=(1+\delta_2)\|\nabla W\|_2^2,
$$
which settles the first claim in this case.

To prove the last claim, we use the hypothesis and the first claim to estimate
\begin{align*}
\|\nabla u_0\|_2^2-\|u_0\|_{\tdt}^{\tdt}
&=\tfrac{2d}{d-2} E(u_0)-\tfrac 2{d-2}\|\nabla u_0\|_2^2\\
&\leq \tfrac{2d}{d-2} (1-\delta_0) E(W) - \tfrac 2{d-2}(1+\delta_2) \|\nabla W\|_2^2\\
&\leq -\tfrac{2\delta_2}{d-2}C_d^{-d}.
\end{align*}
This finishes the proof of the lemma.
\end{proof}

Combining Lemma~\ref{coercive} with a continuity argument and the conservation of energy, one easily deduces

\begin{corollary}[Coercivity II, \cite{kenig-merle}]\label{trap}
Let $u:\ird \to \C$ be a solution to \eqref{nls} with initial data $u(t_0)=u_0\in \dot H^1_x(\R^d)$ for some $t_0\in I$. Assume
$E(u_0)\le (1-\delta_0) E(W)$ for some $\delta_0>0$.  Then there exist two positive constants $\delta_1$ and $\delta_2$
(depending on $\delta_0$) such that

a) If $\|\nabla u_0\|_2\le \|\nabla W\|_2$, then for all $t\in I$
\begin{align*}
\|\nabla u(t)\|_2^2\le (1-\delta_1)\|\nabla W\|_2^2.
\end{align*}

b) If $\|\nabla u_0\|_2\ge \|\nabla W\|_2$, then for all $t\in I$
\begin{align*}
&1)\ \|\nabla u(t)\|_2^2\ge (1+\delta_2)\|\nabla W\|_2^2\\
&2)\ \int_{\R^d}\bigl(|\nabla u(t,x)|^2-|u(t,x)|^{\tdt}\bigr)\,dx \le -\tfrac{2\delta_2}{(d-2)C_d^d}.
\end{align*}
\end{corollary}

Finally, we show that if $u$ is a solution to \eqref{nls} with kinetic energy (at all times) less than that of the ground state
$W$, then the kinetic energy of $u$ does not vary greatly over the interval of existence.

\begin{lemma}[Coercivity III, \cite{kenig-merle}]\label{equal}
Let $u:\ird\to \C$ be a solution to \eqref{nls}.  Assume also that $\sup_{t\in I}\|\nabla u(t)\|_2\leq (1-\delta) \|\nabla
W\|_2$ for some $\delta>0$.  Then for all $t\in I$
\begin{align}\label{ke sim}
E(u(t))\sim \|\nabla u(t)\|_2^2 \sim \|\nabla u_0\|_2^2
\end{align}
and
$$
\int_{\R^d}\bigl(|\nabla u(t,x)|^2-|u(t,x)|^{\tdt}\bigr)\,dx \gtrsim \|\nabla u_0\|_2^2,
$$
where the implicit constants depend only upon $\delta$ (and the dimension $d$).
\end{lemma}

\begin{proof}
Using \eqref{embedding}, \eqref{W ke}, and the hypothesis,
\begin{align*}
\frac 12 \|\nabla u(t)\|_2^2 \geq E(u(t)) &=\frac12 \|\nabla u(t)\|_2^2 \Bigl[ 1- \frac{d-2}d \Bigl(\frac {\|\nabla
u(t)\|_2}{\|\nabla W\|_2} \Bigr)^{\frac4{d-2}} \Bigr] \gtrsim \|\nabla u(t)\|_2^2.
\end{align*}
The second relation in \eqref{ke sim} follows from the conservation of energy.

By \eqref{embedding}, \eqref{W ke} and the hypothesis,
\begin{align*}
\int_{\R^d}\bigl(|\nabla u(t,x)|^2-|u(t,x)|^{\tdt}\bigr)\,dx \geq \Bigl[1-\Bigl(\frac{\|\nabla u(t)\|_2}{\|\nabla
W\|_2}\Bigr)^{\frac4{d-2}} \Bigr]\|\nabla u(t)\|_2^2 \gtrsim \|\nabla u(t)\|_2^2.
\end{align*}
The second claim now follows from \eqref{ke sim}.
\end{proof}

%
%
%
%


\begin{thebibliography}{10}
\newcommand{\msn}[1]{\href{http://www.ams.org/mathscinet-getitem?mr=#1}{\sc MR#1}}

\bibitem{aubin} T. Aubin, \emph{Probl\'emes isop\'erim\'etriques et espaces de Sobolev}, J. Diff. Geom. \textbf{11} (1976), 573--598.
\msn{0448404}


\bibitem{BegoutVargas} P. Begout and A. Vargas,
\emph{Mass concentration phenomena for the $L^2$-critical nonlinear Schr\"odinger equation},
Trans. Amer. Math. Soc. \textbf{359} (2007), 5257--5282. \msn{2327030}

\bibitem{bourg.2d}
J. Bourgain, \emph{Refinements of Strichartz inequality and applications to 2d-NLS with critical nonlinearity,} Int. Math. Res.
Not. (1998), 253--283. \msn{1616917}

\bibitem{borg:scatter}
J. Bourgain, \emph{Global wellposedness of defocusing critical nonlinear Schr\"odinger equation in the radial case},
J. Amer. Math. Soc. \textbf{12} (1999), 145--171. \msn{1626257}

\bibitem{borg:book}
J. Bourgain, \emph{Global solutions of nonlinear Schr\"odinger equations.}
American Mathematical Society Colloquium Publications, \textbf{46}. American Mathematical Society, Providence, RI, 1999.
\msn{1691575}


\bibitem{cwI}
T. Cazenave and F.B. Weissler, \emph{The Cauchy problem for the critical nonlinear Schr\"odinger equation in $H\sp s$},
Nonlinear Anal. \textbf{14} (1990), 807--836.
\msn{1055532}

\bibitem{cazenave:book}
T. Cazenave, \textit{Semilinear Schr\"odinger equations.}
Courant Lecture Notes in Mathematics, \textbf{10}. American Mathematical Society, 2003.
\msn{2002047}

\bibitem{ChW:fractional chain rule}
M.~Christ and M.~Weinstein, \textit{Dispersion of small amplitude
solutions of the generalized Korteweg-de Vries equation,}
J. Funct. Anal. \textbf{100} (1991), 87--109.
\msn{1124294}

\bibitem{ckstt:gwp}
J. Colliander, M. Keel, G. Staffilani, H. Takaoka, and T. Tao, \emph{Global well-posedness and scattering for the
energy-critical nonlinear Schr\"odinger equation in $\R^3$}, to appear in Annals of Math.

\bibitem{DHR}
T. Duyckaerts, J. Holmer, and S. Roudenko, \emph{Scattering for the non-radial 3D cubic nonlinear Schr\"odinger equation,}
prerint \texttt{arXiv:0710.3630}.

\bibitem{DuckyMerle}
T. Duyckaerts and F. Merle, \emph{Dynamic of threshold solutions for energy-critical NLS}, \texttt{arXiv:0710.5915}.

\bibitem{gv:strichartz}
J. Ginibre and G. Velo, \emph{Smoothing properties and retarded estimates for some dispersive evolution equations,}
Comm. Math. Phys. \textbf{144} (1992), 163--188.
\msn{1151250}

\bibitem{glassey}
R. T. Glassey, \emph{On the blowing up of solutions to the Cauchy problem for nonlinear Schr\"odinger equations},
J. Math. Phys. \textbf{18} (1977), 1794--1797.
\msn{0460850}

\bibitem{grillakis}
G. Grillakis, \emph{On nonlinear Schr\"odinger equations}, Comm. PDE \textbf{25} (2000), 1827--1844.
\msn{1778782}

\bibitem{HolmerRoudenko}
J. Holmer and S. Roudenko, \emph{A sharp condition for scattering of the radial 3D cubic nonlinear Schr\"odinger equation},
preprint \texttt{arXiv:0703235}.

\bibitem{tao:keel}
M. Keel and T. Tao, \emph{Endpoint Strichartz Estimates}, Amer. J. Math. \textbf{120} (1998), 955--980.
\msn{1646048}

\bibitem{Evian}
C. E. Kenig, \emph{Global well-posedness and scattering for the energy critical focusing nonlinear Schr\"odinger and wave
equations.} Lecture notes for a mini-course given at ``Analyse des \'equations aux d\'eriv\'ees partielles,''
Evian-les-Bains, July 2007.

\bibitem{kenig-merle}
C. E. Kenig and F. Merle, \emph{Global well-posedness, scattering and blow up for the energy-critical, focusing, nonlinear
Schr\"odinger equation in the radial case}, Invent. Math. \textbf{166} (2006), 645--675.
\msn{2257393}

\bibitem{kenig-merle:wave}
C. E. Kenig and F. Merle, \emph{Global well-posedness, scattering and blow-up for the energy critical focusing non-linear
wave equation}, \texttt{arXiv:0610801}.

\bibitem{kenig-merle:1/2}
C. E. Kenig and F. Merle, \emph{Scattering for $H^{1/2}$ bounded solutions to the cubic, defocusing NLS in 3 dimensions},
\texttt{arXiv:0712.1834}.

\bibitem{keraani-h1}
S. Keraani, \emph{On the defect of compactness for the Strichartz estimates for the Schr\"odinger equations},
J. Diff. Eq. \textbf{175} (2001), 353--392.
\msn{1855973}

\bibitem{keraani-l2}
S. Keraani, \emph{On the blow-up phenomenon of the critical nonlinear Schr\"odinger equation},
J. Funct. Anal. \textbf{235} (2006), 171--192.
\msn{2216444}

\bibitem{KTV}
R. Killip, T. Tao, and M. Visan, \emph{The cubic nonlinear Schr\"odinger equation in two dimensions with radial data,} preprint
\texttt{arXiv:0707.3188}.

\bibitem{KVZ-Q}
R. Killip, M. Visan, and X. Zhang, \emph{Energy-critical NLS with quadratic potentials,} to appear Comm. PDE.

\bibitem{KVZ}
R. Killip, M. Visan, and X. Zhang, \emph{The mass-critical nonlinear Schr\"odinger equation with radial data in dimensions three
and higher,} preprint \texttt{arXiv:0708.0849}.

\bibitem{RadialHD}
R. Killip, M. Visan, and X. Zhang, \emph{The focusing energy-critical nonlinear Schr\"odinger equation with radial data},
unpublished manuscript Sept. 2007.

\bibitem{ll}
L. D. Landau and E. M. Lifshitz, \emph{Course of theoretical physics. Vol. 1. Mechanics.} Third edition.
Pergamon Press, Oxford-New York-Toronto, Ont., 1976.

\bibitem{LiebLoss}
E. H. Lieb and M. Loss, \emph{Analysis. Second edition.} Graduate Studies in Mathematics, \textbf{14}.
American Mathematical Society, Providence, RI, 2001.
\msn{1817225}


\bibitem{merle-vega}
F. Merle and L. Vega, \emph{Compactness at blow-up time for $L^2$ solutions of the critical nonlinear Schr\"odinger equation in 2$D$},
Internat. Math. Res. Notices \textbf{8} (1998), 399--425.
\msn{1628235}

\bibitem{nakanishi}
K. Nakanishi, \emph{Scattering theory for nonlinear Klein-Gordon equation with Sobolev critical power},
Internat. Math. Res. Notices \textbf{1} (1999), 31--60.
\msn{1666973}

\bibitem{OgawaTsutsumi}
T. Ogawa and Y. Tsutsumi, \emph{Blow-up of $H^1$ solution for the nonlinear Schr\"odinger equation},
J. Diff. Eq. \textbf{92} (1991), 317--330.
\msn{1120908}

\bibitem{RV}
E.~Ryckman and M.~Visan, \emph{Global well-posedness and scattering for the defocusing energy-critical nonlinear Schr\"odinger
equation in $\R^{1+4}$,} Amer. J. Math. \textbf{129} (2007), 1--60.
\msn{2288737}


\bibitem{stein:large}
E.~M. Stein, \emph{Harmonic analysis: real-variable methods, orthogonality, and oscillatory integrals},
Princeton Mathematical Series, 43. Princeton University Press, Princeton, NJ, 1993.
\msn{1232192}

\bibitem{strichartz}
R. S. Strichartz, \emph{Restriction of Fourier transform to quadratic surfaces and deay of solutions of wave equations},
Duke Math. J. \textbf{44} (1977), 705--774.
\msn{0512086}

\bibitem{talenti}
G. Talenti, \emph{Best constant in Sobolev inequality}, Ann. Mat. Pura. Appl. \textbf{110} (1976), 353--372.
\msn{0463908}

\bibitem{tao: gwp radial}
T. Tao, \emph{Global well-posedness and scattering for the higher-dimensional energy-critical non-linear Schr\"odinger equation
for radial data}, New York J. of Math. \textbf{11} (2005), 57--80.
\msn{2154347}

\bibitem{tao:book}
T. Tao, \emph{Nonlinear dispersive equations. Local and global analysis.} CBMS Regional Conference Series in Mathematics, 106.
American Mathematical Society, Providence, RI, 2006.
\msn{2233925}

\bibitem{tao:attractor}
T. Tao, \emph{A (concentration-)compact attractor for high-dimensional non-linear Schr\"odinger equations},
Dyn. Partial Differ. Equ. \textbf{4} (2007), 1--53.
\msn{2304091}

\bibitem{tao:pseudo}
T. Tao, \emph{A pseudoconformal compactification of the nonlinear Schr\"odinger equation and applications},
\texttt{arXiv:math/0606254}.

\bibitem{TV}
T. Tao and M.~Visan, \emph{Stability of energy-critical nonlinear Schr\"odinger equations in high dimensions,}
Electron. J. Diff. Eqns. \textbf{118} (2005), 1--28.
\msn{2174550}

\bibitem{tvz:cc}
T. Tao, M. Visan, and X. Zhang, \emph{Minimal-mass blowup solutions of the mass-critical NLS}, to appear in Forum Math.

\bibitem{Taylor:book}
M. E. Taylor, \textit{Tools for PDE.} Mathematical Surveys and Monographs, \textbf{81}.
American Mathematical Society, Providence, RI, 2000.
\msn{1766415}

\bibitem{thesis:art} M. Visan, \emph{The defocusing energy-critical nonlinear Schr\"odinger equation
in higher dimensions,} Duke Math. J. \textbf{138} (2007), 281--374.
\msn{2318286}

\bibitem{Monica:thesis} M. Visan, \emph{The defocusing energy-critical nonlinear Schr\"odinger equation
in dimensions five and higher,} Ph.D. Thesis, UCLA, 2006.

\bibitem{weinstein}
M. Weinstein, \emph{Nonlinear Schr\"odinger equations and sharp interpolation estimates,}
Comm. Math. Phys. \textbf{87} (1983), 567--576.
\msn{0691044}

\bibitem{zack} V.E. Zakharov, E.A. Kuznetsov, \emph{Quasi-classical theory for three-dimensional wave collapse},
Sov. Phys. JETP \textbf{64} (1986), 773--380.

\end{thebibliography}
\end{document}